\documentclass[10pt]{articleHJ}

\usepackage[centertags]{amsmath}
\usepackage{amsfonts,amssymb,amsthm}
\usepackage{caption,graphicx}
\usepackage{subcaption}
\usepackage{mathrsfs}
\usepackage{tikz}
\usepackage{url}
\usepackage[text={6.0in,8.6in},centering,letterpaper]{geometry}
\usepackage[numbers,comma,square,sort&compress]{natbib}

\usepackage{mathtools}

\newlength{\defbaselineskip}
\newcommand{\setlinespacing}[1]%
           {\setlength{\baselineskip}{#1 \defbaselineskip}}

\newcommand{\singlespacing}{\setlength{\baselineskip}{\defbaselineskip}}

\newcommand{\e}{\ensuremath{\varepsilon}}

\newcommand{\sref}[1]{(\ref{#1})}                       

\newtheorem{remark}{Remark}
\newtheorem*{defn*}{Definition}

 {\begin{trivlist}\item[]\textbf{Acknowledgments }}{\end{trivlist}}

\begin{document}

\begin{frontmatter}
\title{Waves in a Stochastic Cell Motility Model}
\journal{....}

\author[LD1]{C. H. S. Hamster\corauthref{coraut}},
\corauth[coraut]{Corresponding author. }
\author[LD2]{P. van Heijster},
\address[LD1]{
  Biometris - Wageningen University and Research \\
 Wageningen; The Netherlands \\
 Email:  {\normalfont{\texttt{christian.hamster@wur.nl}}}
}
\address[LD2]{
  Biometris - Wageningen University and Research \\
 Wageningen; The Netherlands  \\ Email:  {\normalfont{\texttt{peter.vanheijster@wur.nl}}}
}

\date{\today}

\begin{abstract}
\singlespacing
In Bhattacharya et al. (Science Advances, 2020), a set of chemical reactions involved in the dynamics of actin waves in cells was studied. Both at the microscopic level, where the individual chemical reactions are directly modelled using Gillespie-type algorithms, and on a macroscopic level where a deterministic reaction-diffusion equation arises as the large-scale limit of the underlying chemical reactions.
In this work, we derive, and subsequently study, the related mesoscopic stochastic reaction-diffusion system, or Chemical Langevin Equation, that arises from the same set of chemical reactions. We explain how the stochastic patterns that arise from this equation can be used to understand the experimentally observed dynamics from Bhattacharya et al. In particular, we argue that the mesoscopic stochastic model better captures the microscopic behaviour than the deterministic reaction-diffusion equation, while being more amenable for mathematical analysis and numerical simulations than the microscopic model.
\end{abstract}


\end{frontmatter}

\numberwithin{equation}{section}
\numberwithin{figure}{section}
\renewcommand{\theequation}{\thesection.\arabic{equation}}

\section{Introduction}
\begin{figure}[t]
    \centering
    \includegraphics[width=1\columnwidth]{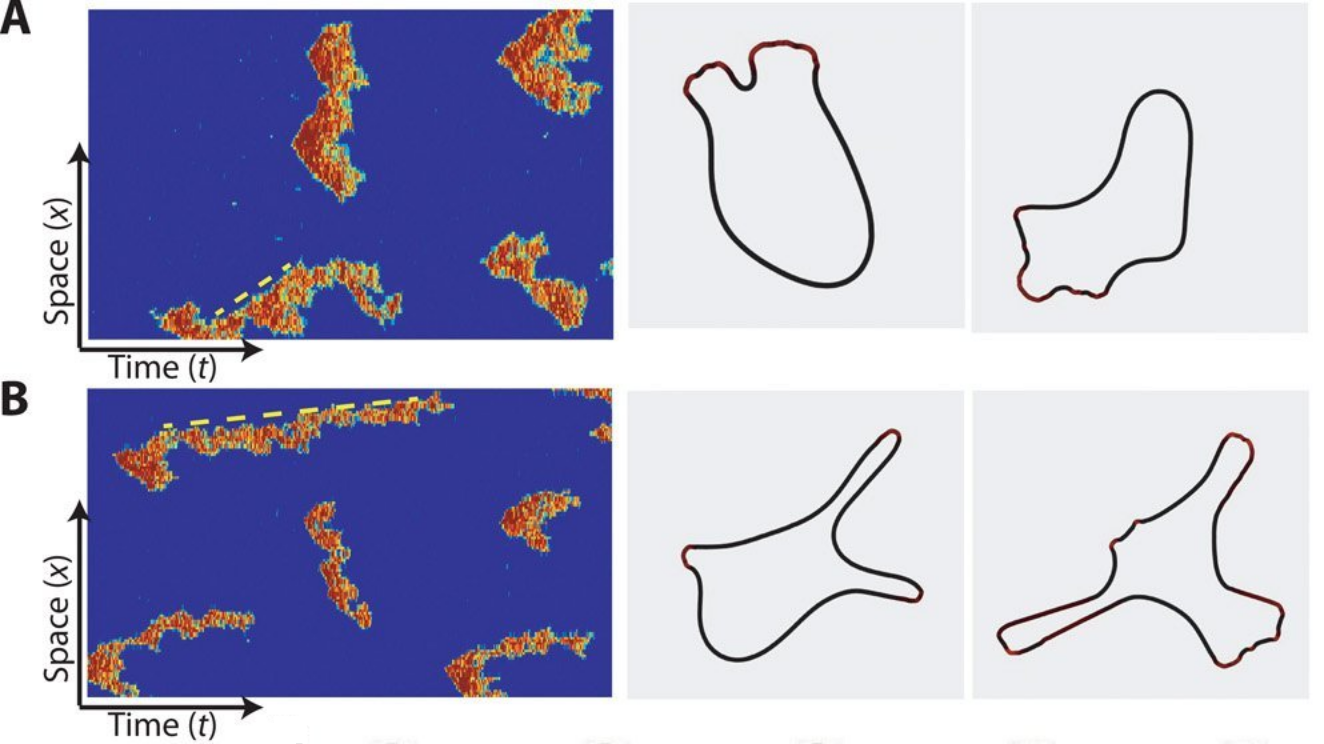}
    \caption{Stochastic simulations of the microscopic Gillespie-type model from~\cite{bhattacharya2020traveling}. The figures on the left show stochastic simulations of the Ras activity for parameter values applicable to (A) wild-type cells and to (B) genetically modified cells, where the phosphatase PTEN has been switched off. The figures on the right show typical cell shapes corresponding to the dynamics in the left figures. This shows that mutations in the gene that codes for PTEN lead to elongated protrusions typically associated with cancer. The dotted yellow line is an indicator of the wave speed, i.e. the actin waves in (B) are slower and live longer than in (A). Reproduced from~\cite{bhattacharya2020traveling} under creative commons license 4.0. }
    \label{fig:FromScienceSim}
\end{figure}

In order to move around, an amoeboid cell can change its shape by polymerising actin to curve the cell membrane. The actin polymerisation is controlled by signalling molecules and experiments in \textit{Dictyostelium discoideum} have shown that activation of these signalling molecules happens at localised patches that can move along the membrane like a wave~\cite{inagaki2017actin,bhattacharya2020traveling}. In wild-type (WT) cells, these waves move fast and die out, creating familiar-shaped pseudopods, while in cancerous cells these waves stick to a point, creating elongated protrusions~\cite{bhattacharya2020traveling}, see Figure~\ref{fig:FromScienceSim}. In absence of a signal, the formation of pseudopods happens at random places on the cell membrane, resulting in random motion. In contrast, when a cell senses a chemical signal, it can concentrate the random protrusions at the side of the cell where the signal comes from, leading to movement in the direction of the signal~\cite{deng2022introduction}. As cells are small, the difference in signal strength between the front and the back of the cell (the gradient) is small as well. Furthermore, the cell can only use discrete points at the membrane where the receptors are to estimate the direction of the signal~\cite{deng2022introduction}. Therefore, one of the main questions is ``How can a cell use a small gradient in the signal to concentrate the actin activity in the front?". This question has been studied intensively, but no complete description of all the microscopic chemical processes involved has been given yet, see~\cite{devreotes2017excitable} for a review. 

In~\cite{bhattacharya2020traveling}, the choice is made to describe the highly complex actin dynamics with a conceptual activator $u$ and inhibitor $v$ that diffuse and react with each other as summarised in Table~\ref{tab:my_label}. The species $u$ and $v$ are an abstraction of the dozens of components that regulate the actual cell movement, but the activator $u$ can be thought of as Ras activity~\cite{bhattacharya2020traveling}, which plays an important role in cell growth and differentiation~\cite{lodish2008molecular}. In particular, $u$ is being activated by Reaction $\#3$ and Reaction $\#4$, while being inhibited by Reaction $\#1$ and Reaction $\#2$, with propensities as indicated in the table. In addition, $v$ is inhibited by Reaction $\#5$, while Reaction $\#6$ activates the inhibitor. 

The information on the chemical reactions, in combination with the diffusion of both species, is generally used in one of two ways. First, there is a Gillespie-type algorithm~\cite{gillespie1976general, gillespie1977exact} which can be used to simulate the involved chemical reactions on a microscopic level. For these simulations, $(u_k(t_n),v_k(t_n))$ (the solution at time $t_n$ at grid cell $k$) is treated as the number of molecules of type $u$ and $v$ at time $t_n$ in a grid cell with finite size. For all these individual molecules the probabilities of diffusing to other grid cells or taking part in a chemical reaction are prescribed as by Table~\ref{tab:my_label}. To be precise, Reaction $\#1$ implies that the time to the next reaction that degrades a $u$ molecule in grid cell $k$ is exponentially distributed with rate parameter  $(a_1u_k(t_n))^{-1}$. See the panels on the left of Figure~\ref{fig:FromScienceSim} for examples of these simulations. This Gillespie-type algorithm approach takes the stochastic nature of a single cell into account. However, it is computationally very expensive and difficult to analyse mathematically. Hence, it is hard to use this type of modelling approach to make valuable predictions.

A second way to use the reactions in Table~\ref{tab:my_label} is to derive an average large-scale limit macroscopic equation. Hence, we assume that $u$ and $v$ are densities on a continuous domain, described by a reaction-rate equation with diffusion, also known as a Reaction-Diffusion Equation (RDE). In particular, the RDE\footnote{Note that~\sref{eq:int:MainPDE} can also be obtained from a quasi-steady-state approximation from a more complex three-component model introduced in~\cite{biswas2021enhanced}, commonly referred to as a Signal Transduction Excitation Network.} related to the chemical reactions in Table~\ref{tab:my_label} is given by  
\begin{align}
\begin{split}
\label{eq:int:MainPDE}
    \partial_t u&=D_u \partial_{xx} u -a_1 u-a_2uv+\frac{a_3u^2}{a_4+u^2}+a_5\,,\\
    \partial_t v&=D_v \partial_{xx} v +\e(-c_1v+c_2u),
    \end{split}
\end{align}
which is a specific version of the general RDE we will encounter in \S\ref{sec:der}.
This model is a variation on the classic FitzHugh-Nagumo model for neuron spiking~\cite{fitzhugh1961impulses, nagumo1962active}. Protrusions are formed at places with high activator $u$ and $u$ is inhibited by the terms $-a_1 u$ and $-a_2 u v$, see Reaction $\#1$ and Reaction $\#2$ in Table~\ref{tab:my_label}. This implies that an increase in $u$ or $v$ leads to a decrease in $u$, unless the increase is high enough such that activation from Reaction $\#3$, modelled by a nonlinear Hill function $a_3u^2/(a_4+u^2)$, takes over and negates the inhibiting effects. Effectively, this means that a small increase in $u$ can lead to a much larger increase in $u$, that is, the system is locally activated. Once $u$ is large and the Hill function levels off at a fixed value $a_3$, the amount of inhibitor $v$ increases via the term $\e c_2u$ (related to Reaction $\#6$), leading to a fast decay in $u$ by the $-a_2uv$ term (related to Reaction $\#2$). The inhibitor $v$ then decays via Reaction $\#5$ to the rest state and activation can happen again. In addition, both species diffuse with diffusion coefficient $D_u$, respectively $D_v$, where it is assumed that $D_u< D_v$.
It is important to realise that, in both approaches, the modelled actin waves happen on the surface of the cell, and, as in~\cite{bhattacharya2020traveling}, we only study a slice of this surface. Therefore, the spatial domain must be thought of as an (approximate) circle.

\begin{table}[]
    \centering
    \begin{tabular}{c c c c c }
    \hline \hline
            No. & Reaction &Propensity &$u$ & $v$   \\ \hline \\[-3mm]
         1& $u\to\emptyset$ & $a_1u$&$-1$&$0$ \\
         2& $u\to\emptyset$ & $a_2uv$&$-1$&$0$\\
         3& $\emptyset\to u$& $a_3u^2/(a_4+u^2)$ &
         $1$ &$0$\\
         4& $\emptyset\to u$&$a_5$ &$1$&$0$\\
         5& $v\to\emptyset$ & $\e c_1v$&$0$&$-1$\\
         6&$\emptyset\to v$&$\e c_2u$&$0$&$1$\\
    \hline \hline
    \end{tabular}
    \caption{The chemical reactions that determine the actin wave dynamics from~\cite{bhattacharya2020traveling}.}
    \label{tab:my_label}
\end{table}

For deterministic RDEs like \sref{eq:int:MainPDE}, a plethora of analytical tools are available (see, for instance, Appendix~\ref{sec:analysis}) and numerical simulations are relatively straightforward. However, being a deterministic equation, this RDE does not show the same stochastic dynamics as the Gillespie simulations and experiments. A crucial difference between the macroscopic RDE model~\sref{eq:int:MainPDE} and the Gillespie simulations revolves around the duration of the patterns. In the RDE, an established pattern, e.g. a standing or travelling wave, will, if uninterrupted, remain there for a very long time, while these patterns are destroyed quickly both in stochastic simulations and experiments. Furthermore, when the rest state of the RDE~\sref{eq:int:MainPDE} is stable, activation cannot come from the RDE itself, but it needs an external signal large enough to activate the nonlinear term $a_3 u^2/(a_4 + u^2)$ related to Reaction $\#3$. We generally refer to the activation of these patterns as activation events.

\begin{figure}[t]
\begin{subfigure}{.49\textwidth}
  \centering
 		\def\svgwidth{\columnwidth}
    		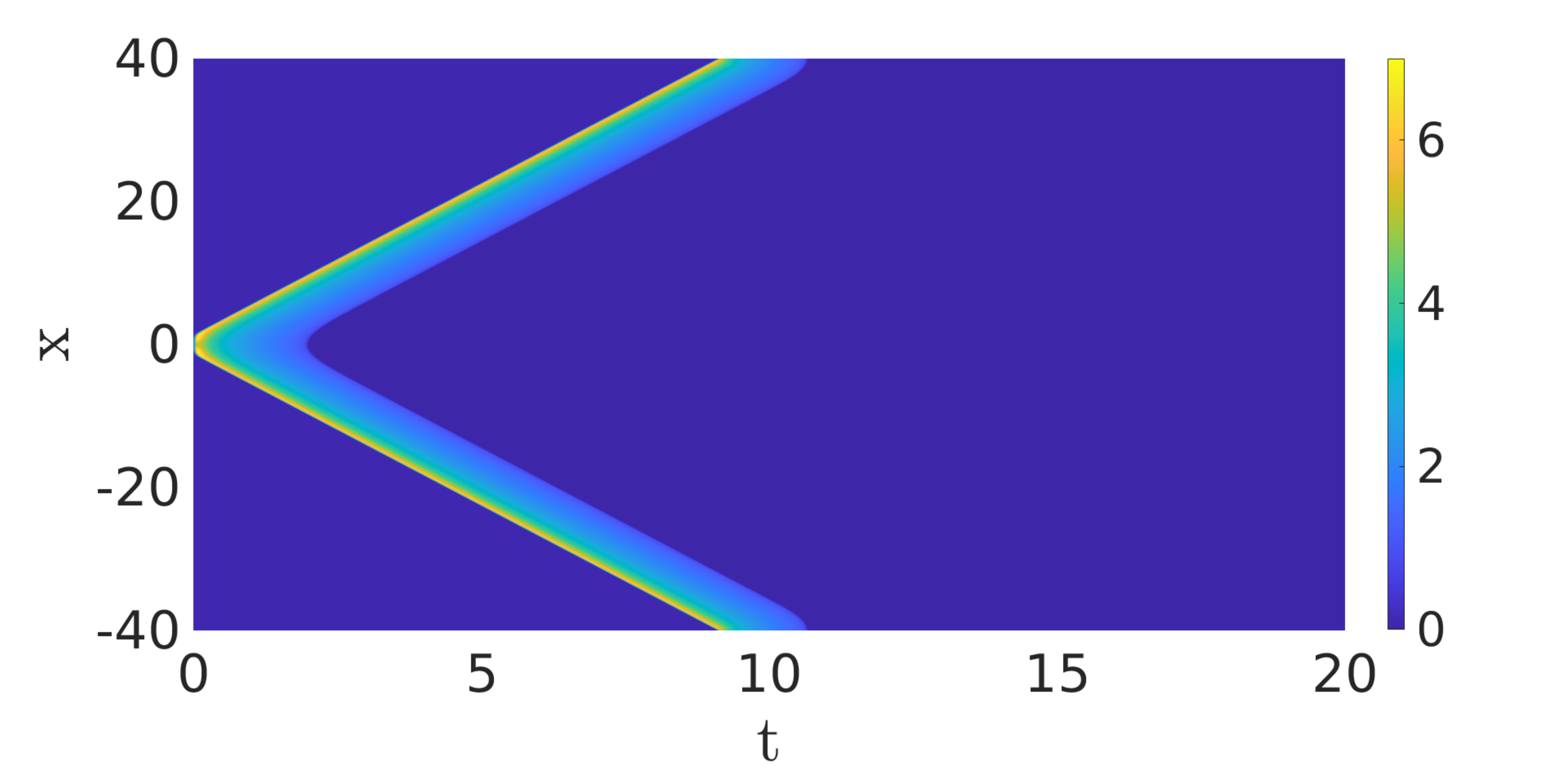
  \caption{}
    \label{fig:Int:DPulseU}
\end{subfigure}
\begin{subfigure}{.49\textwidth}
  \centering
 		\def\svgwidth{\columnwidth}
    		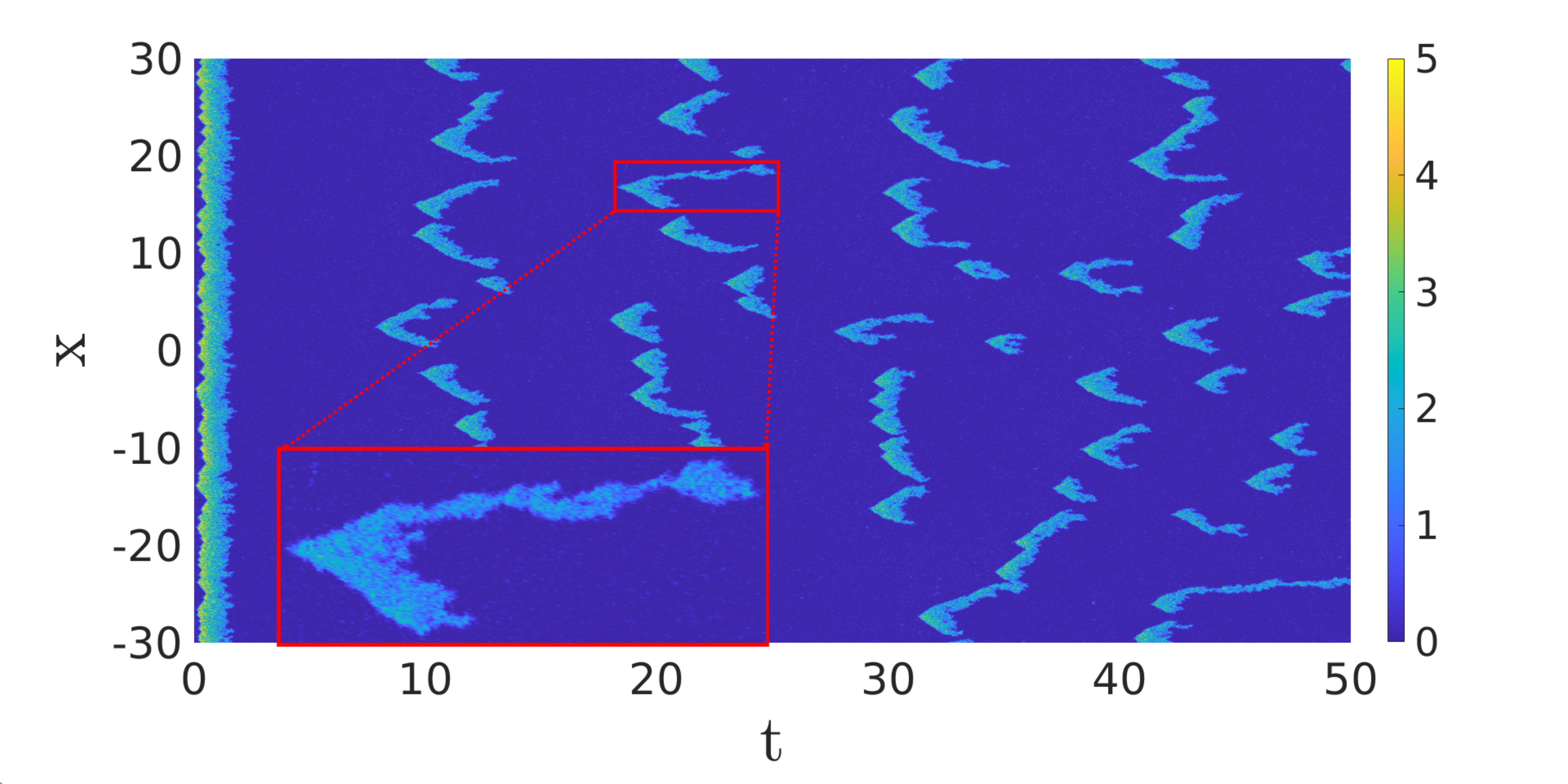
  \caption{}
    \label{fig:Int:DpulseV}
\end{subfigure}
\caption{Comparison of the deterministic model~\sref{eq:int:MainPDE} and its stochastic counterpart~\sref{eq:int:FullSPDE}. In Figure~(a) we show a simulation of~\sref{eq:int:MainPDE}, which is excited at $t=0$, resulting in two counterpropagating travelling waves. In the stochastic simulation in Figure~(b), the influence of the initial excitation quickly disappears and new pulses appear constantly. The same parameters are used as in the simulations shown in the second row of Figure~\ref{fig:FromScienceSim}. Observe the similarities in the shape of the pattern. In Figure~(a), the waves travel around the cell where they cancel each other, while in Figure~(b) the waves cancel each other at a much shorter scale. See \S\ref{sec:WvsP} for more details.}
\label{fig:Int}
\end{figure}

It is important to realise that the dynamics of the different chemical processes in the cell are inherently stochastic and at the size of a single cell chemical reactions are not well approximated by large-scale approximations, as Figures~\ref{fig:FromScienceSim} and \ref{fig:Int} show. In other words, treating the relevant enzymes and receptors like a continuous medium of infinitely many, infinitely small, particles is invalid, and the stochastic nature of reactions between individual molecules becomes important. This so-called \emph{internal noise} can serve as a signal to activate the dynamics if it is large enough at a certain point in space and time. As we noted before, the cell hence executes a random walk in the absence of a signal\footnote{Describing the motion of free cells is a very subtle problem and random motion does not necessarily mean Brownian motion~\cite{li2008persistent,selmeczi2008cell}.}. This implies that an external signal does not necessarily activate the dynamics at a certain point on the membrane, but rather changes the random walk of the cell into a biased random walk in the direction of the signal. Using a more extended model than presented here, it is shown in~\cite{biswas2021enhanced} that coupling an external signal to the stochastic dynamics of the cell indeed can lead to movement in the direction of that signal. 

Instead of studying the complex internal dynamics of the cell, it can be advantageous to perturb the deterministic RDE \sref{eq:int:MainPDE}. For instance, in~\cite{bhattacharya2020traveling}, an external source of noise is applied to the RDE~\sref{eq:int:MainPDE}, turning it into a Stochastic RDE (or Stochastic Partial Differential Equation (SPDE)). While this approach can indeed activate the dynamics and make long-term deterministic waves collapse, it is inherently \emph{ad hoc} and not \emph{a priori} based on any of the involved biologically relevant processes.

In between the macroscopic level of the RDE and the microscopic level of the chemical reactions, one can derive a mesoscopic SPDE, known as a Chemical Langevin Equation (CLE)~\cite{gillespie2000chemical}, that also incorporates the internal noise of the cell. 
In \S\ref{sec:der}, we will show that the SPDE associated with the chemical reactions as described in Table~\ref{tab:my_label} plus diffusion is given by
\begin{align}
\label{eq:int:FullSPDE}
    \begin{split}
        du&=\left(D_u \partial_{xx} u-(a_1+a_2v)u+\frac{a_3u^2}{a_4+u^2}+a_5\right)dt+\sigma\sqrt{(a_1+a_2v)u+\frac{a_3u^2}{a_4+u^2}+a_5}\,dW^1_t\\
        &\qquad\qquad+\sigma\partial_x\sqrt{2D_uu}\,d\tilde W^1_t,\\
        dv&=\left(D_v \partial_{xx}v+\e(-c_1v+c_2u)\right)dt+\sigma\sqrt{\e(-c_1v+c_2u)}\,dW^2_t+\sigma\partial_x\sqrt{2D_vv}\,d\tilde W^2_t.
    \end{split}
\end{align}
Here, $(d W^{1}_t, d W^{2}_t)$ and $(d \tilde{W}^{1}_t,d \tilde{W}^{2}_t)$ are two independent noise vectors with space-time white noise (each component is also independent of the other) and $\sigma$ is a measure for the strength of the noise. Indeed, in the no-noise limit $\sigma \to 0$ the mesoscopic SPDE~\sref{eq:int:FullSPDE} reduces to the macroscopic RDE~\sref{eq:int:MainPDE}. In that sense, $\sigma$ serves as a scale parameter. 

The main advantage of the SPDE description is, on one hand, that the solutions still show the rich dynamics of the Gillespie models, i.e. the activation and destruction of waves, but are computationally significantly less expensive. On the other hand, since the SPDE in the no-noise limit reduces to the deterministic RDE model~\sref{eq:int:MainPDE}, we can use well-developed Partial Differential Equation (PDE) theory to gain insight into the dynamics of the RDE~\sref{eq:int:MainPDE} and use this to study the closely related SPDE, see for instance \cite{hamster2020travelling, kuehnreview}. To give an idea of the differences between the deterministic and stochastic models we plot two simulations in Figure~\ref{fig:Int} that will be discussed later in \S\ref{sec:sim}. It is clear that the simulation of the SPDE paints a much more dynamic picture than the deterministic one, which is more in line with the inherently noisy nature of the cell's chemical processes. Hence, SPDEs are an invaluable tool in unravelling the dynamics of a cell. 

This article is now organised as follows. In $\S$\ref{sec:der} we explain how to derive the SPDE~\sref{eq:int:FullSPDE} from Table~\ref{tab:my_label}. Subsequently, in $\S$\ref{sec:sim} we study both the SPDE~\sref{eq:int:FullSPDE} and the RDE~\sref{eq:int:MainPDE} numerically in different parameter regimes and qualitatively compare the observed dynamics to the Gillespie simulations from~\cite{bhattacharya2020traveling}. In $\S$\ref{sec:dis}, we discuss the results and how they relate to the questions posed in this introduction. 


\section{Derivation of the SPDE}
\label{sec:der}
Our starting point to derive \sref{eq:int:FullSPDE} is the set of chemical reactions as laid out in Table~\ref{tab:my_label}. First, we introduce the column vector $X(t)=(u(t),v(t)))^T$, where $T$ indicates that we transpose the row vector, and the column vector $\mathcal{R}(X(t))$ with the propensities of the six reactions:
\begin{align*}
\mathcal{R}(X(t))=\left(a_1u(t),a_2u(t)v(t),\dfrac{a_3u(t)^2}{a_4+u(t)^2},a_5,\e c_1v(t),\e c_2u(t)\right)^T.    
\end{align*}
The associated stoichiometric matrix $\mathcal{S}$, which describes the change in $X(t)$ for each reaction, is then given by 
\begin{align}
\label{eq:Stoi}
\mathcal{S}=\begin{pmatrix}
-1&-1&1&1&0&0\\
0&0&0&0&-1&1\\
\end{pmatrix},   
\end{align}
see the last two columns of Table~\ref{tab:my_label}.
On top of these reactions, we assume that both variables also diffuse, so for a well-mixed solution in a large container we find the classic PDE
\begin{align}
\label{eq:RDE}
    \partial_t X=D\partial_{xx} X+\mathcal{S}\mathcal{R}(X),
\end{align}
where $D$ is a diagonal diffusion matrix with coefficients $D_u$ and $D_v$ on the diagonal~\cite{bressloff2014stochastic}. This PDE is identical to the RDE~\sref{eq:int:MainPDE} and describes the dynamics of $X(t)$, averaged over many individual reactions. When the number of reacting molecules is large enough, and when we zoom out far enough such that all individual molecules become effectively a density, the macroscopic PDE gives a good approximation of the microscopic behaviour. Statistically speaking, this means that the probability distribution of all possible states must be very sharply peaked around the average value described by the PDE, so the deviations from the mean can be ignored. 

\subsection{Motivating Example}
The assumption that we can ignore deviations from the mean is not always valid. For example, in population dynamics, we can write down birth-death models for several hundred individuals and with this number of individuals, random deviations from the mean are actually significant. To further exemplify, and to set the stage for the upcoming derivation, let us study such a simple discrete birth-death process: suppose a population is at time $t$ in state $X(t)$. In the next timestep $dt$, there are three possible outcomes: (i) the population grows by one individual with probability $b(X(t))dt$, (ii) the population decreases by one individual with probability $d(X(t))dt$, or (iii) nothing happens to the population with probability $1-b(X(t))dt-d(X(t))dt$.

Now, assume we have a continuous Stochastic Differential Equation (SDE)
\begin{align}
\label{eq:SDE}
    dx(t)=f(x(t))dt+g(x(t)) d  \beta_t,
\end{align}
where $\beta_t$ is Brownian motion, i.e. we can think of $d\beta_t$ as a random step with average zero and variance $dt$. We now ask the question: ``When is this continuous SDE a good approximation of the described discrete birth-death process?''. Or, more precisely, ``What should $f(x)$ and $g(x)$ be such that~\sref{eq:SDE} is a good approximation of the described discrete process?''. Given a solution $x$ of the SDE, we see that the average expected value at $x(t+dt)$ is approximated, at lowest order in $dt$, by
\begin{align*}
    E[x(t+dt)]=x(t)+f(x(t))dt+\mathcal{O}(dt^2).
\end{align*}  
For the described birth-death process, we have that the expectation is
\begin{align*}
E[X(t+dt)]=X(t)+[b(X(t))-d(X(t))]dt.    
\end{align*}
Hence, the average expected jump size in population is identical for the SDE~\sref{eq:SDE} and the birth-death process if we take $f(x):=b(x)-d(x)$. 

Next, we compute the deviation from the mean of the SDE~\sref{eq:SDE}
 \begin{align*}
     \text{Var}[x(t+dt)]=\text{Var}[g(x(t))d\beta_t]+\mathcal{O}(dt^2)=g(x(t))^2dt+\mathcal{O}(dt^2)\,,
 \end{align*}
 while this deviation for the birth-death process is
 \begin{align*}
 \text{Var}[X(t+dt)]=b(x(t))+d(x(t))+\mathcal{O}(dt^2).    
 \end{align*}
Therefore, to make these deviations coincide at first order in $dt$, we must take $g(x):=\sqrt{b(x)+d(x)}$. Hence, the process $x(t)$ described in~\sref{eq:SDE}, which is continuous in population size and time, is a good approximation of the discrete process $X(t)$ when
\begin{align}
    \label{eq:bdp}
    dx(t)=(b(x(t))-d(x(t)))dt+\sqrt{b(x(t))+d(x(t))}d\beta_t.
\end{align}
The stochastic process $x(t)$ shares the average and variance with $X(t)$ but differs in other points. Higher order moments of $x(t)$ and $X(t)$ will not be identical and $x(t)$ can become negative, even when $b$ and $d$ are chosen such that this is not possible in the discrete model. 

In order to link the SDE above to chemical reactions, we make the following observation. The birth of an individual can be thought of as the chemical reaction $\emptyset\to X$ with propensity $b(X)$ and stoichiometric value $1$, while the death of an individual can be seen as the chemical reaction $X\to\emptyset$ with propensity $d(X)$ and stoichiometric value $-1$. Next, we make an assumption which is called the \textit{leap condition}~\cite{bressloff2014stochastic}. That is, we assume that, given a state $X(t)$, enough reactions happen in the interval $[t,t+dt]$ to describe the average jump size in $[t,t+dt]$ by a Poisson process whose parameters depend on $X(t)$. With this leap condition assumption, we implicitly also assume that $X(t)$ is a good approximation of the solution in the whole time interval $[t,t+dt]$. We now turn the discrete process $X(t)$ into a continuous process $x(t)$ by approximating the discrete Poisson process by a continuous Gaussian, see~\cite{kim2017stochastic} for details. This approach results in an SDE similar to the SDE~\sref{eq:bdp}:
\begin{align}
    \label{eq:bdp2}
    dx(t)=(b(x(t))-d(x(t)))dt+\sqrt{b(x(t))}d\beta^1_t-\sqrt{d(x(t))}d\beta^2_t,
\end{align}
for two independent Brownian motions $\beta^1_t$ and $\beta_t^2$. Although visually different from~\sref{eq:bdp}, both SDEs have a noise term that is Gaussian with identical average and variance. Therefore, both SDEs describe the same stochastic process and hence we can say that \sref{eq:bdp} and~\sref{eq:bdp2} are equivalent. 

\subsection{Derivation of the CLE}
We have now gained some intuition for linking more general discrete chemical reactions to continuous S(P)DEs:
if we have $M$ different molecules in a vector $X(t)$ with diffusion matrix $D$, $N$ reactions given by a vector $\mathcal{R}(X(t))$ and a stoichiometric matrix $\mathcal{S}$, then the continuous SPDE for $X(t)$ is given by 
\begin{align}
\label{eq:cle}
    dX(t)=\left(D\partial_{xx} X(t)+\mathcal{S}\mathcal{R}(X(t))\right)dt+\frac{1}{\sqrt{\Omega}}\mathcal{S}\sqrt{\text{diag}(\mathcal{R}(X(t)))}dW_t+\frac{1}{\sqrt{\Omega}}\partial_x\sqrt{2DX(t)}d\tilde W_t,
\end{align}
see~\cite{bressloff2014stochastic, kim2017stochastic}. The equation is made of two parts, a local equation that describes the kinetics as in SDE \sref{eq:bdp}
\begin{align}
    dX(t)=\mathcal{S}\mathcal{R}(X(t))dt+\frac{1}{\sqrt{\Omega}}\mathcal{S}\sqrt{\text{diag}(\mathcal{R}(X(t)))}dW_t
\end{align}
and a stochastic diffusion equation 
\begin{align}
    dX(t)=D\partial_{xx} X(t)+\frac{1}{\sqrt{\Omega}}\partial_x\sqrt{2DX(t)}d\tilde W_t,
\end{align}
as derived in \cite{dogan2011derivation}. Here, $dW_t$ and $d\tilde W_t$ are two independent vectors with space-time white noise. The vector $dW_t$ has $N$ components coming from the $N$ reactions, while $d\tilde W_t$ has the dimension $M$ of $X(t)$. SPDE~\sref{eq:cle} is known as the Chemical Langevin Equation (CLE)~\cite{gillespie2000chemical}. The vector $X(t)$ now describes the densities of the molecules involved, not the actual number of molecules. How well the discrete number of molecules is approximated by a density is determined by the scale parameter $\Omega$ and is in that sense a measure for the \emph{noisiness} of the system. In the no-noise limit $\Omega\to\infty$, we recover the classic RDE~\sref{eq:RDE}. In contrast, for small $\Omega$ the dynamics of the discrete process is dominated by random events and the discrete process should be described in full detail by a \emph{chemical master equation}~\cite{gillespie1992rigorous}. The CLE can be understood as the lowest order approximation of the  chemical master equation for large $\Omega$, see for more details~\cite{bressloff2014stochastic}. For an overview of all different paths leading from molecular kinetics to (S)PDEs, see~\cite[Fig. 3.4]{lei2021systems}. It is important to realise that SPDE \sref{eq:cle} does not necessarily inherit all the statistical properties of the chemical master equation, only averages and variances. Another potential issue is that it does not necessarily ensures positivity of the solutions. 

Just as \sref{eq:bdp} and~\sref{eq:bdp2} are identical, 
we can rewrite \sref{eq:cle} in the following way:
\begin{align}
\label{eq:cle2}
    dX(t)=(D\partial_{xx} X(t)+\mathcal{S}\mathcal{R}(X(t)))dt+\frac{1}{\sqrt{\Omega}}\sqrt{\mathcal{S}\text{diag}(\mathcal{R}(X(t)))\mathcal{S}^T}dW_t+\frac{1}{\sqrt{\Omega}}\partial_x\sqrt{2DX(t)}d\tilde W_t.
\end{align}
This time, the noise vector $dW_t$ has just $M$ components, reducing the number of random vectors that must be generated (when $M<N$). The downside is that the computation of $\sqrt{\mathcal{S}\text{diag}(\mathcal{R}(X))\mathcal{S}^T}$ is in general numerically more expensive then the computation of $\mathcal{S}\sqrt{\text{diag}(\mathcal{R}(X))}$. However, in the present setting, there are no connections between the two variables in the stoichiometric matrix $\mathcal{S}$~\sref{eq:Stoi} and the matrix $\mathcal{S}\text{diag}(\mathcal{R}(X))\mathcal{S}^T$ is thus diagonal, making the computation of the square root trivial. 

Note that once we have the CLE~\sref{eq:cle2}, it can be applied to any set of chemical reactions and can therefore have widespread use. For example, we can now return to Table~\ref{tab:my_label} and apply the CLE to these reactions, which results in
\begin{align}
\label{eq:FullSPDE}
    \begin{split}
        du&=\left(D_u\partial_{xx}u-(a_1+a_2v)u+\frac{a_3u^2}{a_4+u^2}+a_5\right)dt+\sigma\sqrt{(a_1+a_2v)u+\frac{a_3u^2}{a_4+u^2}+a_5}dW^1_t\\
        &\qquad\qquad+\sigma\partial_x\sqrt{2D_uu}d\tilde W^1_t,\\
        dv&=\left(D_v\partial_{xx}v+\e(-c_1v+c_2u)\right)dt+\sigma\sqrt{\e(-c_1v+c_2u)}dW^2_t+\sigma\partial_x\sqrt{2D_vv}d\tilde W^2_t.
    \end{split}
\end{align}
For notational convenience,  we replaced $1/\sqrt{\Omega}$ by a small parameter $\sigma$, resulting in the SPDE~\sref{eq:int:FullSPDE} from the introduction. In the remainder of this work, we will study the SPDE above, mainly using numerical techniques. 

\begin{remark}
\label{rem:ill}
It is important to realise that the SPDE above does not have a function-valued solution in general. The term $\partial_x\sqrt{2DX(t)}d\tilde W_t$ makes the equation ill-posed and solutions can only be understood in terms of distributions. Therefore, it is not a priori clear if the numerical solutions shown in the next section converge to a solution of the SPDE when the spatio-temporal discretisations $dx$ and $dt$ are sent to 0. In $\S$\ref{sec:WvsP}, we will discuss the implications of omitting this term on the wave dynamics.
\end{remark}

\section{Simulations}
\label{sec:sim}
In this section, we will numerically investigate the PDE~\sref{eq:int:MainPDE} and SPDE~\sref{eq:FullSPDE}. We investigate three of the main building blocks of the PDE dynamics: localised standing waves, localised travelling waves and time-periodic solutions, together with their counterparts in the SPDE. However, before we can investigate the dynamics, we must first establish some basic properties of the (S)PDE, like the existence, uniqueness and stability of the background state(s).

\begin{figure}
\begin{subfigure}{0.49\textwidth}
        \centering
		\def\svgwidth{\columnwidth}
    		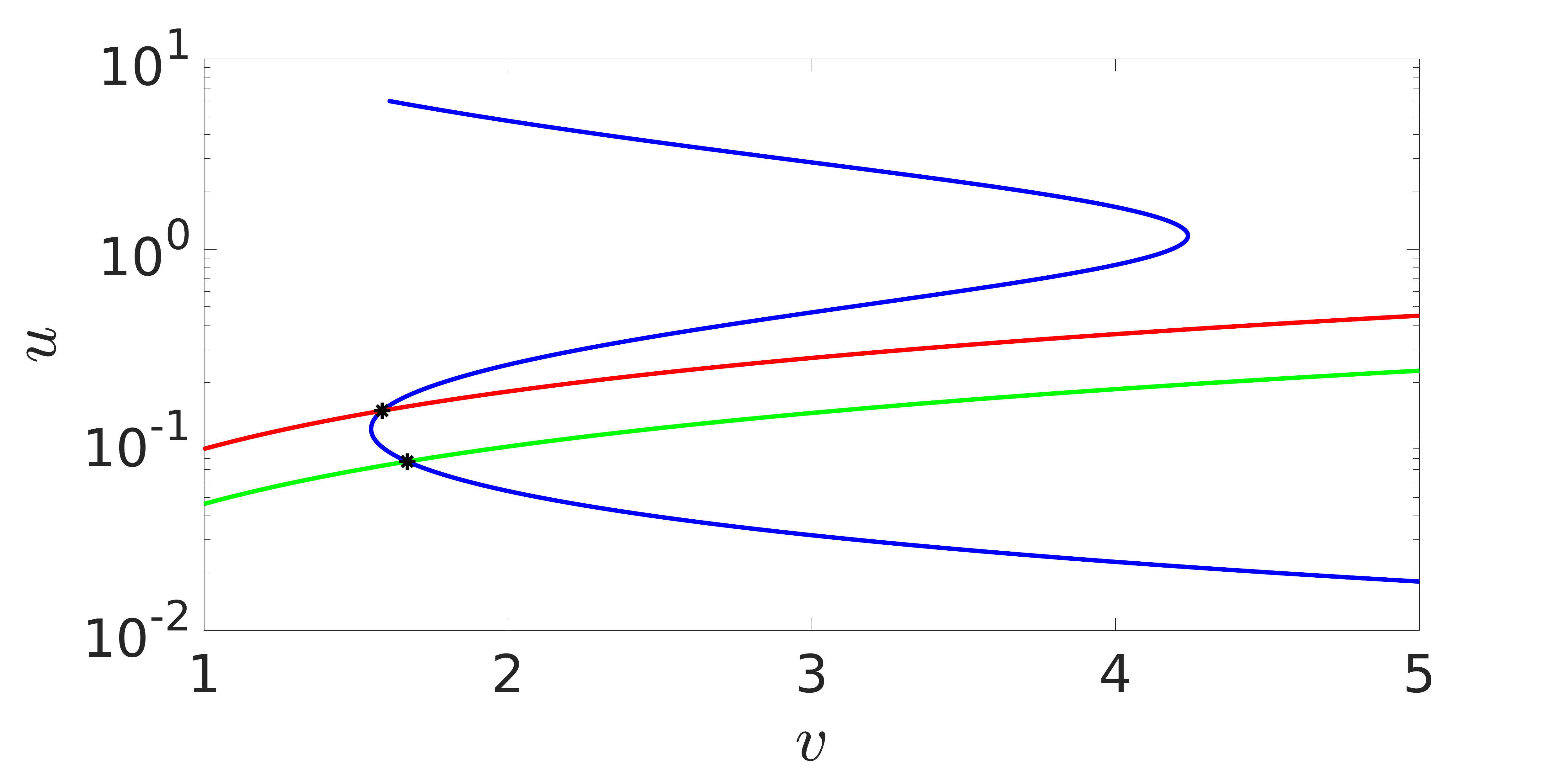
      \caption{}
    \label{fig:PPc1}
\end{subfigure}
\begin{subfigure}{0.49\textwidth}
    \centering
 		\def\svgwidth{\columnwidth}
    		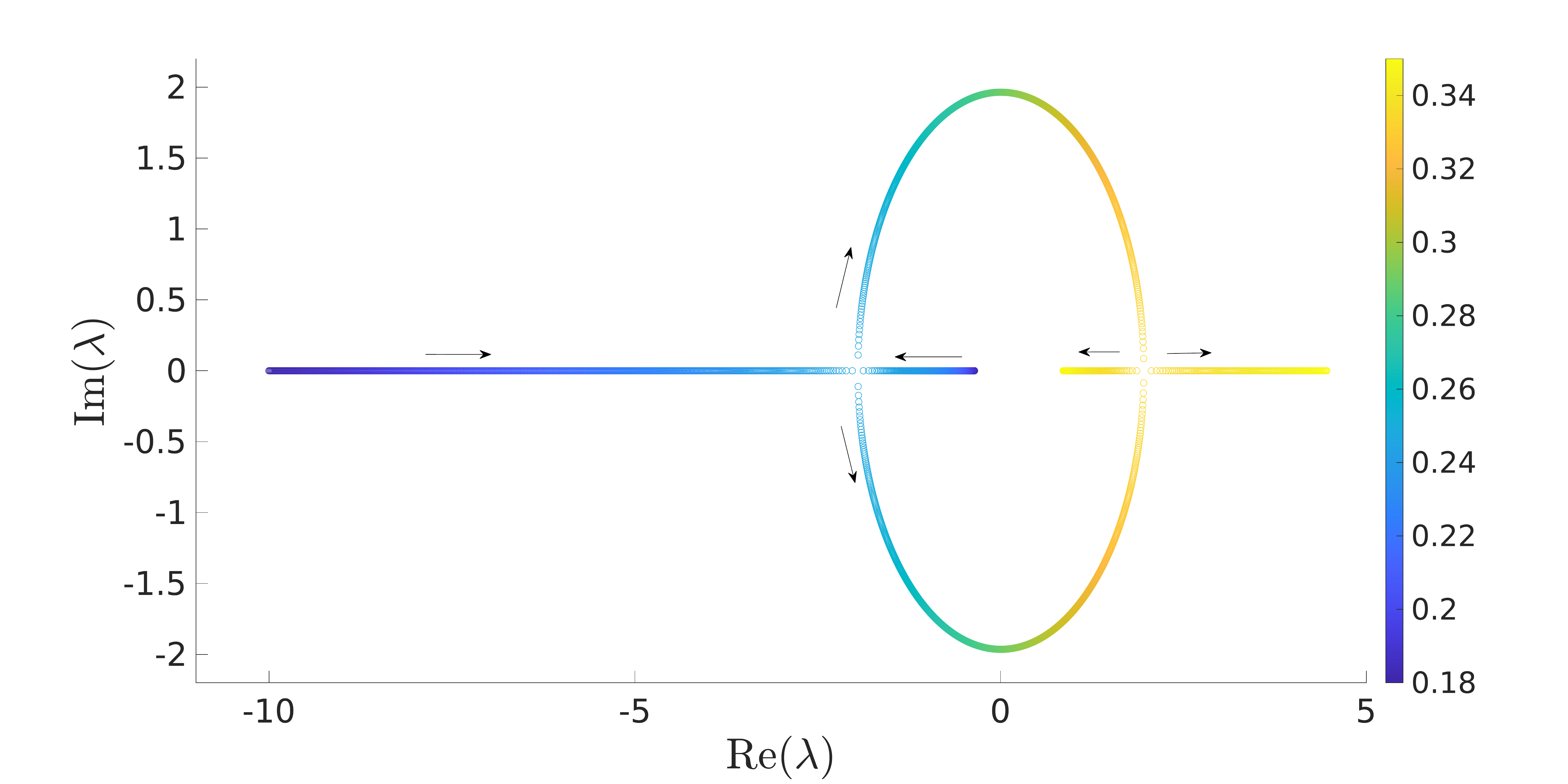
    \caption{}
    \label{fig:bifdiag}    
\end{subfigure}
    \caption{(a) The green line is the $v$-nullcline for $c_1=0.18$, while the red line is the nullcline for $c_1=0.35$. The blue line is the $u$-nullcline, independent of $c_1$. The $u$-axis is plotted logarithmically to better highlight the shape of the nullcline for small $u$.  Note how the background state moved around the fold. (b) Visual representation of the evolution of the two (complex) eigenvalues of the Jacobian matrix~\sref{eq:sim:jac} for $c_1$ varying from $0.18$ (dark blue) to $0.35$ (yellow), following the black arrows. The other parameters are fixed at $a_1=0.167$, $a_2=16.67$, $a_3=167$, $a_4=1.44$, $a_5=1.47$, $\e=0.52$ and $c_2=3.9$. }
\label{fig:PPbifdiag}
\end{figure}

Since we are interested in localised waves and expect the activator to be \emph{in rest} otherwise, we need for the existence of these localised waves that the spatially homogeneous background state is stable. In contrast, for the time-periodic solutions, we expect the background state to be unstable such that continuous excitations of the background state can happen. The possible background states $(u^*,v^*)$ of~\sref{eq:int:MainPDE} are given by the positive real solutions of the $u$-nullcline  and $v$-nullcline
\begin{align}
\label{eq:sim:null}
    0= -a_1 u-a_2uv+\frac{a_3u^2}{a_4+u^2}+a_5\,, \qquad 
    0=\e(-c_1v+c_2u) \,.
\end{align}
See Figure~\ref{fig:PPc1} for a typical representation of the shape of the nullclines. 
Since the system parameters are all assumed to be positive, this is equivalent to finding the positive solutions $u^*$ of
\begin{align*}
   -\frac{a_2 c_2}{c_1}u^4 -a_1u^3 + \left(a_3 + a_5 -\frac{a_2 a_4 c_2}{c_1}\right)u^2  -a_1 a_4 u +a_5a_4=0\,, 
\end{align*}
with $v^*= c_2 u^*/c_1$. Due to the complexity of the general solution formula for quartic polynomials, it is not feasible to write down its solutions explicitly.
However, by Descartes' rule of signs~\cite{descartes2020discours} we know that there is only one positive real root if $c_1(a_3 + a_5) < a_2 a_4 c_2$ and one or three positive real roots otherwise\footnote{Note that the origin $(0,0)$ is only a background state if $a_5=0$.}. The stability of a background state $(u^*,v^*)$ is then determined by the eigenvalues of the associated Jacobian matrix
\begin{align}
\label{eq:sim:jac}
    J(u^*,v^*)=
    \begin{pmatrix}
    -a_1-a_2v^*-\dfrac{2a_3a_4u^*}{(a_4+(u^*)^2)^2}& -a_2u^*\\
    \e c_2&-\e c_1
    \end{pmatrix}.
\end{align}
Since we do not have an explicit formula for $(u^*,v^*)$, we must compute these eigenvalues numerically. For example, when we allow one free parameter, e.g. $c_1$, and fix the other values, then we can compute the background states and the associated eigenvalues of the Jacobian matrix. Taking the parameter values $a_1=0.167$, $a_2=16.67$, $a_3=167$, $a_4=1.44$, $a_5=1.47$, $\e=0.52$ and $c_2=3.9$ from~\cite{bhattacharya2020traveling} and letting $c_1$ range from $0.18$ to $0.35$, such that $c_1(a_3 + a_5) < a_2 a_4 c_2$, results in one admissible positive background state ranging from $(u^*,v^*)\approx(0.077, 1.669)$ to $(u^*,v^*)\approx(0.142, 1.586)$. Initially, for the lower values of $c_1$, the eigenvalues are real and negative, resulting in a stable background state. Increasing the value of $c_1$ to approximately $0.25$ results in complex eigenvalues, still with negative real parts. When we further increase the value of $c_1$ to approximately $0.29$, both eigenvalues cross the imaginary axis, i.e. the background state undergoes a Hopf bifurcation and we expect to see time-periodic solutions. See Figure~\ref{fig:bifdiag} for a visual representation of the evolution of the eigenvalues. In Figure~\ref{fig:PPc1} we show the nullclines for $c_1=0.18$ and $c_1=0.35$. The unique background state moved along the fold in the $u$-nullcline and as long as the background state is in between the two folds, the fixed point is unstable. 

In the next sections, we will study localised standing and travelling waves for the same parameter set with $c_1<0.25$ and for time-periodic solutions with $c_1>0.29$. The complex dynamics of pulse adding for $c_1$-values in the intermediate regime between these two boundary values, where the eigenvalues of the Jacobian are stable but complex-valued, is outside the scope of this work, see for example~\cite{carter2016stability} for more information.

So far, we only looked at background states, which are spatially homogeneous. However, we are interested in spatially nonhomogeneous patterns. By definition, a localised wave is a fixed profile $(\Phi_u,\Phi_v)$ that moves with a fixed speed $c$ (possibly zero). Therefore, when we change the spatial coordinate $x$ to $\xi=x-ct$ using the chain rule, the profile $(\Phi_u,\Phi_v)$ is a stationary solution of the following shifted Ordinary Differential Equation (ODE): 
\begin{align}
   \label{eq:sim:bvp}
    \begin{split}
        0&=D_u\partial_{\xi\xi}\Phi_u+c\partial_{\xi}\Phi_u-(a_1+a_2\Phi_v)\Phi_u+\frac{a_3\Phi_u^2}{a_4+\Phi_u^2}+a_5,\\
        0&=D_v\partial_{\xi\xi}\Phi_v+c\partial_{\xi}\Phi_v+\e(-c_1\Phi_v+c_2\Phi_u).
    \end{split}
\end{align}
This ODE problem can be solved using numerical fixed-point algorithms. For these algorithms, a crude starting point is needed for the profile and the value of $c$,  which can come from a PDE simulation. Note that this problem is translation invariant, meaning that we find a one-dimensional family of travelling waves, all shifted versions of each other. Hence, for the solver to converge, an extra condition to fix the location of the wave is necessary.

\subsection{Standing Waves}
\label{sim:SW}
\begin{figure}
\begin{subfigure}{.49\textwidth}
  \centering
 		\def\svgwidth{\columnwidth}
    		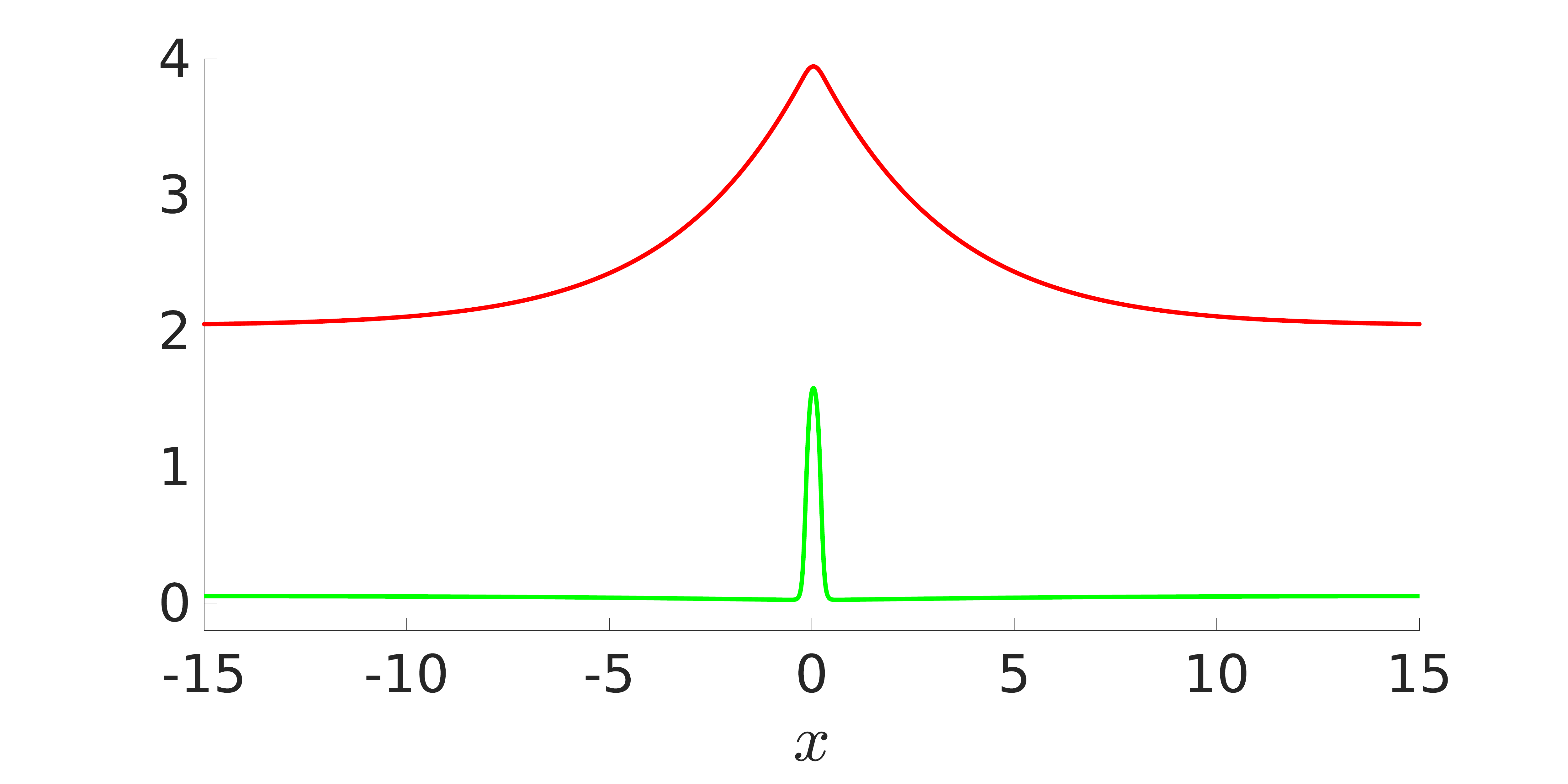
  \caption{}
    \label{fig:Sim:Pulse}
\end{subfigure}
\begin{subfigure}{.49\textwidth}
  \centering
 		\def\svgwidth{\columnwidth}
    		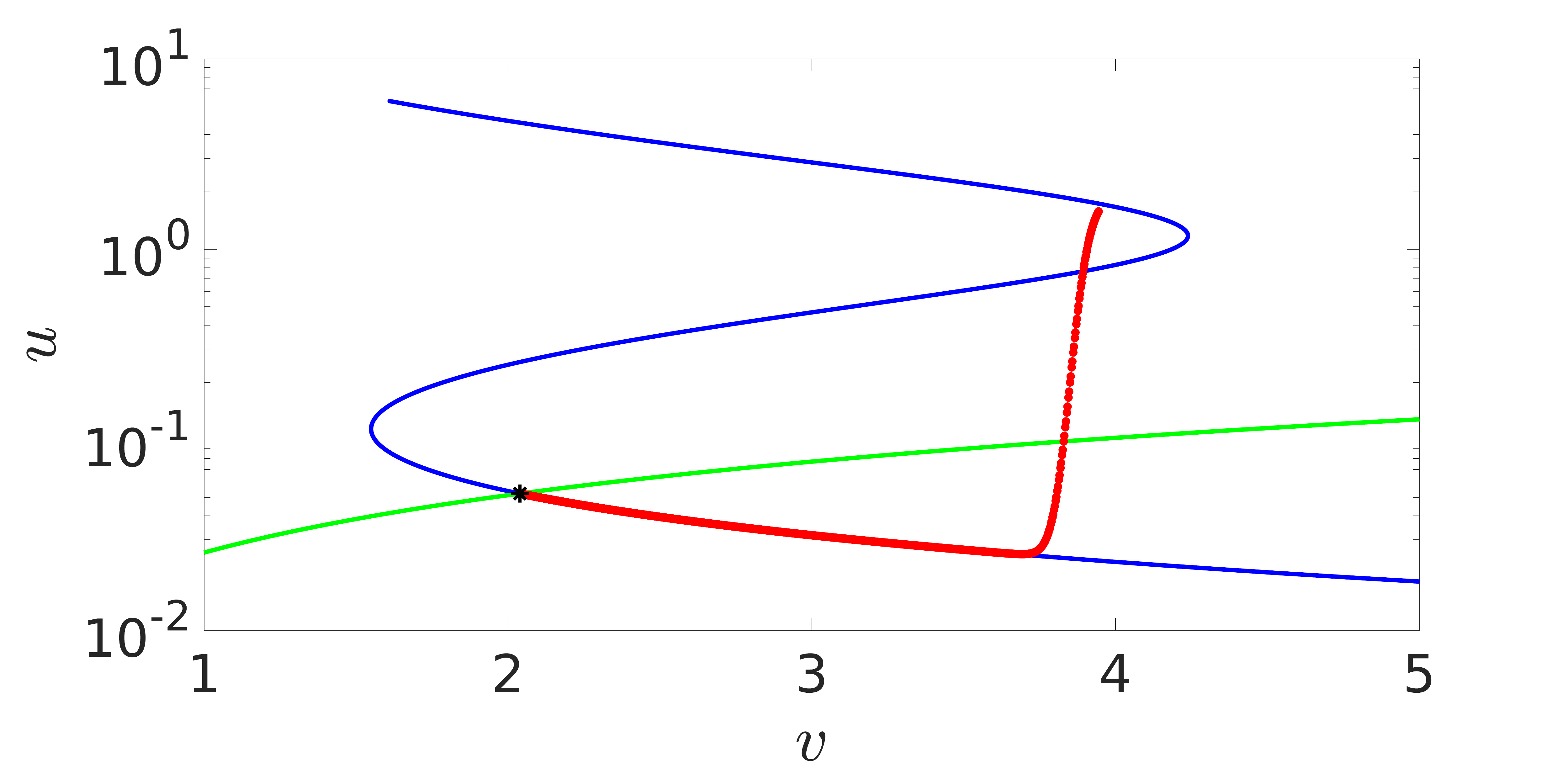
  \caption{}
    \label{fig:Sim:PulsePP}
\end{subfigure}
\caption{Figure~(a) shows a localised standing wave solution to ODE~\sref{eq:sim:bvp}, found numerically with Matlab's \texttt{fsolve}. The green curve is the $u$-component and the red curve is the $v$-component. In Figure~(b), the $u$-nullcline (blue) and $v$-nullcline (green) of~\sref{eq:int:MainPDE} are shown together with the $v$-$u$ phase plane of the standing wave from Figure~(a). The phase plane is plotted on a semi-log scale to better highlight the dynamics for small $u$. We observe that the standing wave starts from the background state (indicated by an asterisk) and initially follows the lower branch of the $u$-nullcline before jumping to the upper branch of the $u$-nullcline and follows the same track back to the background state. The system parameters are taken from~\cite{bhattacharya2020traveling} and set to $D_u=0.1$, $D_v=1$, $a_1=0.167$, $a_2=16.67$, $a_3=167$, $a_4=1.44$, $a_5=1.47$, $\e=0.52$, $c_1=0.1$, and  $c_2=3.9$. 
}
\label{fig:sim:PulseBoth}
\end{figure}

In this section, we will study standing waves, which means we look for solutions of~\sref{eq:sim:bvp} with $c=0$. A solution to this ODE is shown in Figure~\ref{fig:Sim:Pulse}. We observe that both components $u$ and $v$ indeed start at and return to their background state $(u^*,v^*) \approx (0.0523, 2.0394)$. We observe that the activator $u$ changes rapidly in a small region in the spatial domain and we, therefore, call the activator $u$ the {\emph{fast variable}}. On the other hand, the inhibitor $v$ is the {\emph{slow variable}} as it changes more gradually over a larger spatial distance. Figure~\ref{fig:Sim:PulsePP} shows the corresponding phase plane. The majority of the spatial dynamics happens near the lower branch of the $u$-nullcline before it has a fast excursion from the lower branch to the upper branch of this nullcline and, by the symmetry $x \mapsto -x$ of the ODE~\sref{eq:sim:bvp}, it then returns back to the lower branch in a similar fashion. The fact that both components of the standing pulse evolve on a different spatial scale allows us to mathematically analyse this standing pulse, see~Appendix~\ref{sec:analysis}. 
For instance, the value $\bar v$ at which the activator $u$ makes a sharp transition (approximately $3.8$ in Figure~\ref{fig:Sim:PulsePP}), can be approximated by the algebraic relation~\sref{eq:an:barv}. The analysis also explains why the solution trajectory in the phase plane closely follows the lower branch of the $u$-nullcline for the most part of the trajectory. 

\begin{figure}
\begin{subfigure}{.49\textwidth}
  \centering
 		\def\svgwidth{\columnwidth}
    		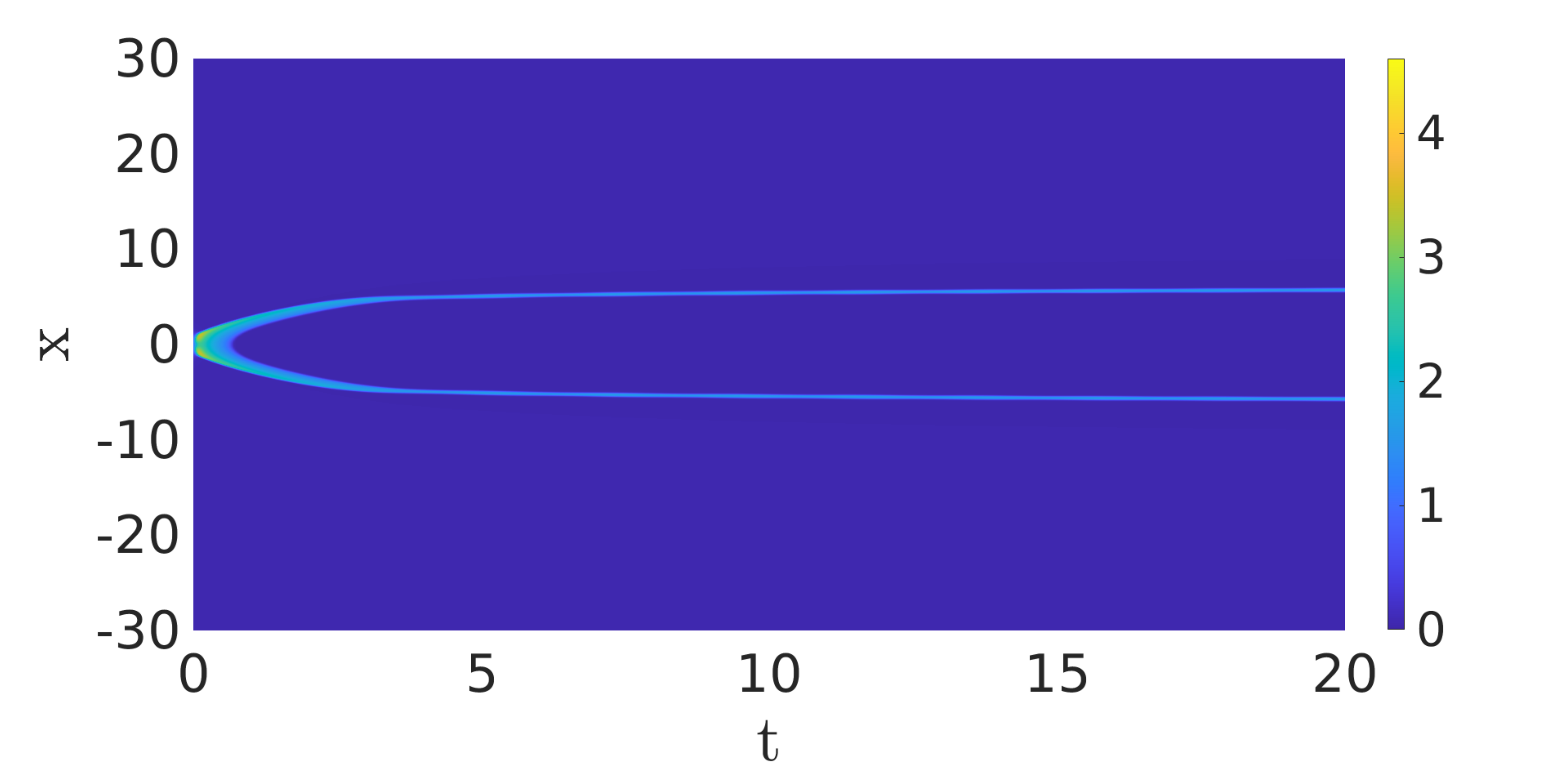
  \caption{}
    \label{fig:Sim:DPulseU}
\end{subfigure}
\begin{subfigure}{.49\textwidth}
  \centering
 		\def\svgwidth{\columnwidth}
    		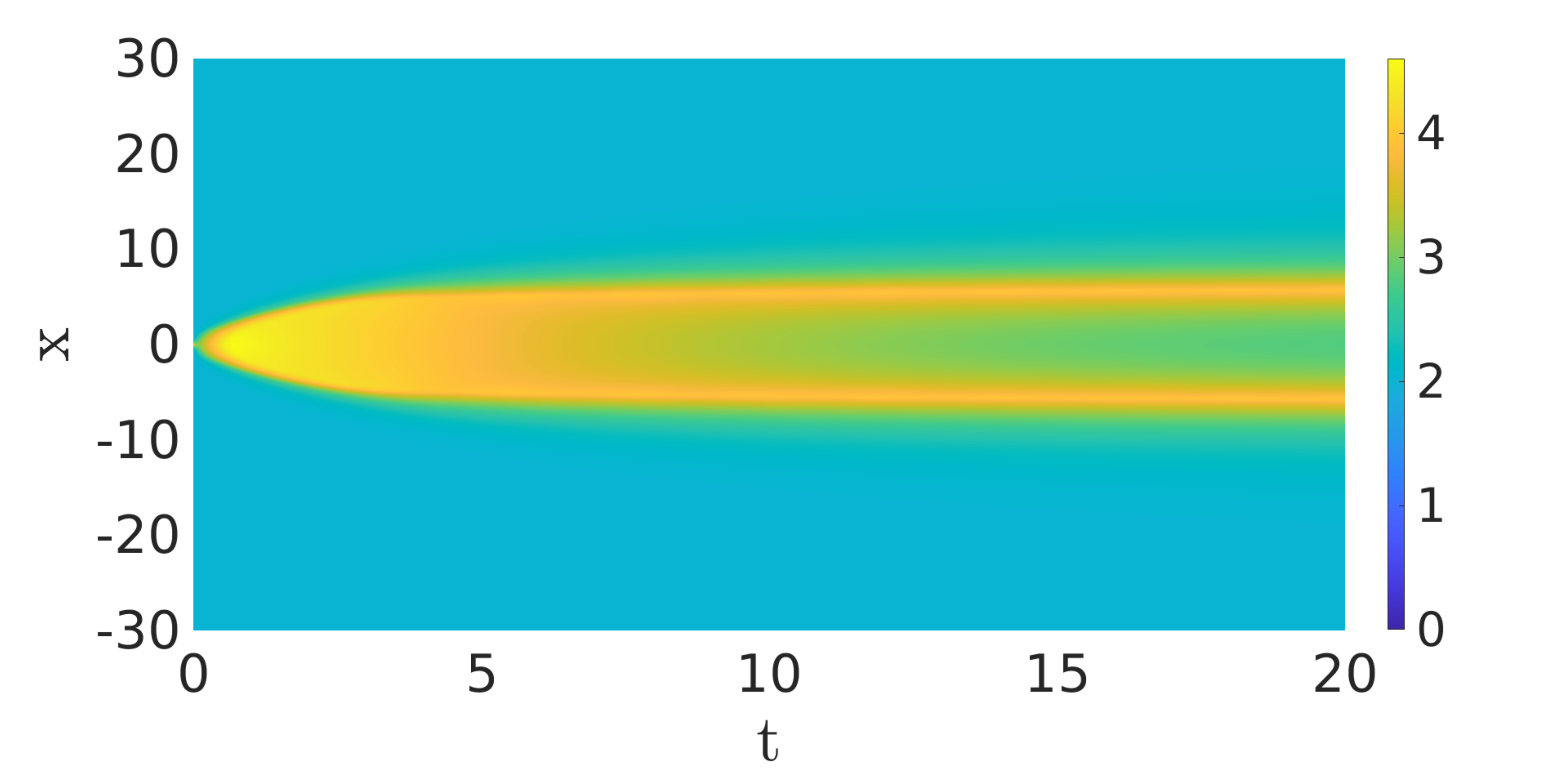
  \caption{}
    \label{fig:Sim:DpulseV}
\end{subfigure}
\caption{Simulation of the PDE~\sref{eq:int:MainPDE}, Figure~(a) shows the activator $u$ and Figure~(b) the inhibitor $v$ with an initial condition as described in the main text. The same parameters are used as in Figure~\ref{fig:sim:PulseBoth}. Note that the $v$-component does not return to its rest state in the region between the two pulses.}
\label{fig:Sim:Dpulse}
\end{figure}

By assumption, the standing wave in Figure~\ref{fig:Sim:Pulse} is a stationary solution of the PDE~\sref{eq:int:MainPDE}. This can be confirmed by using the wave from the ODE as the initial condition for a PDE simulation (not shown). However, we are not likely to find this single standing wave in a PDE simulation without a fine-tuned initial condition. As an example, we use for the simulation the initial condition $u_0=u^*+e^{-x^2}$ and $v_0=v^*+2/\cosh^2(5x)$ as a crude approximation of the wave. The resulting simulation is shown in Figure~\ref{fig:Sim:Dpulse}. This initial condition splits in, what appears to be, two well-separated localised standing waves\footnote{We also observe the evolution of the initial condition back to the stable background state $(u^*,v^*)$, especially for initial conditions with smaller amplitudes. Simulations are not shown.}. However, the plot of the slow $v$-component makes clear that this is not the case, and that the two standing waves are connected through the slow component, i.e. the slow component is not in its rest state in between the two standing waves. For more details on the numerics of the (S)PDE simulations, see Appendix~\ref{sec:app:num}. 

\begin{figure}
\begin{subfigure}{.49\textwidth}
  \centering
 		\def\svgwidth{\columnwidth}
    		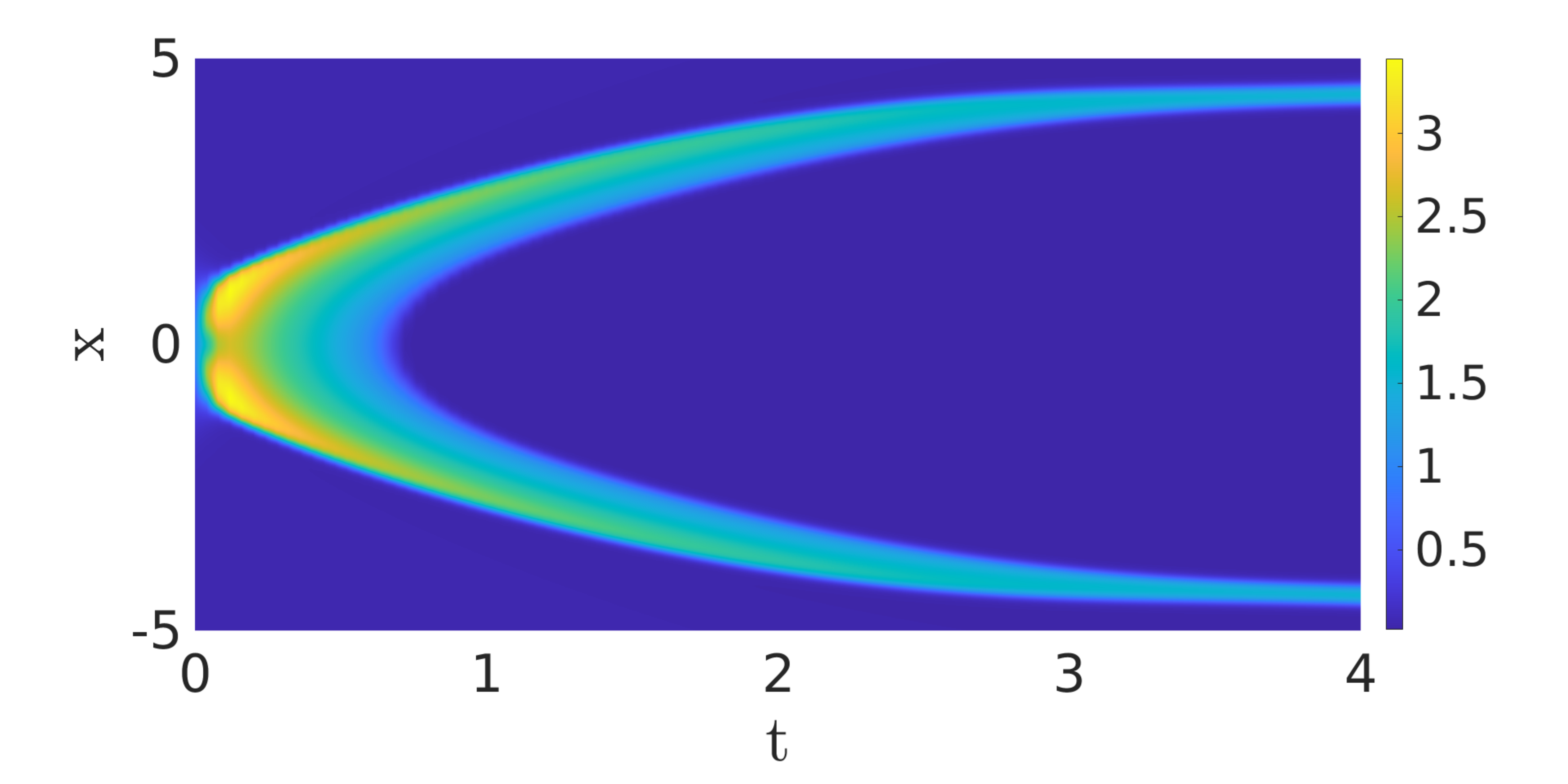
  \caption{}
    \label{fig:Sim:DPulseZoom}
\end{subfigure}
\begin{subfigure}{.49\textwidth}
  \centering
 		\def\svgwidth{\columnwidth}
    		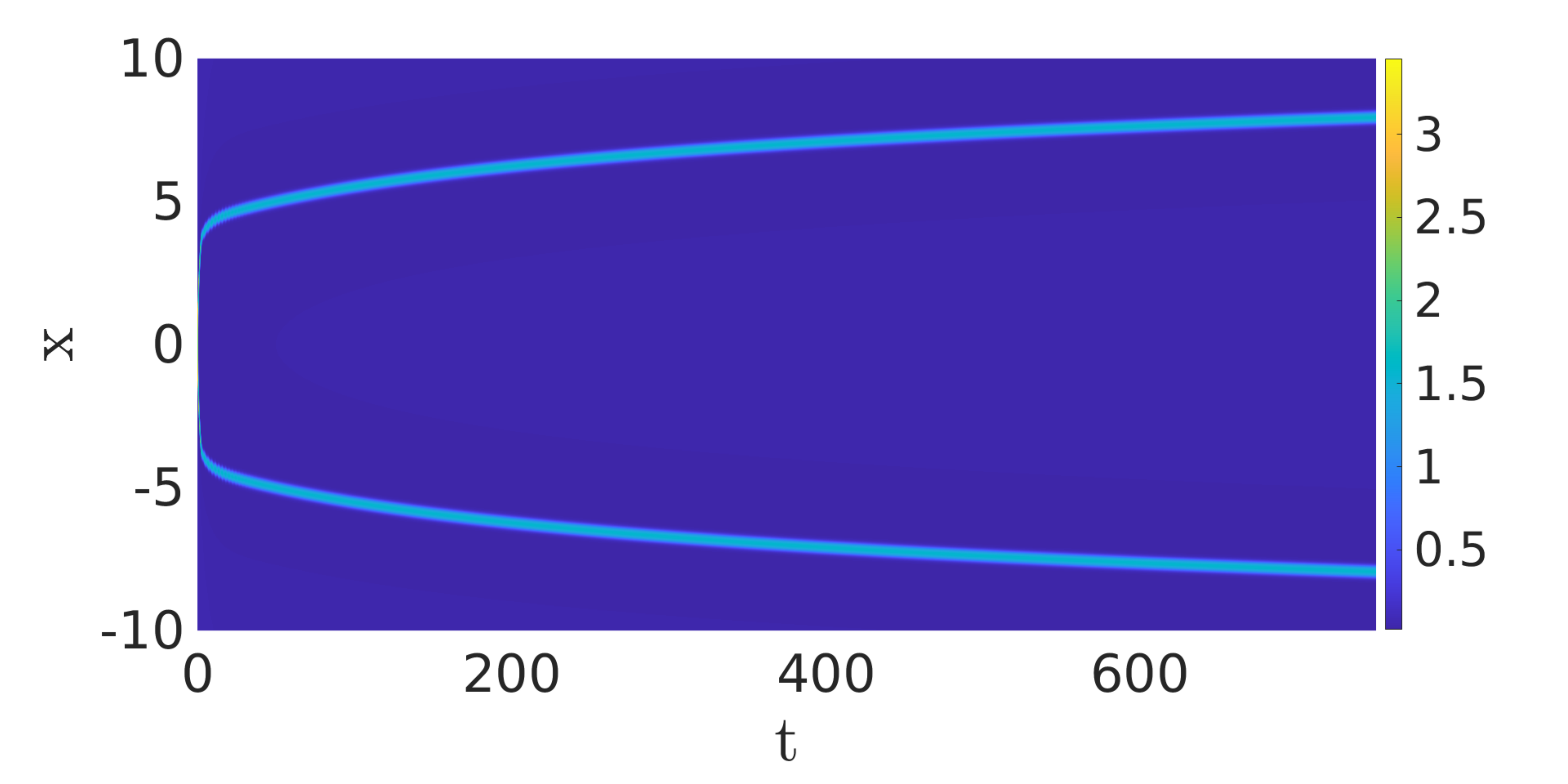
  \caption{}
    \label{fig:Sim:DpulseLong}
\end{subfigure}
\caption{Same simulation as in Figure~\ref{fig:Sim:Dpulse}, but on different time scales. Figure~(a) shows the $u$-component, zoomed in to highlight the short-time dynamics, while Figure~(b) shows the long-time dynamics of $u$ highlighting the pulse splitting phenomenon. Both simulations were done on a larger grid $[-60,60]$, so the waves would not affect each other on the other side of the domain on this large time scale.}
\label{fig:Sim:DPulseLZ}
\end{figure}

The interaction between the two standing waves in Figure~\ref{fig:Sim:Dpulse} through the slow $v$-component makes that the two standing waves repel each other on a very slow timescale as is made clear by taking long integration times, see Figure~\ref{fig:Sim:DpulseLong}. On an infinite domain, the two standing waves slowly drift apart forever, but on a periodic domain,  we can expect them to stabilise once they are at an equal distance on both sides. On the timescales of biological processes, this slow continuous splitting is probably not relevant and on short timescales, the term `standing waves' for the solution at later times in Figure~\ref{fig:Sim:Dpulse} is biologically justifiable. Furthermore, note that for our understanding of the presented dynamics, it is essential to look at both components simultaneously. In other words, for our understanding of Figure~\ref{fig:Sim:DPulseU} it is essential to also look at Figure~\ref{fig:Sim:DpulseV}.

We now take a closer look at the short-time dynamics presented in Figure~\ref{fig:Sim:DPulseZoom}. In~\cite{bhattacharya2020traveling}, this splitting of the initial condition is described as two counter-propagating travelling waves, sometimes called trigger waves \cite{gelens2014spatial}. By the formal mathematical definition, a travelling wave is a fixed profile moving with a fixed speed, i.e. a solution of~\sref{eq:sim:bvp}. Therefore, mathematically speaking, these do not classify as travelling waves. Instead, what we observe here would be classified as transient dynamics and pulse splitting. However, it is clear that at $t=0$, the activity of $u$ is around $x=0$, and after some time it moved to two different places, justifying the term `travelling'. If we adopt the terms `standing' and `travelling', it is clear from Figure~\ref{fig:Sim:DPulseU} that around $t=3$ a transition occurs from travelling to standing. 

\begin{figure}
\begin{subfigure}{.49\textwidth}
  \centering
 		\def\svgwidth{\columnwidth}
    		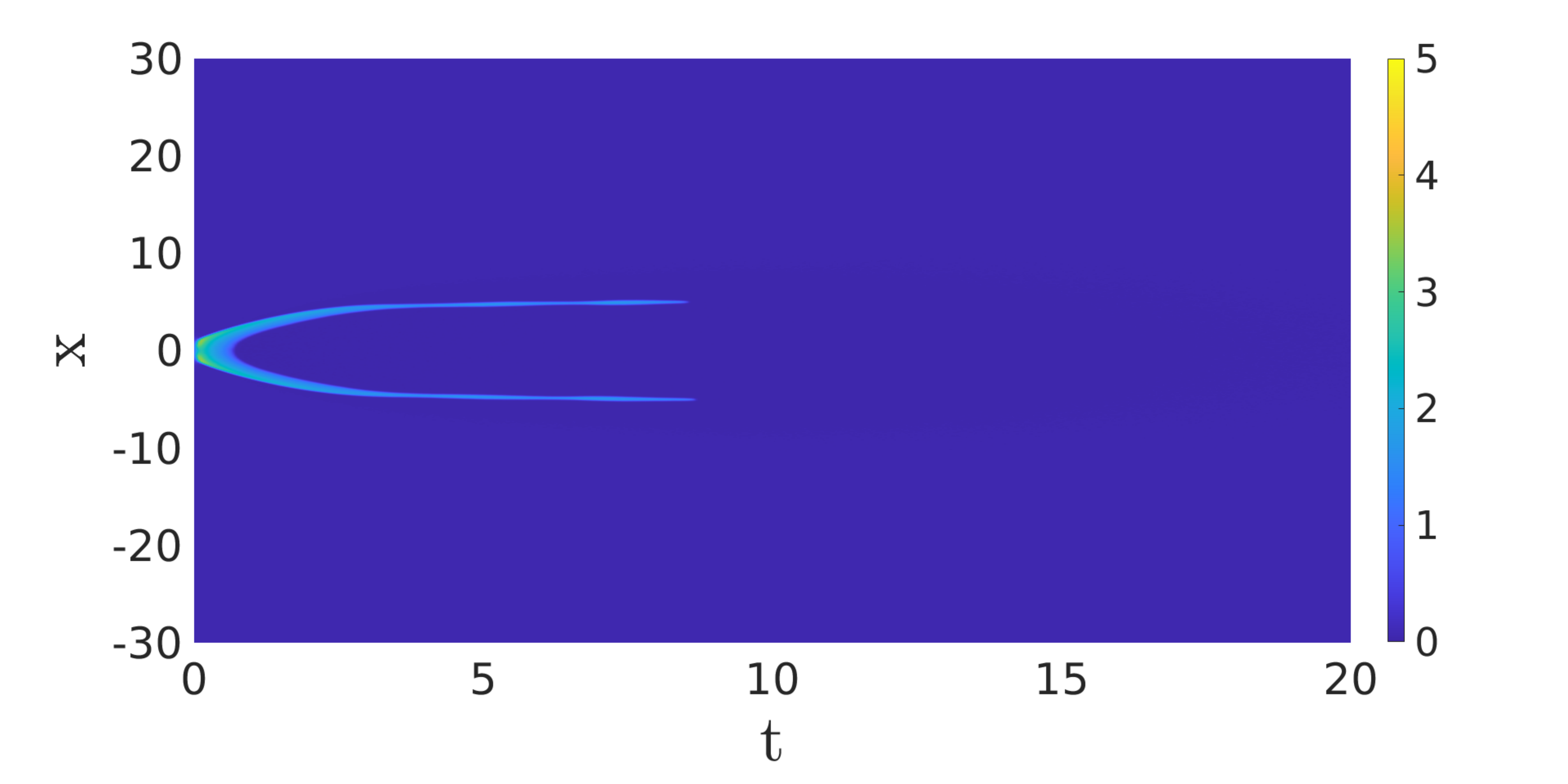
  \caption{$\sigma=0.002$}
    \label{fig:Sim:PulseNoise002}
\end{subfigure}
 \begin{subfigure}{.49\textwidth}
  \centering
 		\def\svgwidth{\columnwidth}
    		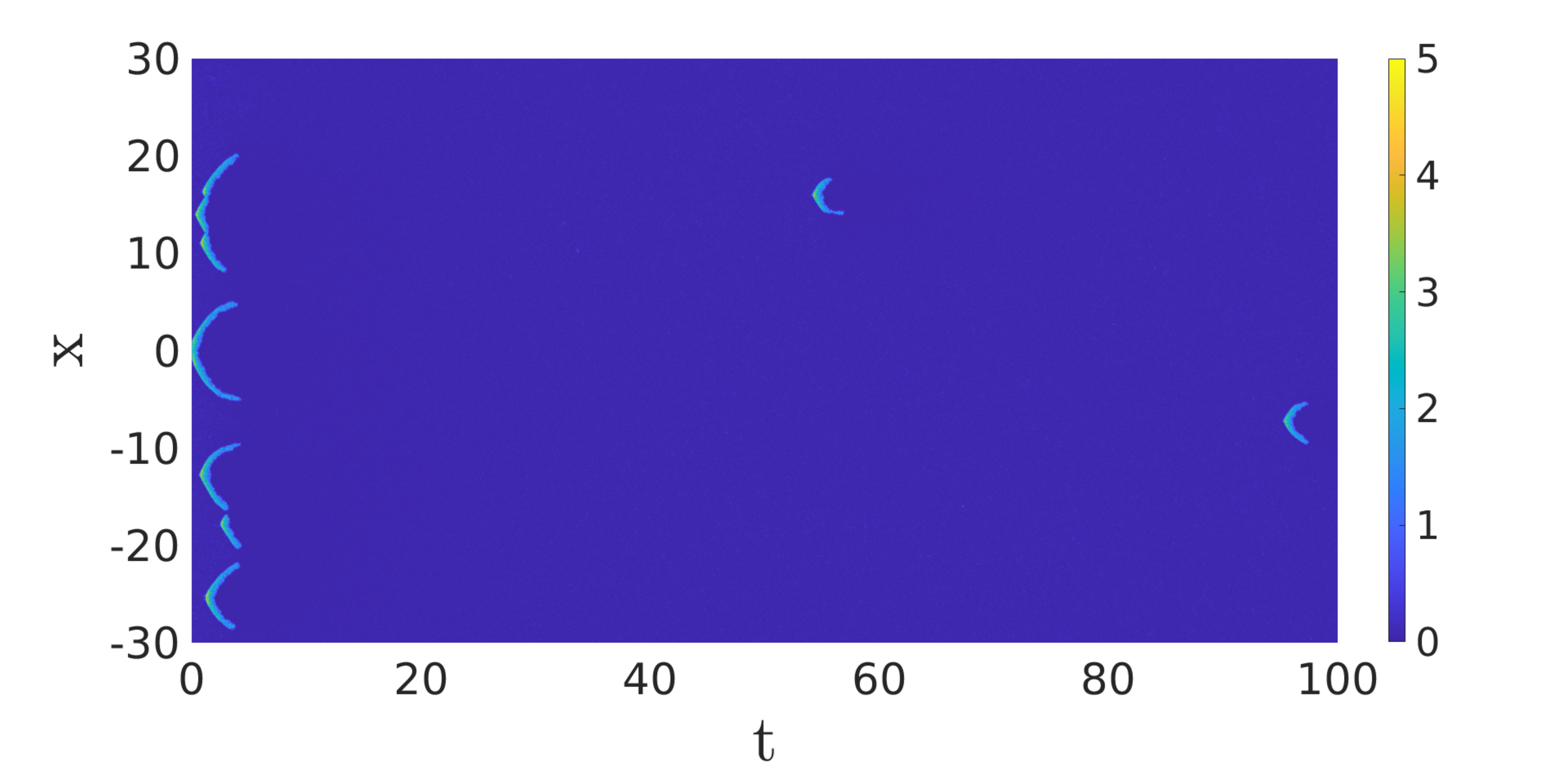
  \caption{$\sigma=0.035$}
    \label{fig:Sim:PulseNoise035}
\end{subfigure}
\begin{subfigure}{.49\textwidth}
  \centering
 		\def\svgwidth{\columnwidth}
    		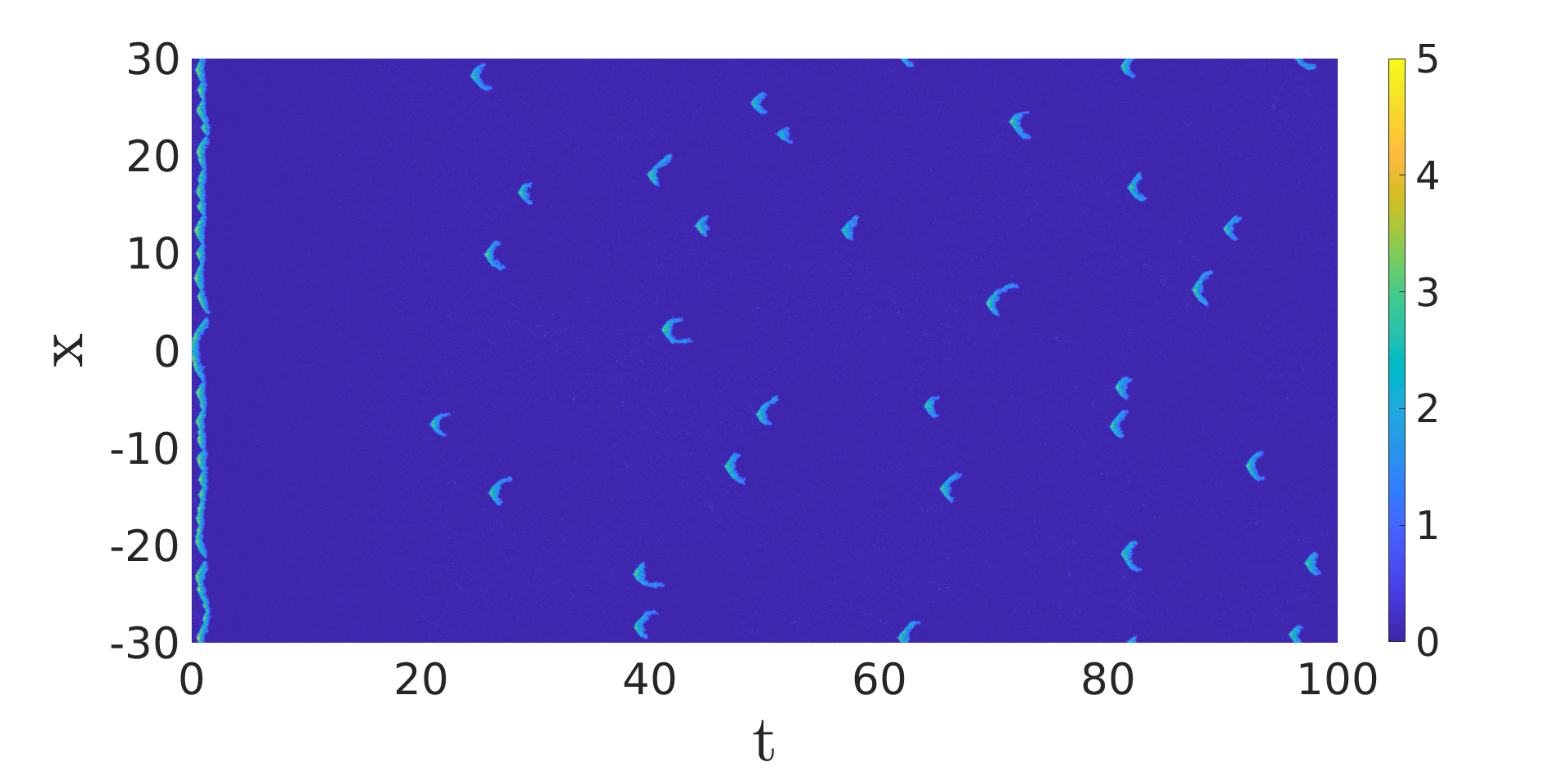
  \caption{$\sigma=0.05$}
    \label{fig:Sim:PulseNoise05}
\end{subfigure}
\begin{subfigure}{.49\textwidth}
  \centering
 		\def\svgwidth{\columnwidth}
    		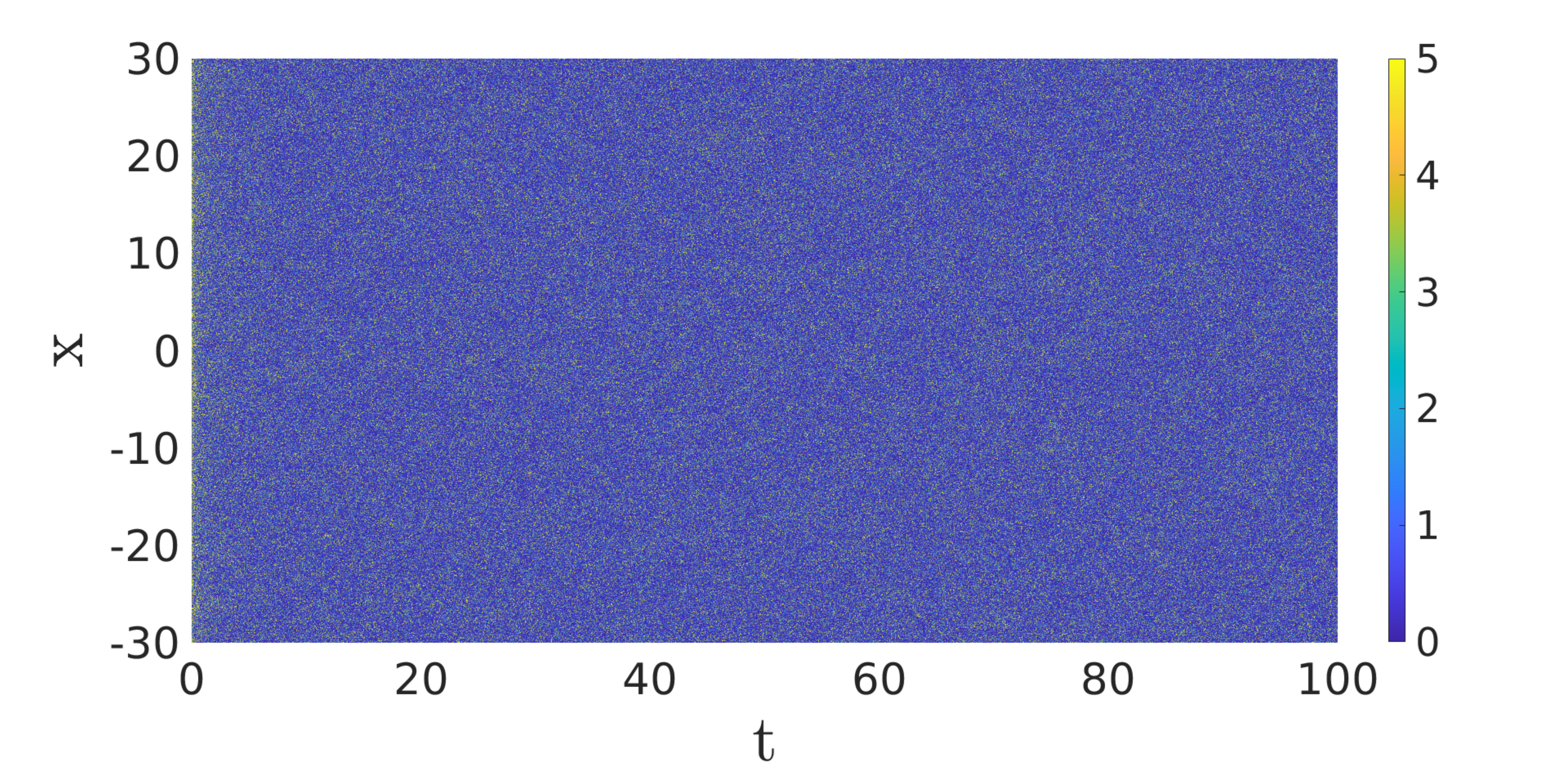
  \caption{$\sigma=0.5$}
    \label{fig:Sim:PulseNoise5}
\end{subfigure}
\caption{The $u$-component of the SPDE~\sref{eq:FullSPDE} for four different values of the noise $\sigma$. The other system parameters and initial conditions are the same as in the previous figures. In Figure~(a), we only show the simulation of wave integrated up to $T=20$ because the solution remains in the background state afterwards, the other three figures are shown up to $T=100$.}
\label{fig:Sim:PulseNoise}
\end{figure}

\paragraph{Standing waves with noise}
For the same parameter values as in the previous paragraph, we now study the full SPDE~\sref{eq:FullSPDE}. In Figure~\ref{fig:Sim:PulseNoise}, we plot realisations of the SPDE for different noise intensities. For low noise levels, we see two quasi-stationary waves appear, like in Figure~\ref{fig:Sim:Dpulse}, before they are destroyed at different points in time by the noise. Since the noise is low, no new activation events happen. When we increase the noise intensity, the noise is able to activate the stable background state, but the waves are also destroyed more quickly, resulting in a constant appearance and disappearance of waves. Note the comparison between Figures~\ref{fig:Sim:PulseNoise05} and the figures in~\cite{biswas2021enhanced}, where a similar model is studied using Gillespie algorithms. This activation of the background state is not possible in the deterministic PDE~\sref{eq:int:MainPDE} without an external force. In Figures~\ref{fig:Sim:PulseNoise035} and \ref{fig:Sim:PulseNoise05}, we see that in the first instances, many patterns are generated, causing the inhibitor to increase everywhere which blocks new activation events. After this initial phase, new activation events appear, and significantly more for higher values of the noise as expected. When we increase the noise even further, it becomes impossible to form patterns as every activation event is destroyed instantly. Therefore, pattern formation happens at intermediate values of the noise. The idea that there is some `optimal' value of the noise resulting in complex dynamics has been observed before in, for instance, the context of nerve impulses~\cite{garcia2001noise}. 

\begin{figure}
\begin{subfigure}{.5\textwidth}
  \centering
 		\def\svgwidth{\columnwidth}
    		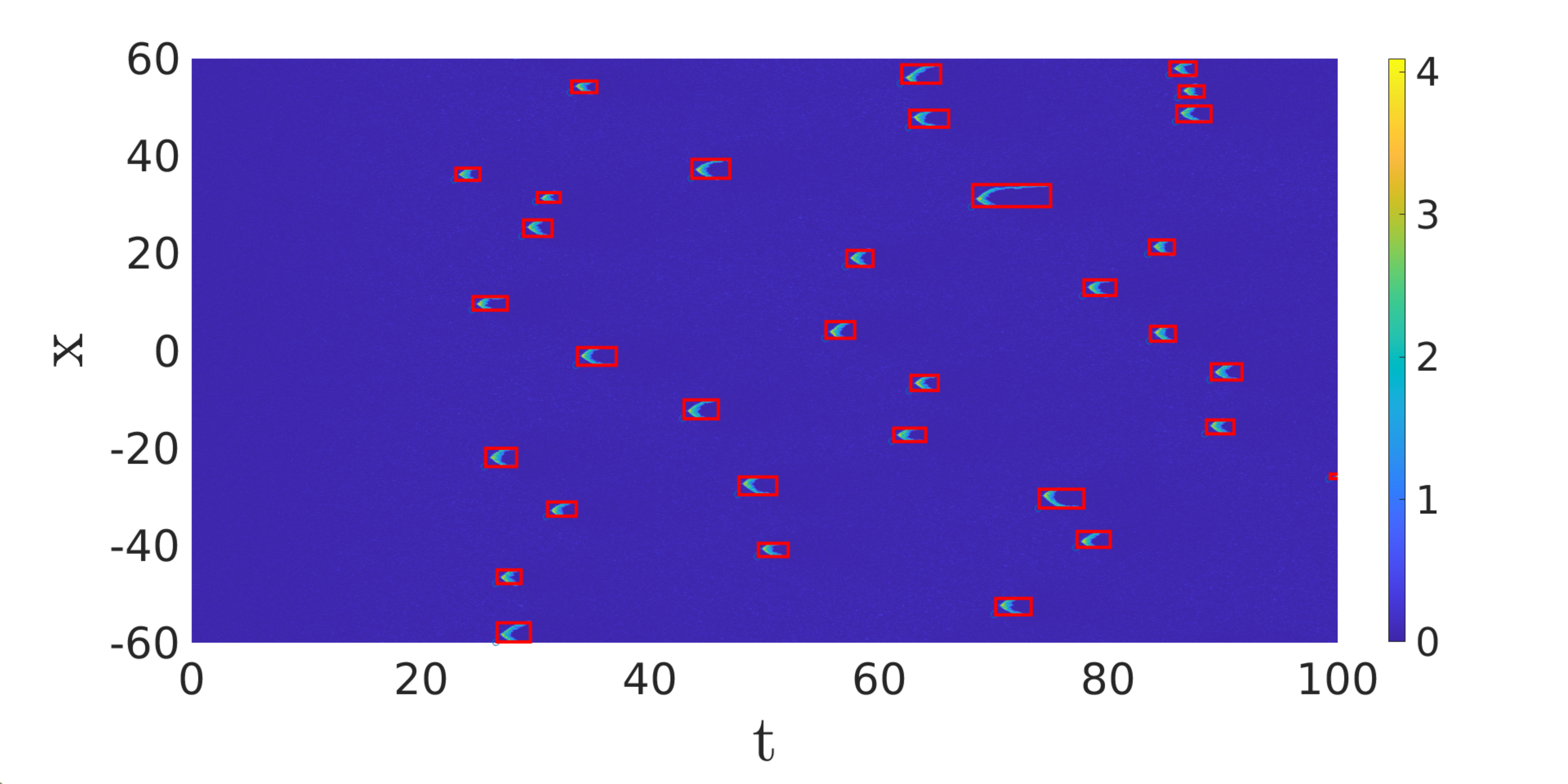
    \caption{}
    \label{fig:Sim:PulseNoiseBox}
\end{subfigure}
\begin{subfigure}{.5\textwidth}
  \centering
 		\def\svgwidth{\columnwidth}
    		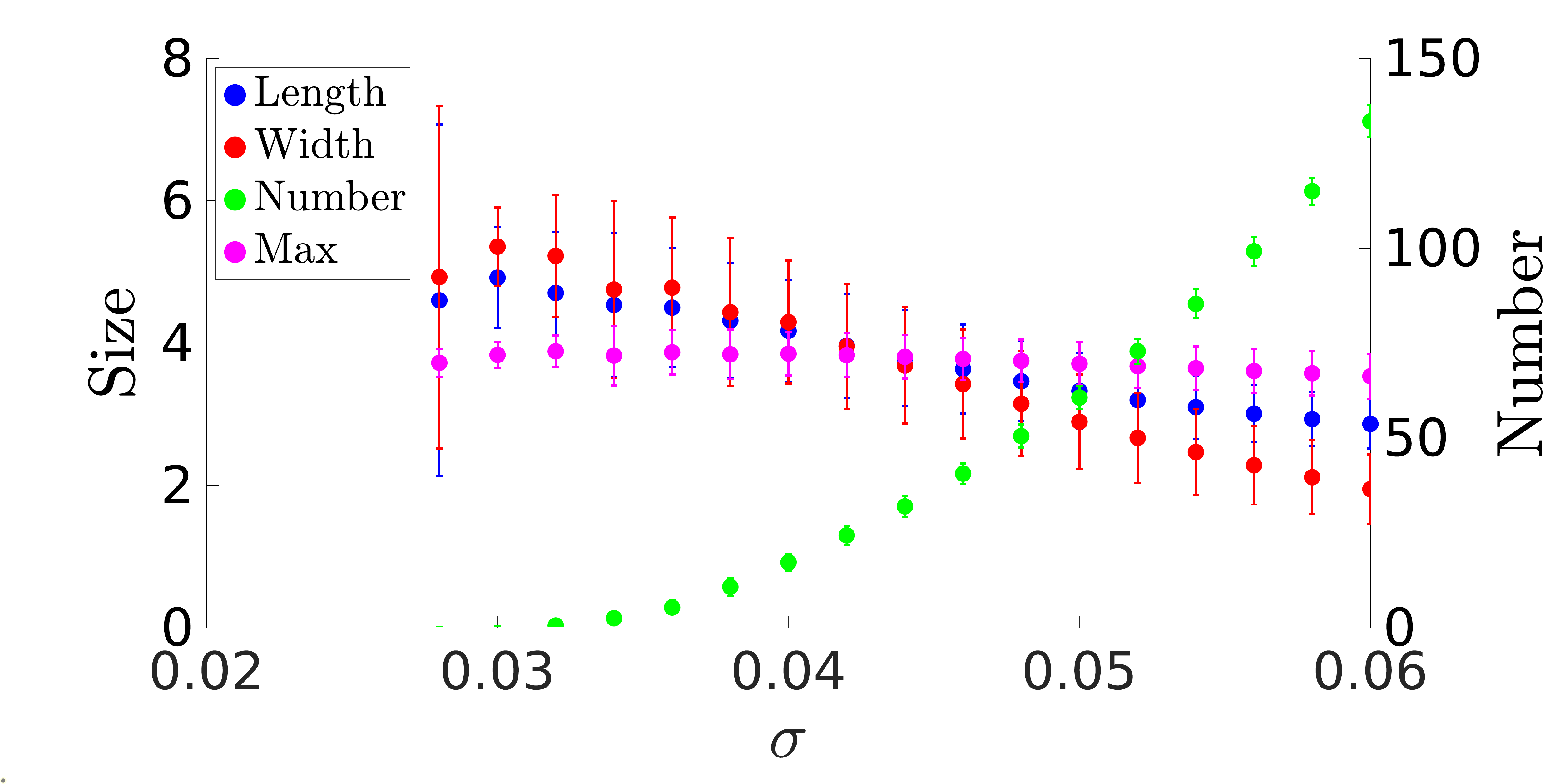
  \caption{}
    \label{fig:Sim:Stat}
\end{subfigure}
\caption{In Figure~(a) we show a simulation similar to those in Figure~\ref{fig:Sim:PulseNoise}, but with $\sigma=0.45$ and a homogeneous initial condition with $(u^*,4v^*)$. The red boxes are the result of the pattern finding algorithm \texttt{regionprops} in Matlab; it identifies all the regions of excitations which we would also find by eye, see Appendix~\ref{sec:app:num} for details. In Figure~(b), we used this algorithm to find the length, width and maximum of these pulses (left axis), as well as the total number of activation events (right axis). For each value of $\sigma$, the number of events is averaged over 100 simulations, and the length, width and maximum are averaged over all events in the 100 simulations. We plot the average together with the standard deviation.}
\label{fig:Sim:Box}
\end{figure}

\begin{figure}
\begin{subfigure}{.33\textwidth}
  \centering
 		\def\svgwidth{\columnwidth}
    		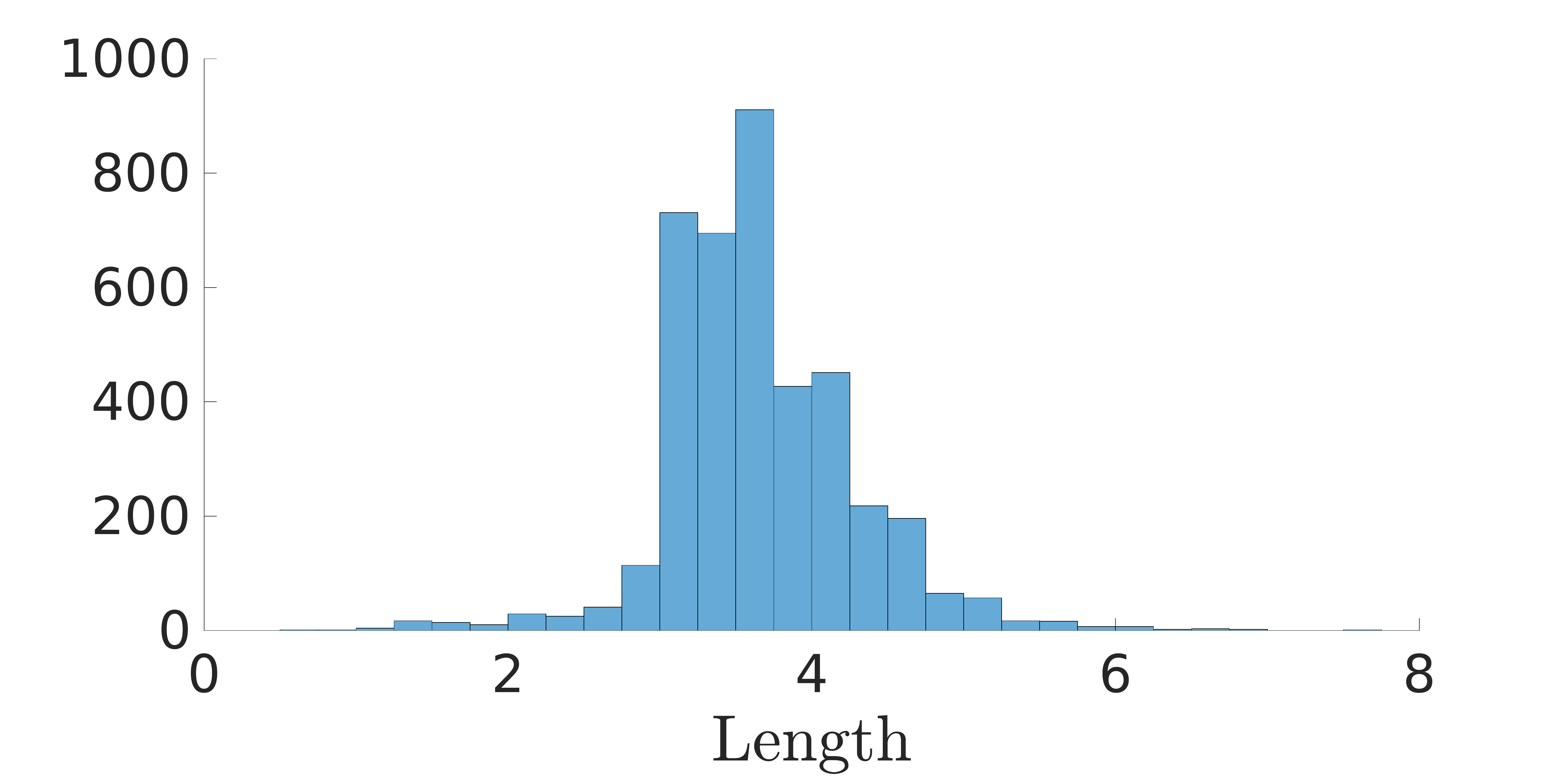
  \caption{}
    \label{fig:Sim:HistL}
\end{subfigure}
\begin{subfigure}{.33\textwidth}
  \centering
 		\def\svgwidth{\columnwidth}
    		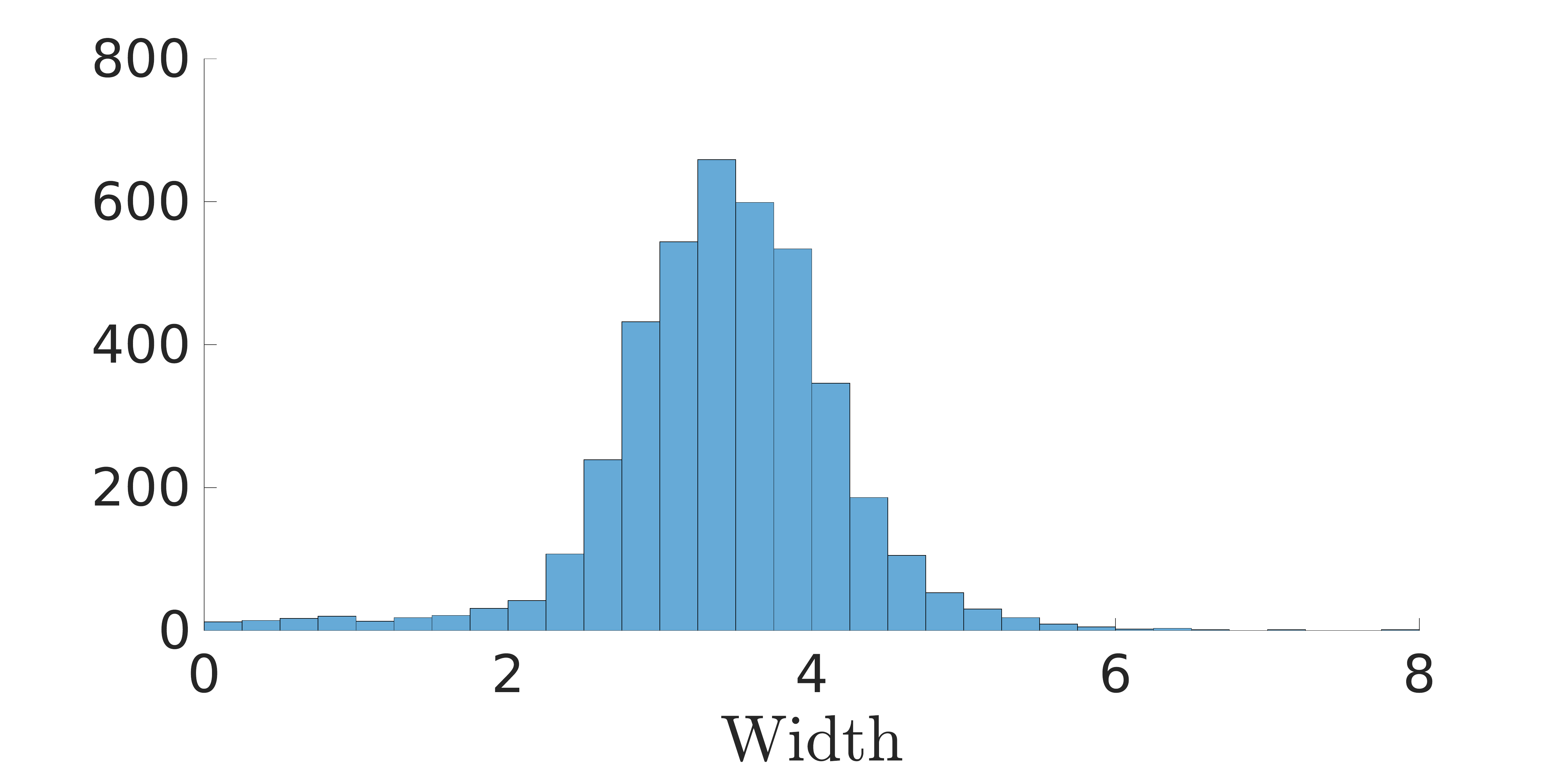
  \caption{}
    \label{fig:Sim:HistW}
\end{subfigure}
\begin{subfigure}{.33\textwidth}
  \centering
 		\def\svgwidth{\columnwidth}
    		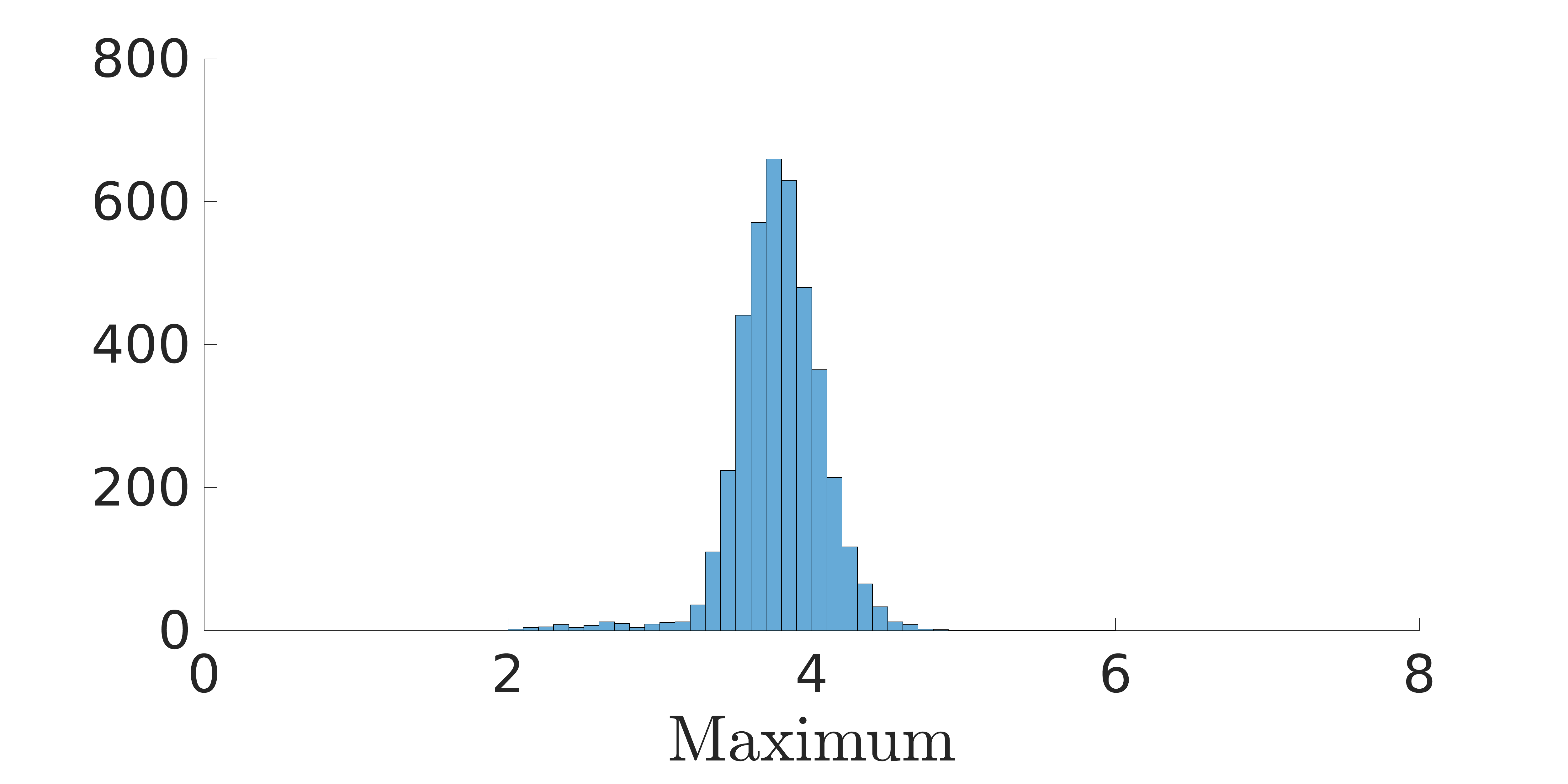
  \caption{}
    \label{fig:Sim:HistMax}
\end{subfigure}
\begin{subfigure}{.33\textwidth}
  \centering
 		\def\svgwidth{\columnwidth}
    		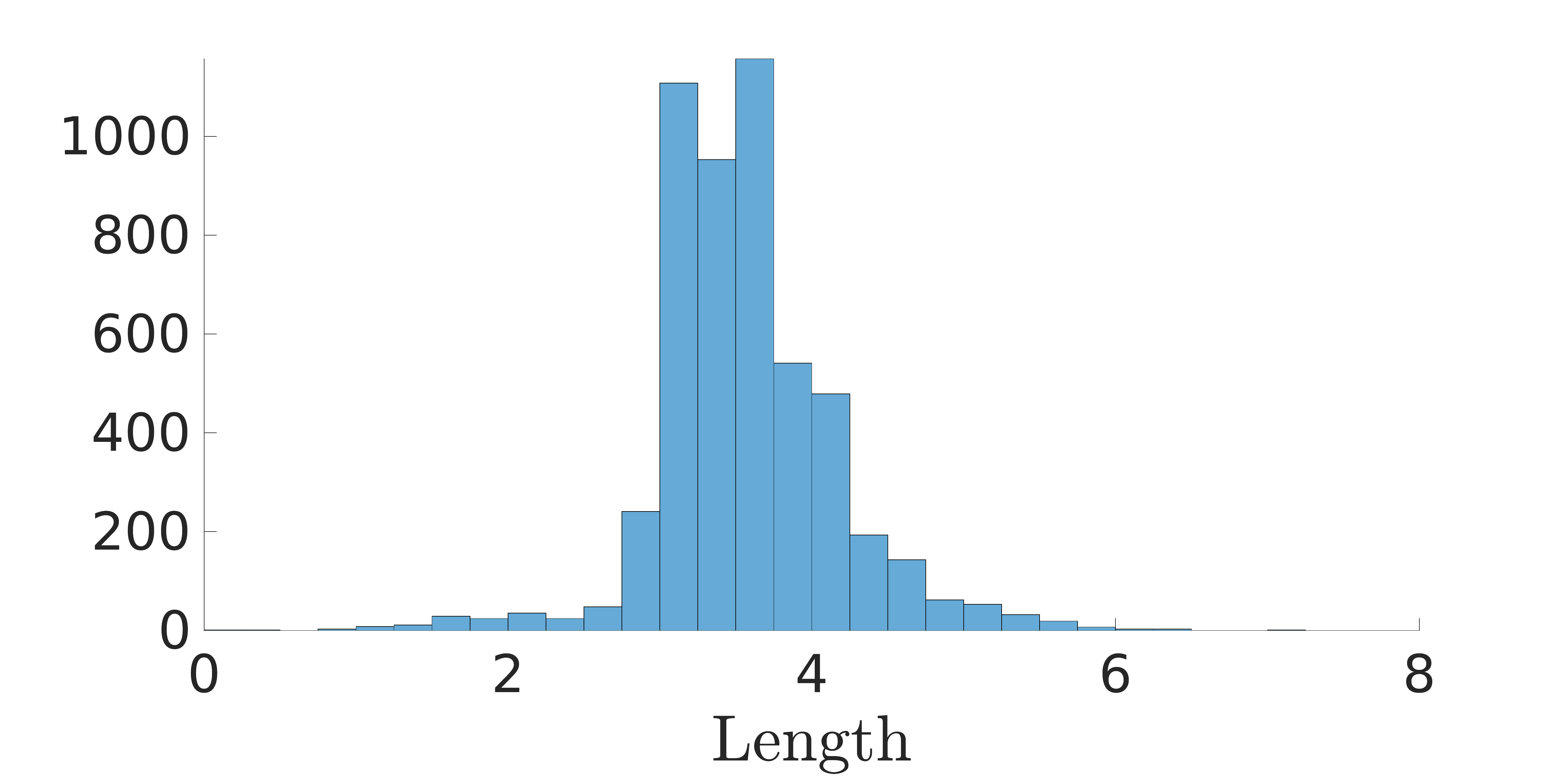
  \caption{}
    \label{fig:Sim:HistLNoDiff}
\end{subfigure}
\begin{subfigure}{.33\textwidth}
  \centering
 		\def\svgwidth{\columnwidth}
    		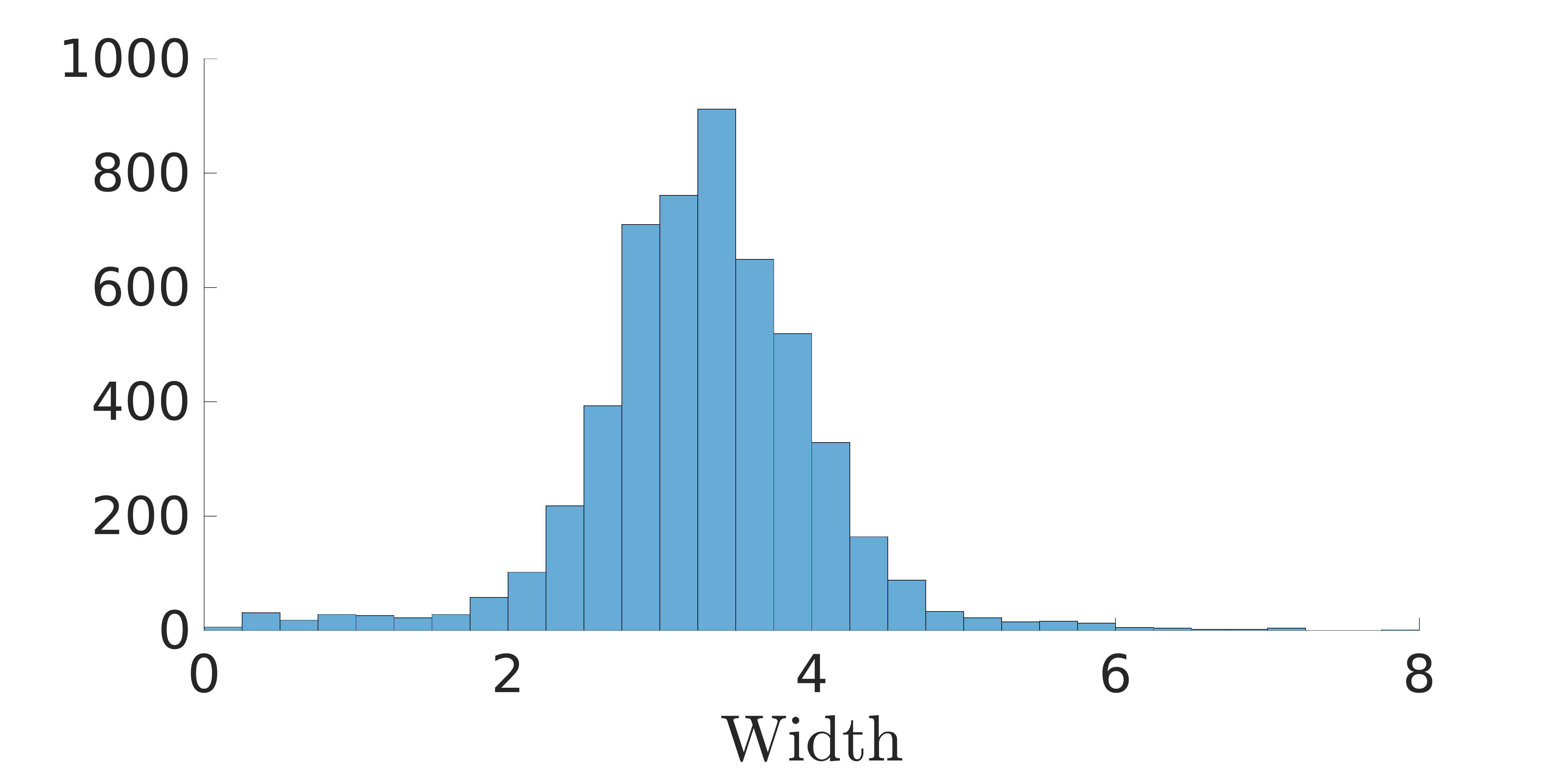
  \caption{}
    \label{fig:Sim:HistWNoDiff}
\end{subfigure}
\begin{subfigure}{.33\textwidth}
  \centering
 		\def\svgwidth{\columnwidth}
    		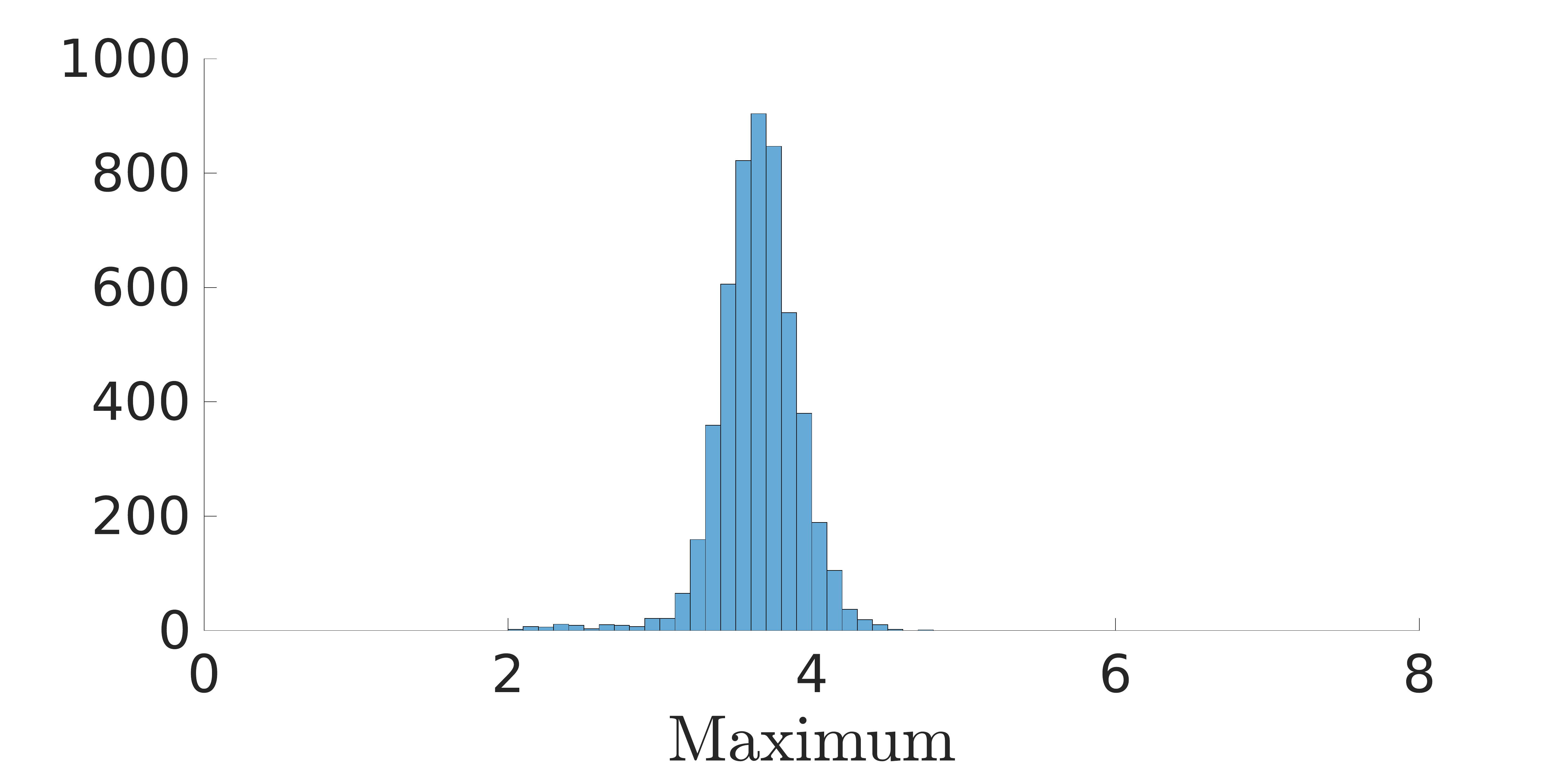
  \caption{}
    \label{fig:Sim:HistMaxNoDiff}
\end{subfigure}
\begin{subfigure}{.33\textwidth}
  \centering
 		\def\svgwidth{\columnwidth}
    		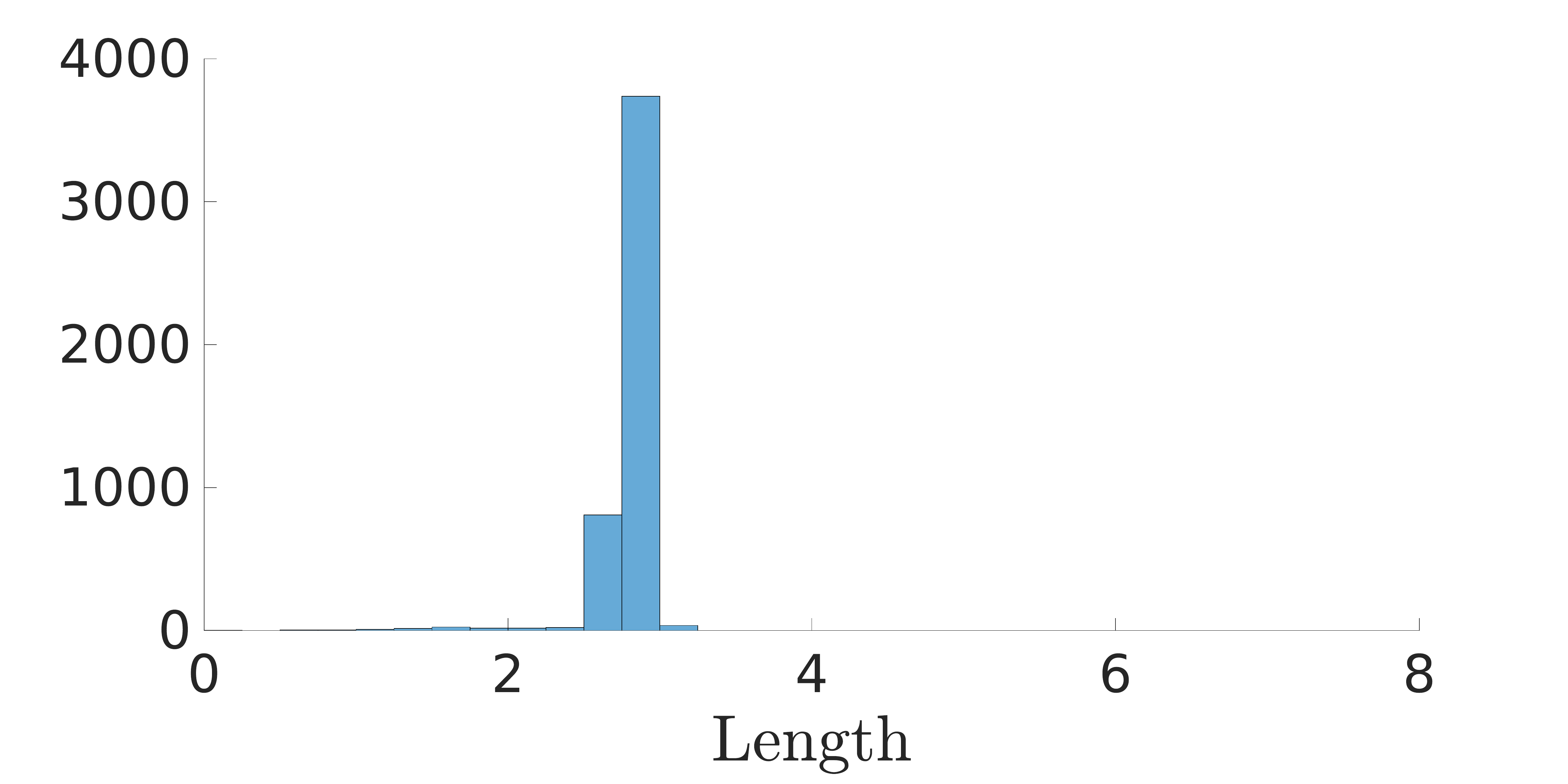
  \caption{}
    \label{fig:Sim:HistLWhite}
\end{subfigure}
\begin{subfigure}{.33\textwidth}
  \centering
 		\def\svgwidth{\columnwidth}
    		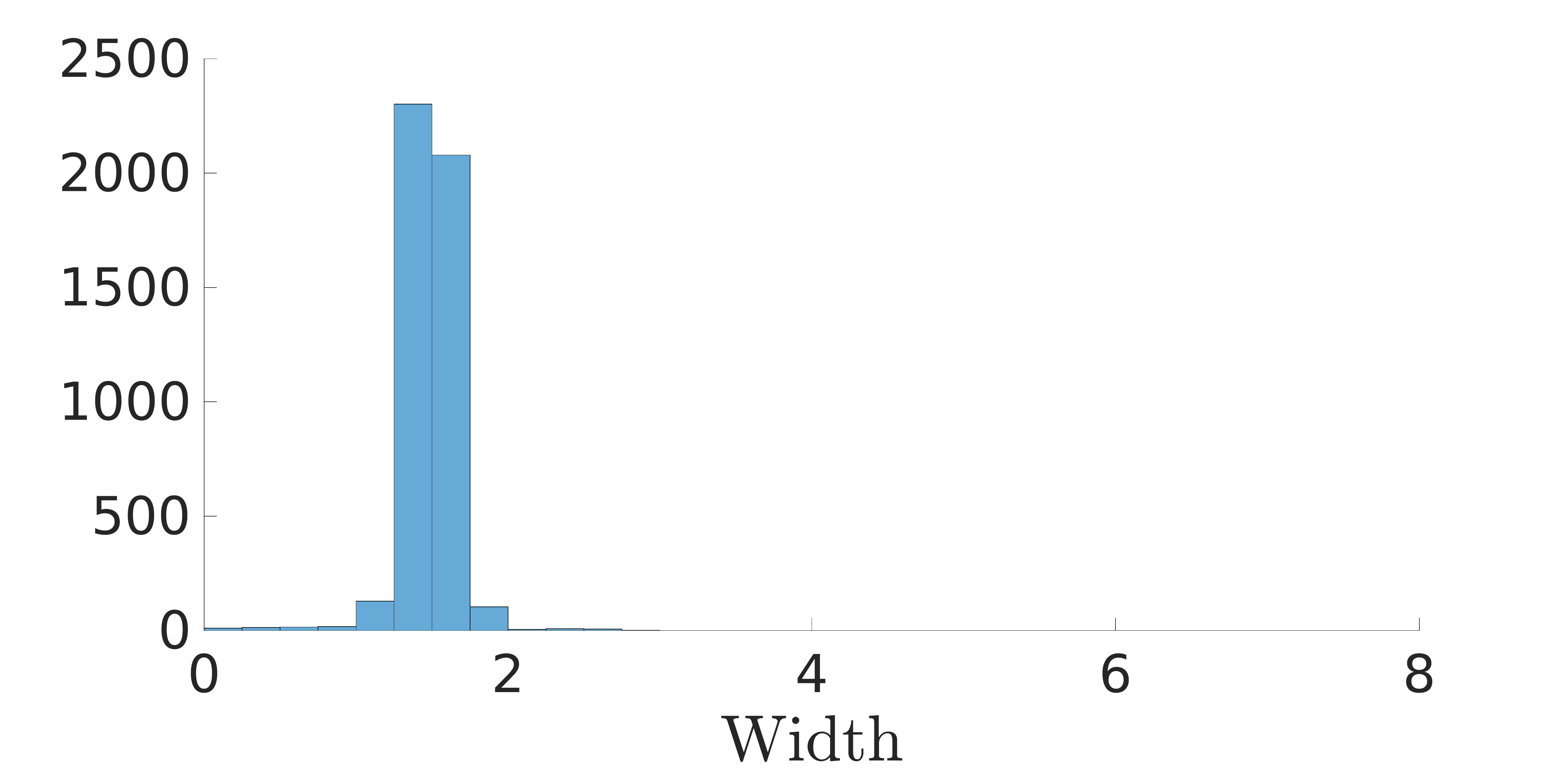
  \caption{}
    \label{fig:Sim:HistWWhite}
\end{subfigure}
\begin{subfigure}{.33\textwidth}
  \centering
 		\def\svgwidth{\columnwidth}
    		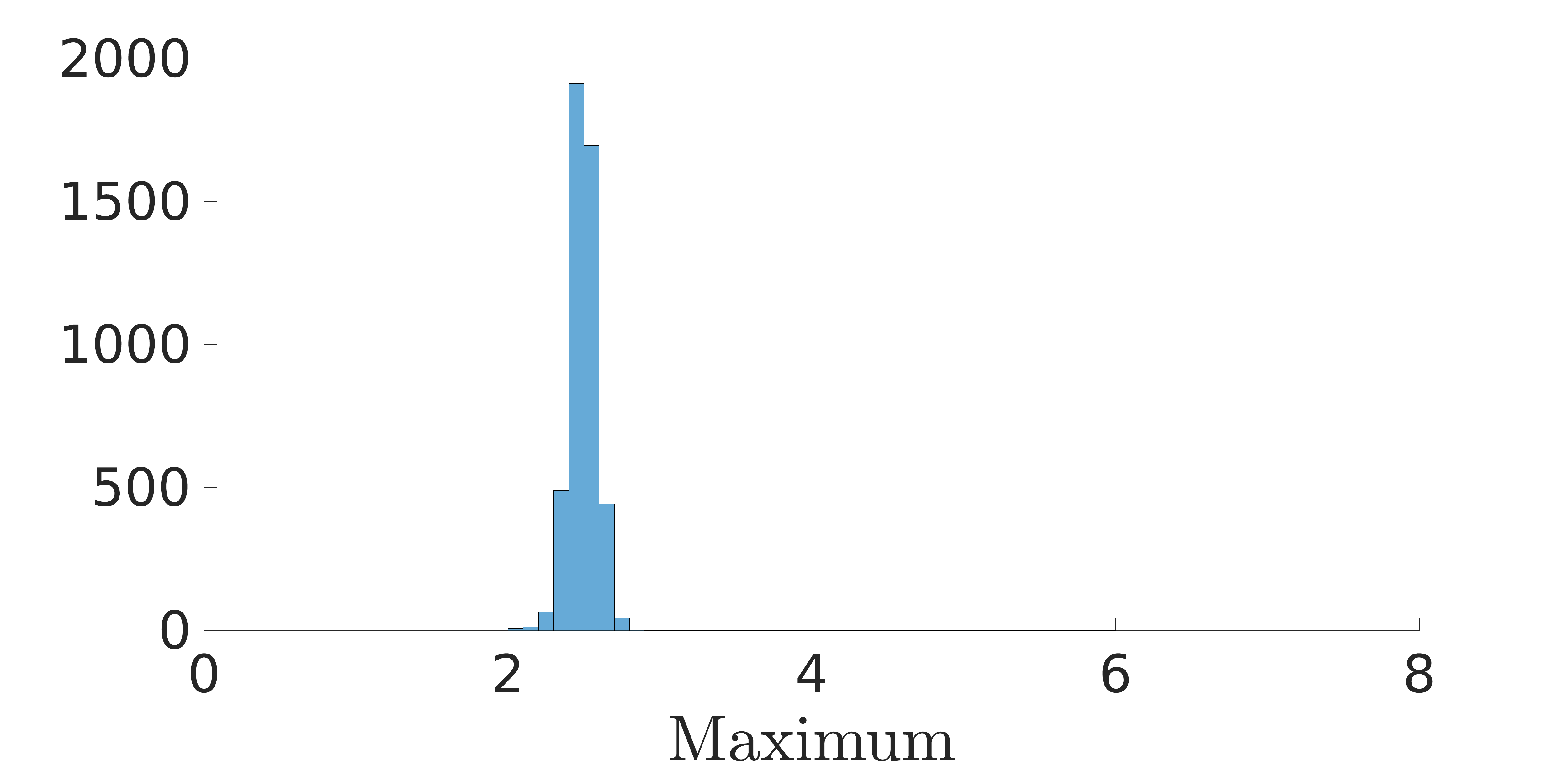
  \caption{}
    \label{fig:Sim:HistMaxWhite}
\end{subfigure}
\caption{Figures~(a)-(c) show the histograms for the width, length and maximum of the pulses for $\sigma=0.046$, for the same data as in Figure~\ref{fig:Sim:Stat}. For Figures~(a) and (b), the bin width is fixed to $0.25$, and for Figure~(c) to $0.1$. Figures (d)-(f) show the same histograms, but in the simulations, the noise coming from the diffusion (the last term in \sref{eq:FullSPDE}) was set to 0. In Figures (g)-(i), we again show the same histograms, but with just white noise on the $u$-component. In order to compare the noise levels, we did not choose the same $\sigma$ value for the three cases but chose $\sigma$ values such that the average number of activation events per simulation is approximately 50. For the Figures (d)-(f), this means $\sigma=0.056$, and for (g)-(i) $\sigma=0.3$.}
\label{fig:Sim:Hist}
\end{figure}

In order to quantify this notion of optimality in the noise intensity we must first quantify the size and shape of the patterns in Figures~\ref{fig:Sim:PulseNoise035} and \ref{fig:Sim:PulseNoise05}. Using Matlab's \texttt{regionprops} algorithm we can automatically detect the patches with a high value for the activator $u$ (see Appendix~\ref{sec:app:num} for details), giving us the possibility to compute the number of activation events and determine the width and duration of each event, see Figure~\ref{fig:Sim:PulseNoiseBox}. In Figure~\ref{fig:Sim:Stat} we show the statistics for a range of $\sigma$ values. This figure shows that there is a clear cutoff for when activation events are likely to happen. For values of $\sigma<0.035$, the average number of events is lower than $1$, and the number of activation events increases sharply after this value. We observe that the width, the length and the maximum height of the events are all higher when the number of excitation events is low, but the variability in these values is also larger. in Figure~\ref{fig:Sim:Hist}, we look at the statistics of the events for the specific value $\sigma=0.046$. The value of the maximum is sharply peaked. This is something we expect, as the maximum is mainly determined by the deterministic dynamics after the excitation. The width and length of the events are much more spread out. Especially for the width, we see a heavy tail towards zero. This is also expected because activation events come in two forms. Most events result in two waves, but a small part of the events has the shape of just a single wave, which has a width of 0.87 in the deterministic case. We checked whether or not these histograms are well approximated by a Gaussian distribution, but this was rejected using a Kolmogorov-Smirnov test $(p\sim 10^{-14})$.

Using the statistics on the width, length, and maximum, we can compare the solutions of SPDE~\sref{eq:FullSPDE} to SPDEs with the same deterministic part but different noise terms. First, we can set the $\partial_x\sqrt{2DX}dW_t$ term coming from the diffusion to zero. As noted in Remark \ref{rem:ill}, this term makes the mathematical analysis of the SPDE~\sref{eq:FullSPDE} significantly harder. Figures~\ref{fig:Sim:HistLNoDiff}-\ref{fig:Sim:HistMaxNoDiff} show that the statistics of the solutions do not change significantly when we delete this term. This indicates that the noise coming from the reaction terms plays a more influential role in determining the shape of the patterns. 

We are now also in the position to compare the CLE approach with the more ad hoc approach of adding additive white noise the to $u$-component to mimic the inherent noisiness of the system, see Figures~\ref{fig:Sim:HistLWhite}-\ref{fig:Sim:HistMaxWhite}. Indeed, the properties of the patterns are significantly different when we compare them to the full SPDE. In particular,  with just white noise, the patterns are all short and narrow and do not reflect the complicated dynamics of the underlying chemical reactions and experiments (not shown).

\subsection{Travelling Waves}
In order to find a travelling wave solution of~\sref{eq:int:MainPDE}, understood as a solution of~\sref{eq:sim:bvp} with $c\neq0$, we must ensure that the dynamics starting from the initial condition does not reach the standing phase or returns to the background state. This can be achieved by increasing the value of $c_1$. Increasing $c_1$ results in a faster exponential decay of $v$ back to the background state after an excitation, see Table \ref{tab:my_label}, preventing the inhibitor from glueing the two waves together like in Figure~\ref{fig:Sim:Dpulse}. Simulations for an increased value of $c_1$, from $0.1$ to $0.2$\footnote{For values in between, say $c_1=0.15$, the numerics becomes very sensitive to the chosen discretisation, see Appendix~\ref{sec:app:num}.}, are shown in Figure~\ref{fig:Sim:TWext}. Note that the PDE~\sref{eq:int:MainPDE} still only has one stable background state $(u^*,v^*) \approx (0.0833, 1.625)$. The initial condition splits into two counter-propagating travelling waves, but opposite to what happened with the standing wave before, they keep separating and move away from each other at a fixed speed until they collide and cancel each other out due to the periodicity of the domain, see Figure~\ref{fig:Sim:TWext}.

\begin{figure}
\begin{subfigure}{.49\textwidth}
  \centering
 		\def\svgwidth{\columnwidth}
    		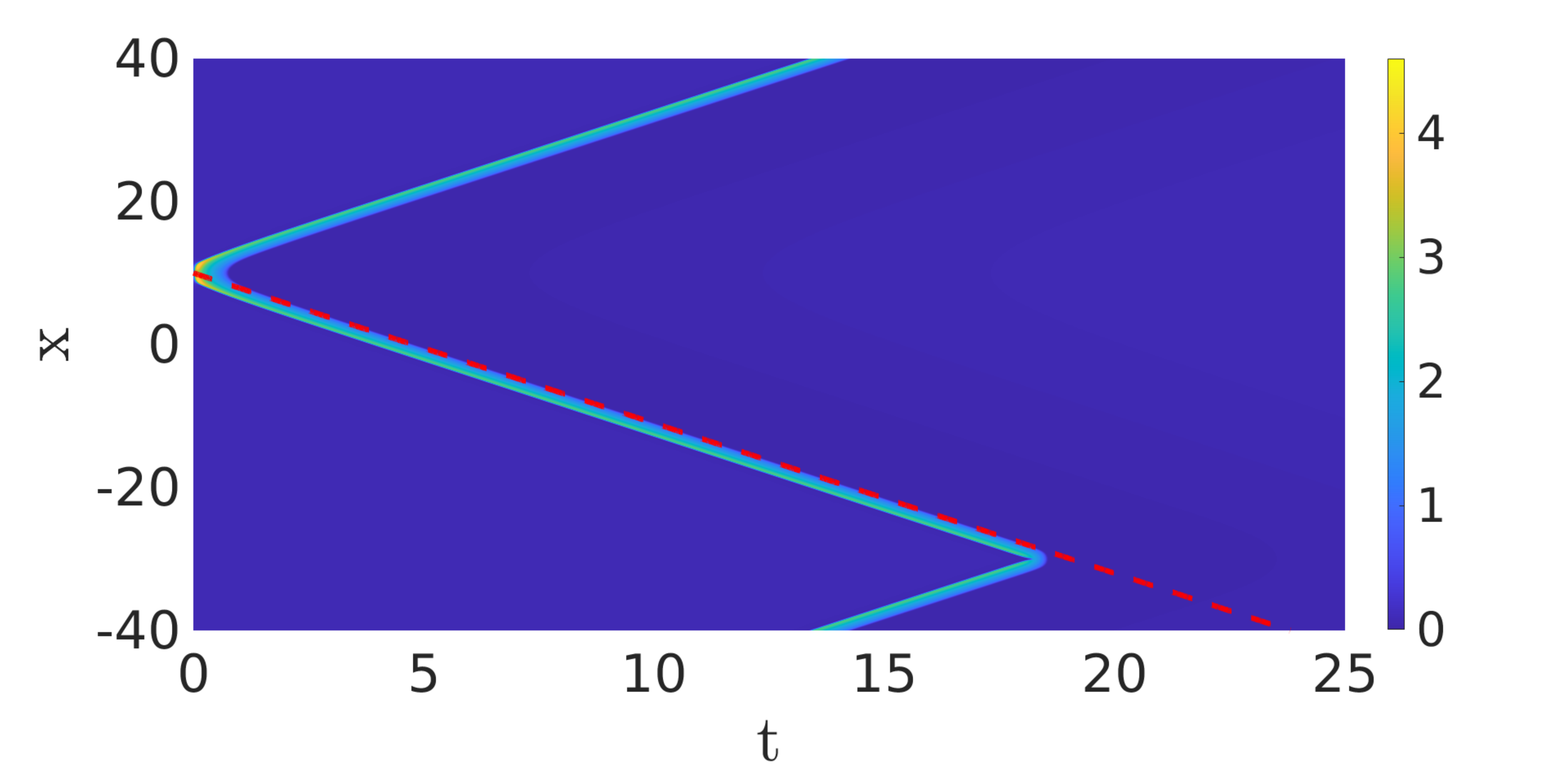
  \caption{}
    \label{fig:Sim:TWuExt}
\end{subfigure}
\begin{subfigure}{.49\textwidth}
  \centering
 		\def\svgwidth{\columnwidth}
    		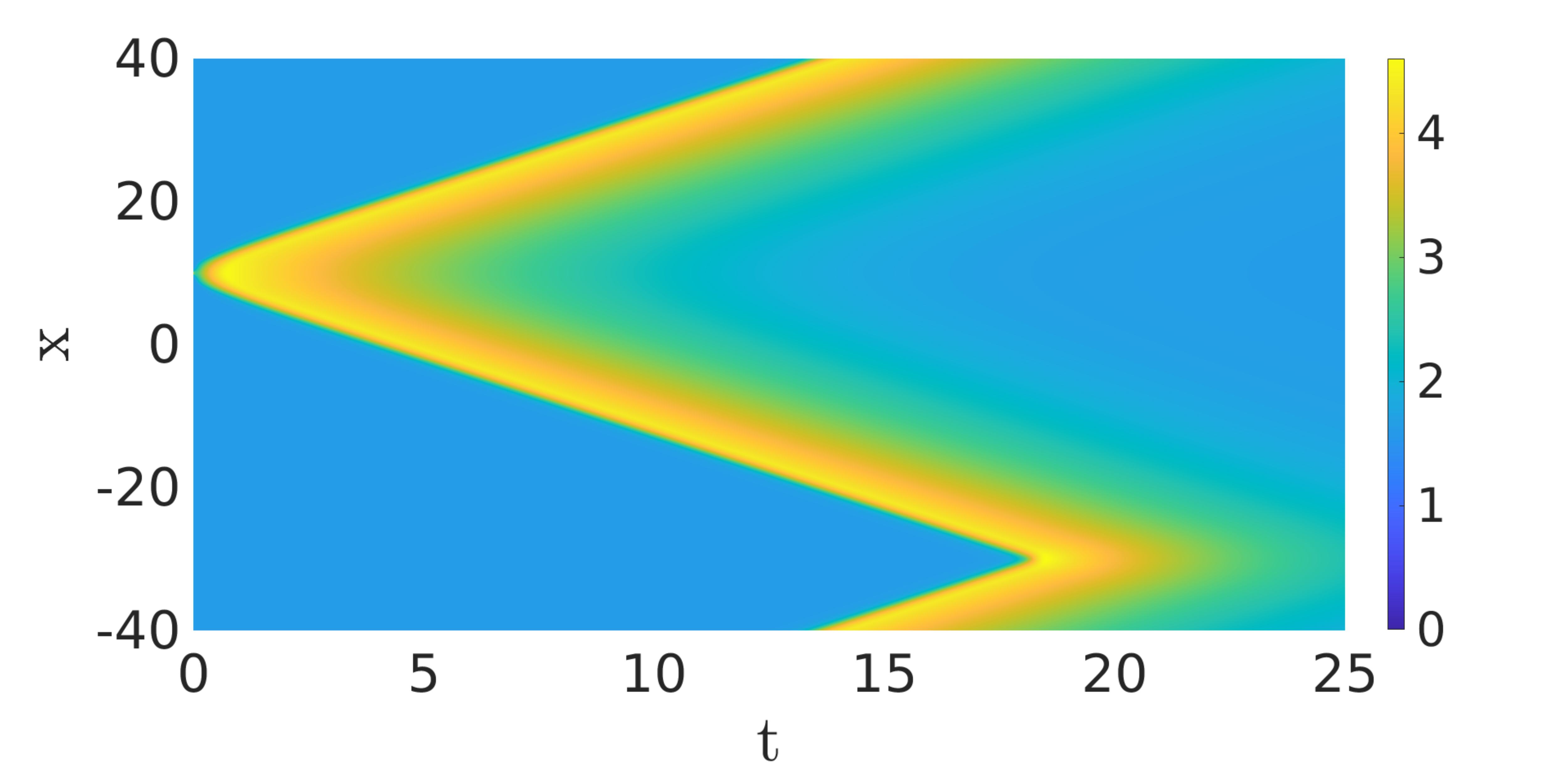
  \caption{}
    \label{fig:Sim:TWvExt}
\end{subfigure}
\caption{Simulation of the PDE~\sref{eq:int:MainPDE}, Figure~(a) shows the activator $u$ and Figure~(b) the inhibitor $v$. We observe the splitting of the initial condition in two counterpropagating travelling waves with a constant speed that exist until they cancel each other out due to the periodicity of the domain. The slow inhibitor $v$ decays back to its rest state in between the pulses. The red dotted line has a speed of $-2.10$, which is close to the value of approximately $-2.17$ found by solving~\sref{eq:sim:bvp} using a fixed point method. The parameters are $D_u=0.1$, $a_1=0.167$, $a_2=16.67$, $a_3=167$, $a_4=1.44$, $a_5=1.47$, $D_v=1$, $\e=0.52$,  $c_1=0.2$ and $c_2=3.9$.}
\label{fig:Sim:TWext}
\end{figure}

\begin{figure}
\begin{subfigure}{.49\textwidth}
  \centering
 		\def\svgwidth{\columnwidth}
    		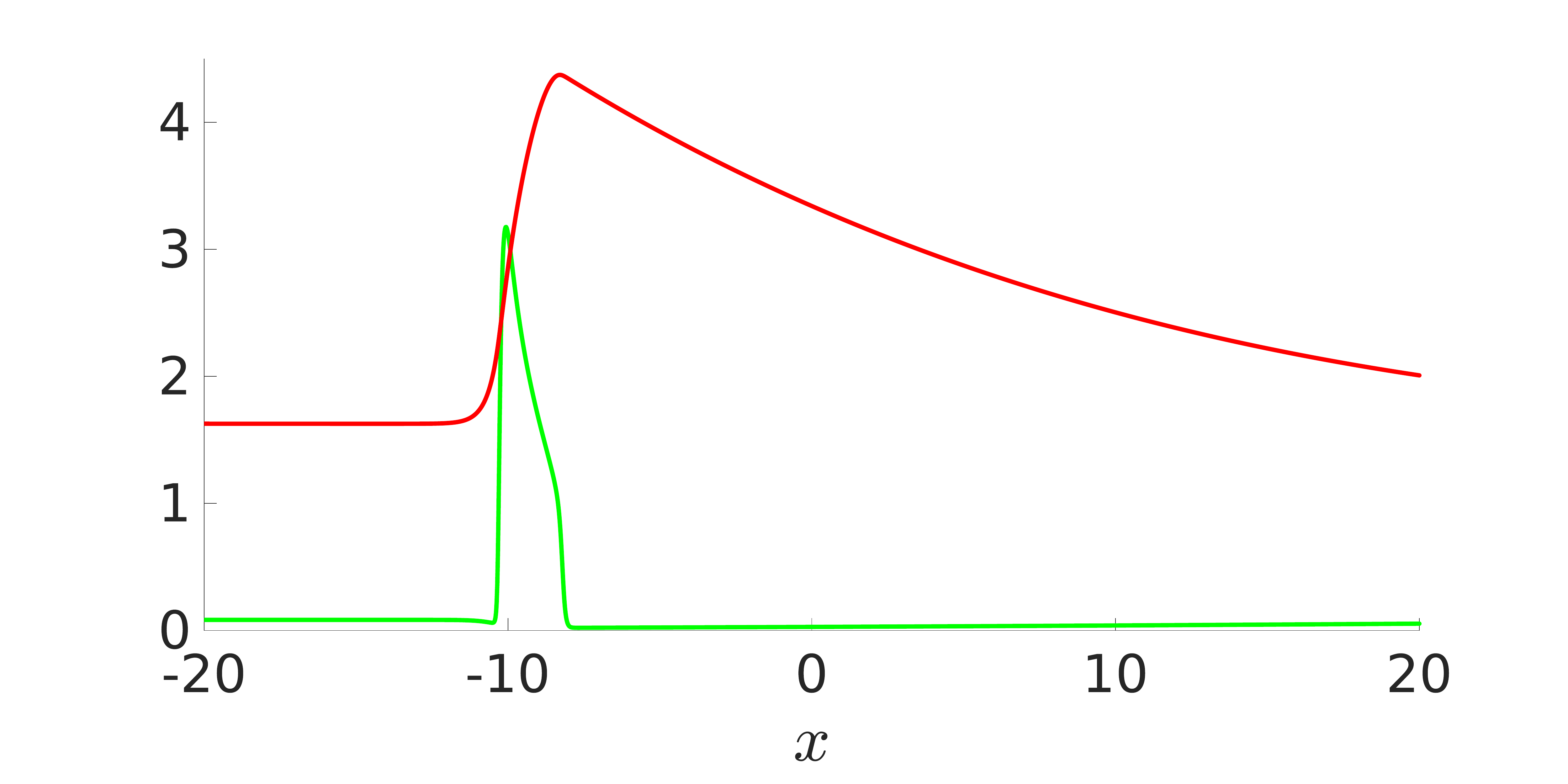
  \caption{}
    \label{fig:Sim:TWprofile}
\end{subfigure}
\begin{subfigure}{.49\textwidth}
  \centering
 		\def\svgwidth{\columnwidth}
    		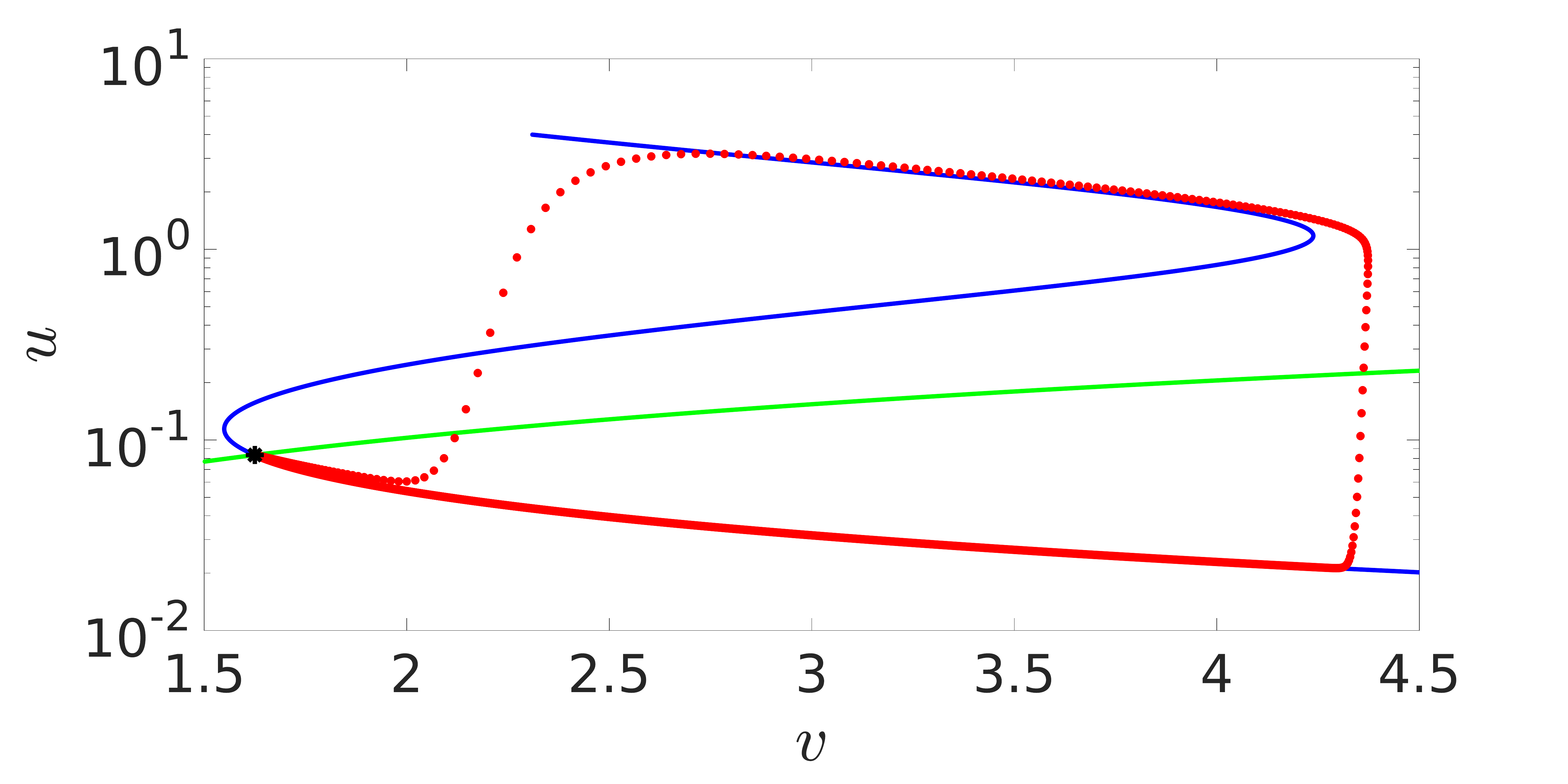
  \caption{}
    \label{fig:Sim:TWPP}
\end{subfigure}
\caption{Profile of a single travelling wave.  Figure~(a) shows both components $u$ (green) and $v$ (red) and Figure~(b) the related phase plane, plotted on a semi-log scale to highlight the dynamics for small $u$, as well as the nullclines. The asterisk indicates the fixed point. This solution is obtained as the endpoint of a PDE simulation (not shown), i.e. similar to Figure~\ref{fig:Sim:TWext}, but with just one of the two waves as initial condition.}
\label{fig:Sim:TW2}
\end{figure}

To find a single travelling wave, we again need to properly tune the initial condition. This can be done by selecting one of the two waves in Figure~\ref{fig:Sim:TWext} and using it as the initial condition of the PDE simulation (not shown). In Figure~\ref{fig:Sim:TW2} we show the travelling wave profile and its associated phase plane. As with the standing pulse, the dynamics around the $u$-nullcline is essential. The solution trajectory starts from near the background state and follows the lower branch of the $u$-nullcline, jumps towards the upper branch of the nullcline and keeps following it until it falls off and returns to the lower branch to slowly evolve back towards the stable background state. In contrast to the standing pulse, see Figure~\ref{fig:sim:PulseBoth}, the travelling wave is no longer symmetric and it jumps back to the lower branch by falling off the edge of the upper branch. These travelling wave solutions could be analysed further using techniques similar to Appendix~\ref{sec:analysis}. 

It is important to realise that we do not expect to see travelling waves in practice as the travelling wave gets destroyed when it collides with another wave. Therefore, in the stochastic simulations, it might not always be clear if we are looking at a travelling wave that collapses or at the transient dynamics towards a double pulse that subsequently gets destroyed by the noise.

\begin{figure}
\begin{subfigure}{.49\textwidth}
  \centering
 		\def\svgwidth{\columnwidth}
    		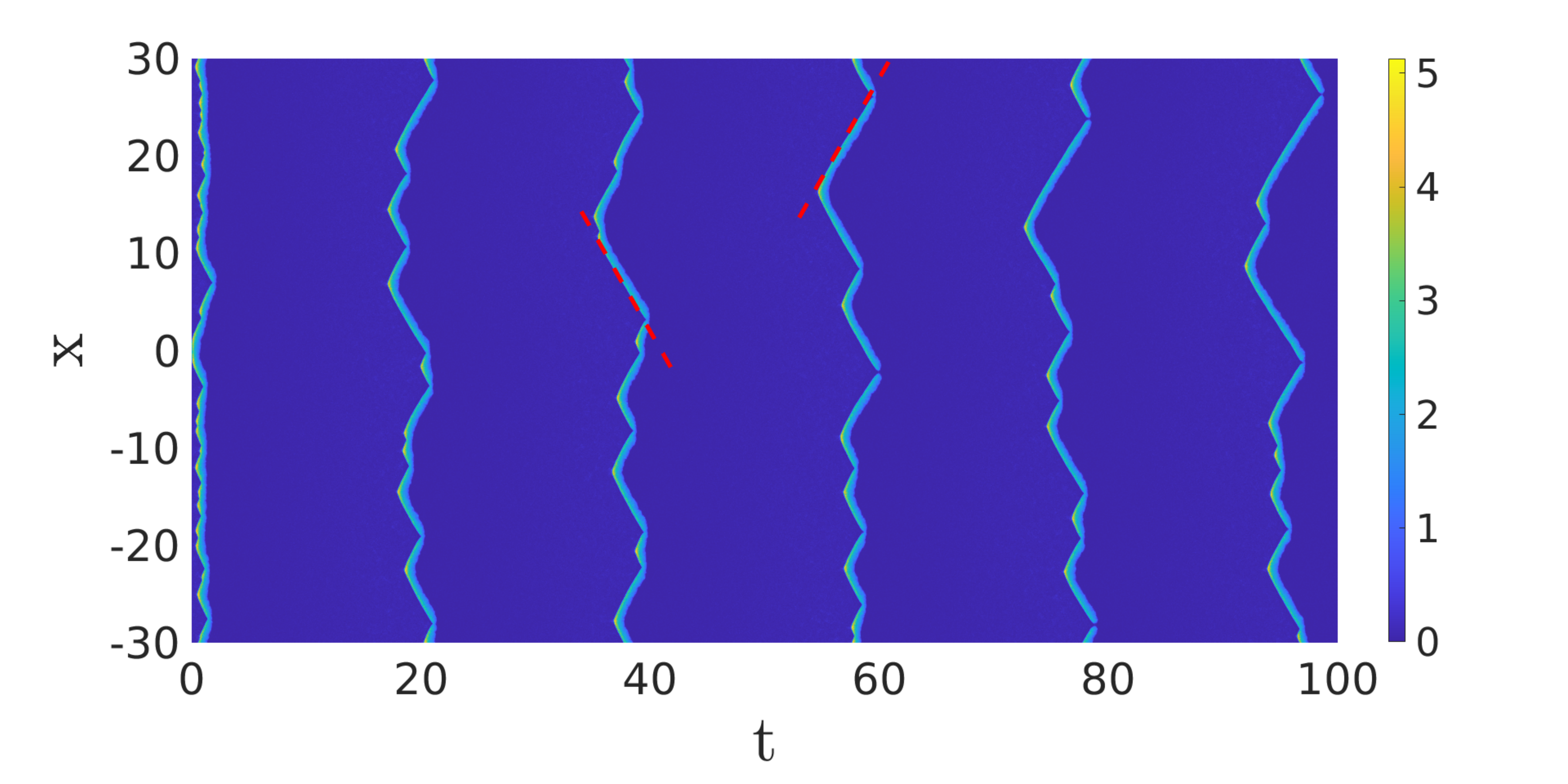
  \caption{}
    \label{fig:Sim:TWnoiseU}
\end{subfigure}
\begin{subfigure}{.49\textwidth}
  \centering
 		\def\svgwidth{\columnwidth}
    		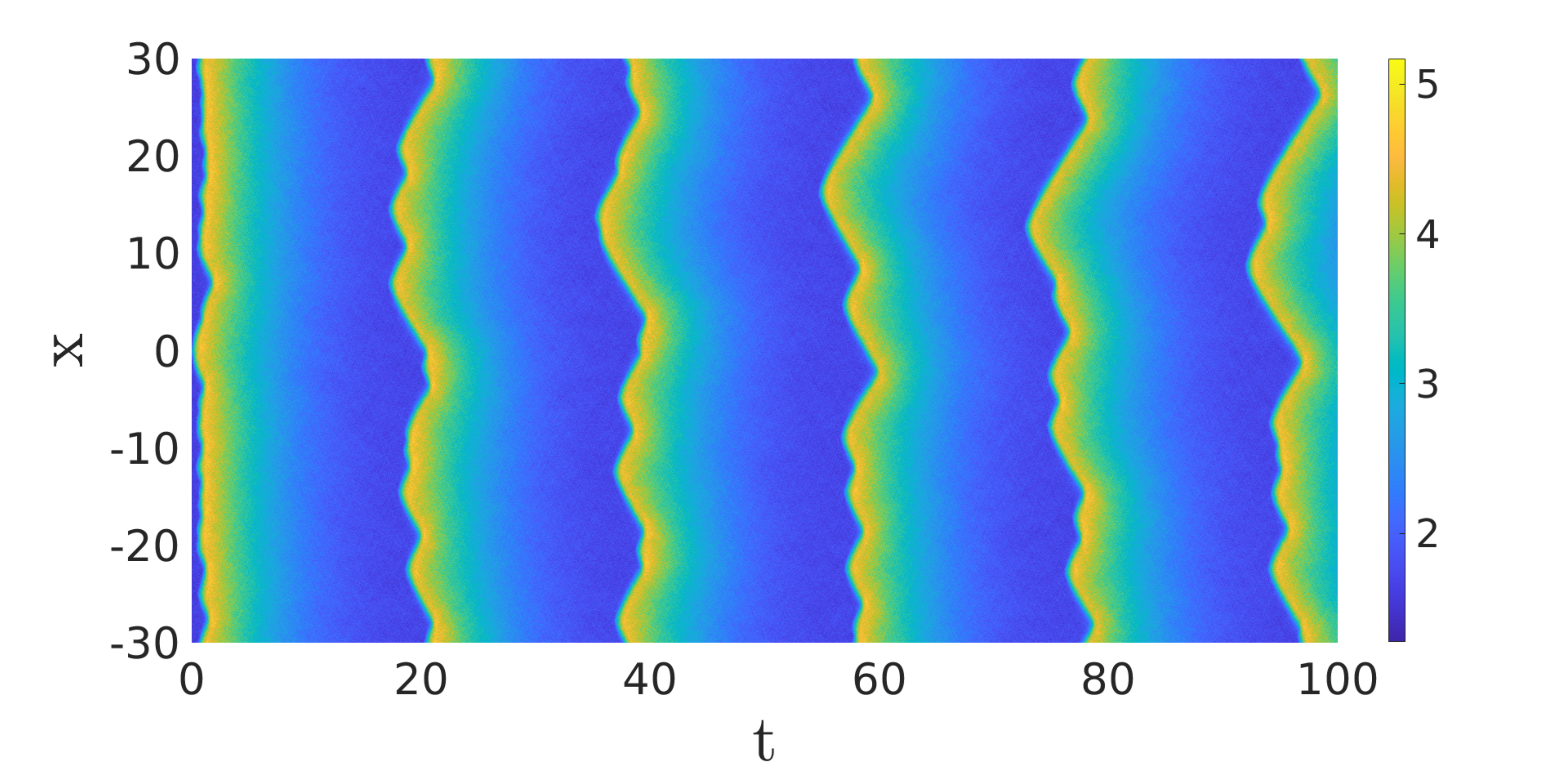
  \caption{}
    \label{fig:Sim:TWnoiseV}
\end{subfigure}
\begin{subfigure}{.49\textwidth}
  \centering
 		\def\svgwidth{\columnwidth}
    		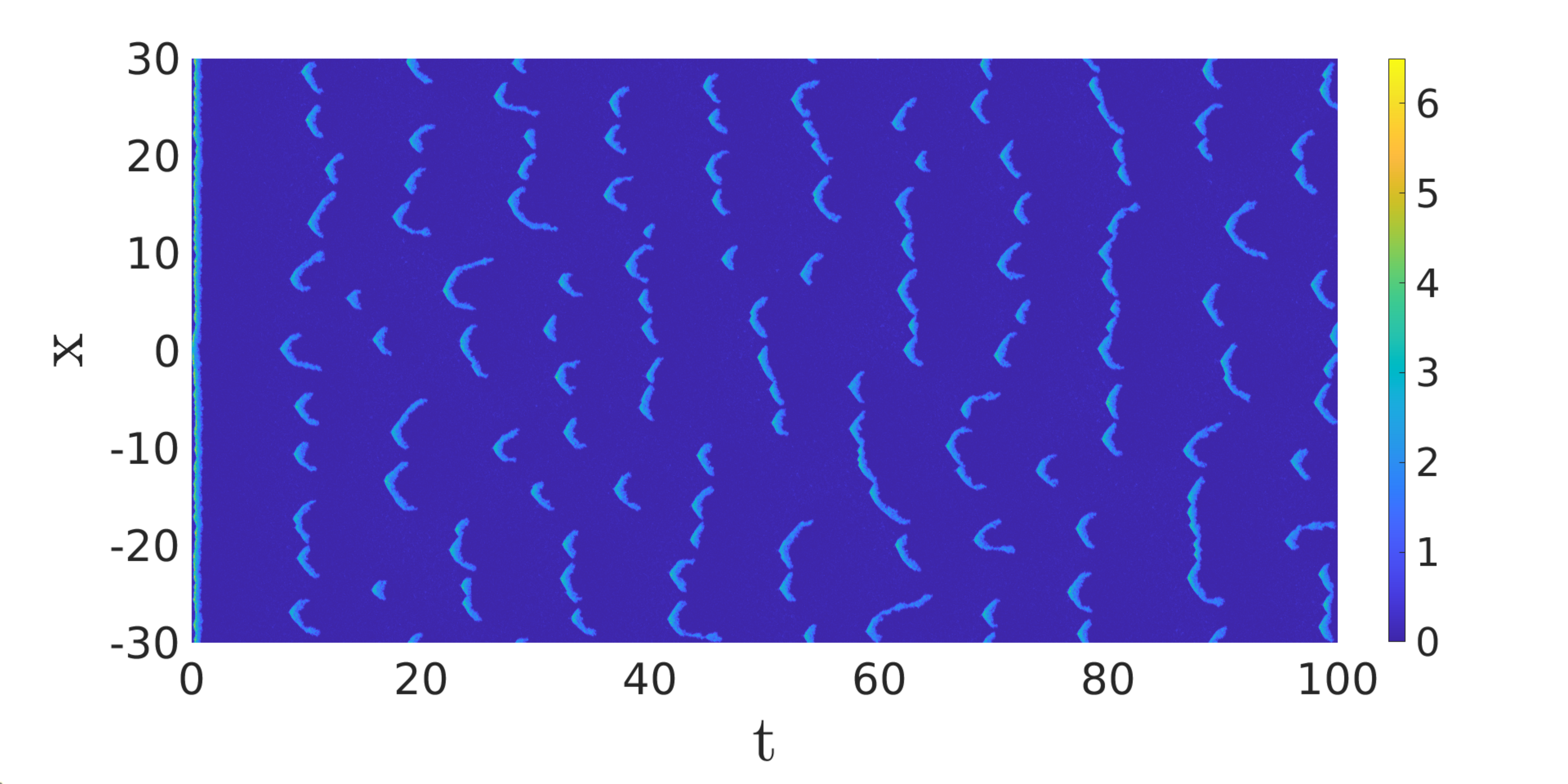
  \caption{}
    \label{fig:Sim:TWnoiseU05}
\end{subfigure}
\begin{subfigure}{.49\textwidth}
  \centering
 		\def\svgwidth{\columnwidth}
    		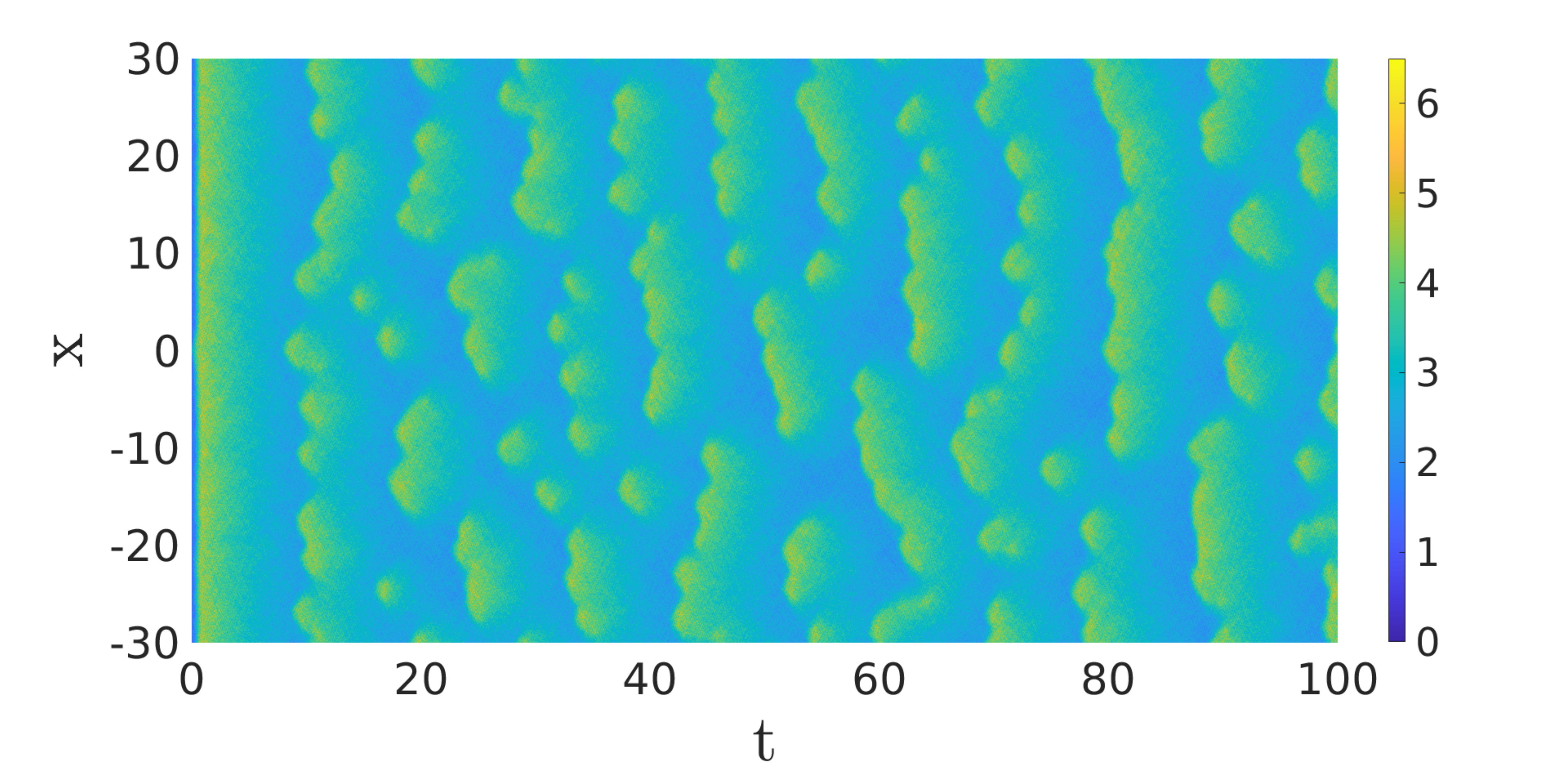
  \caption{}
    \label{fig:Sim:TWnoiseV05}
\end{subfigure}
\caption{Simulation of the SPDE~\sref{eq:FullSPDE} for $\sigma=0.02$, Figures~(a) and (b), and $\sigma=0.05$, Figures~(c) and (d). The red dashed line in (a) has a slope of 2.05, close to the deterministic wave speed, but given the short time interval the wave exists, precise estimates are difficult to obtain. We observe that there is a quasi-periodic behaviour with a period of roughly 20. In Figures~(c) and (d), the quasi-periodic structure is destroyed. The same parameters and initial condition are used as in Figure~\ref{fig:Sim:TWext}.}
\label{fig:Sim:TWnoise}
\end{figure}

\paragraph{Travelling waves with noise}
When we now return to SPDE~\sref{eq:FullSPDE}, there are now four regimes for the same  parameters as in the previous section. For high values of the noise, we, as before, do not observe any patterns (not shown). For low values of the noise, we just find the travelling wave (if the simulation is initiated by an appropriate initial condition) since the noise is not strong enough to destroy the wave, nor to activate another pattern, on the timescales of the simulation (not shown). 
The interesting dynamics happens again at the intermediate levels of the noise. As Figure~\ref{fig:Sim:TWnoiseU} shows, the noise activates the dynamics, resulting in many counter-propagating travelling waves. A travelling wave is subsequently annihilated when it collides with a travelling wave coming from the other direction. Hence, the collision dynamics of Figure~\ref{fig:Sim:TWext} is repeated many times on smaller spatial-temporal scales. We see in Figure~\ref{fig:Sim:TWnoise} that after the annihilation of the travelling waves, the slow inhibitor $v$ initially remains high preventing the activation of new counterpropagating travelling waves. Only when after a certain time the inhibitor has sufficiently decayed, do we see the activation of new counterpropagating travelling waves by the noise. The creation and annihilation of travelling waves happen at a shorter time scale than the decay of the inhibitor, which makes the dynamics look synchronised, or even periodic. In Figure~\ref{fig:QPeriodvSigma} we plot the approximate period versus the intensity of the noise. As expected, the period decreases with the intensity of the noise.  It differs however significantly from the true time periodic motion we will discuss in $\S$\ref{sec:per}. When we increase the noise, the quasi-periodic pattern is broken up, as the counter-propagating travelling waves are destroyed before they collide and annihilate each other, so no synchronised patterns emerge, see Figures~\ref{fig:Sim:TWnoiseU05}and \ref{fig:Sim:TWnoiseV05}. These patterns become relevant when we discuss the comparison between the CLE and the Gillespie simulations in Figure~\ref{fig:FromScienceSim}, see $\S$\ref{sec:WvsP}

\subsection{Time Periodic Solutions}
\label{sec:per}
\begin{figure}
\begin{subfigure}{.49\textwidth}
  \centering
 		\def\svgwidth{\columnwidth}
    		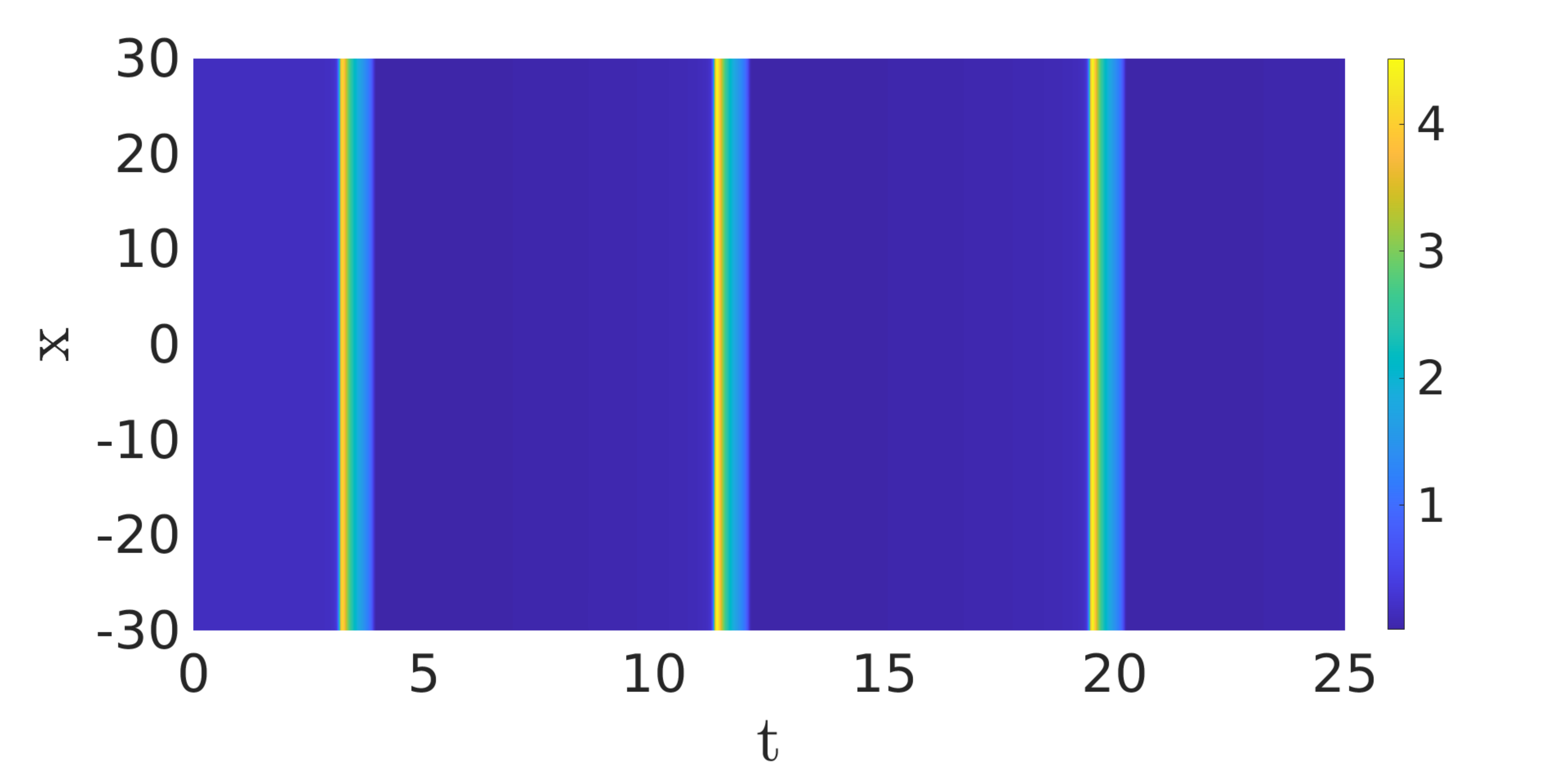
  \caption{}
    \label{fig:Sim:Periodicu}
\end{subfigure}
\begin{subfigure}{.49\textwidth}
  \centering
 		\def\svgwidth{\columnwidth}
    		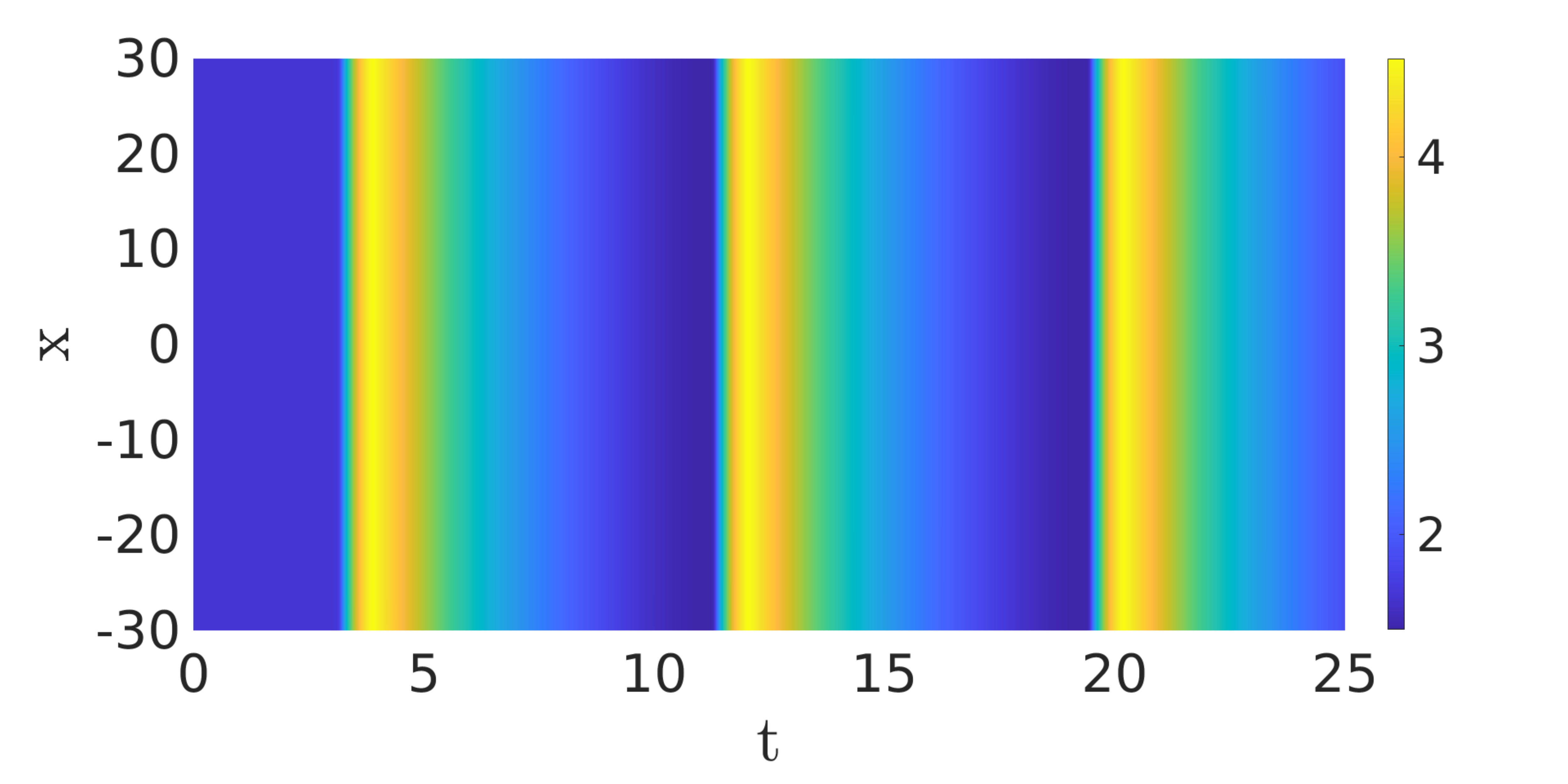
  \caption{}
    \label{fig:Sim:Periodicv}
\end{subfigure}
\caption{Simulation of the PDE~\sref{eq:int:MainPDE}, Figure~(a) shows the activator $u$ and Figure~(b) the inhibitor $v$. By measuring the distances between the maxima of the oscillations we find the estimate $T=8.14$ for the period of the oscillation. Note that this is significantly smaller than the quasi-periodic oscillations in Figure~\ref{fig:Sim:TWnoise}. The parameters are set to $D_u=0.1$, $a_1=0.167$, $a_2=16.67$, $a_3=167$, $a_4=1.44$, $a_5=1.47$, $D_v=1$, $\e=0.52$, $c_1=0.4$ and $c_2=3.9$.}
\label{fig:Sim:Periodic1}
\end{figure}

\begin{figure}
\begin{subfigure}{.49\textwidth}
  \centering
 		\def\svgwidth{\columnwidth}
    		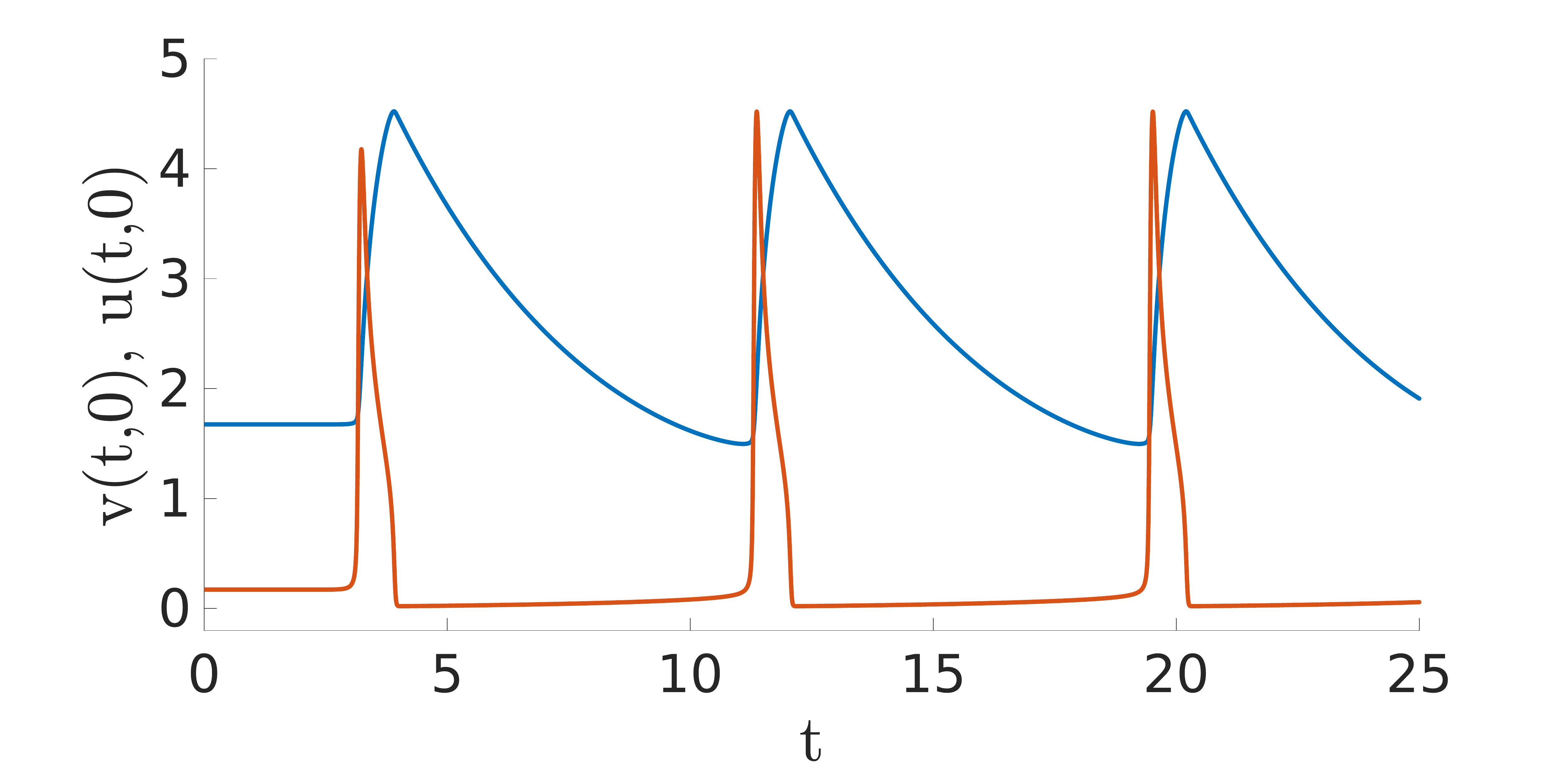
  \caption{}
    \label{fig:Sim:Periodic}
\end{subfigure}
\begin{subfigure}{.49\textwidth}
  \centering
 		\def\svgwidth{\columnwidth}
    		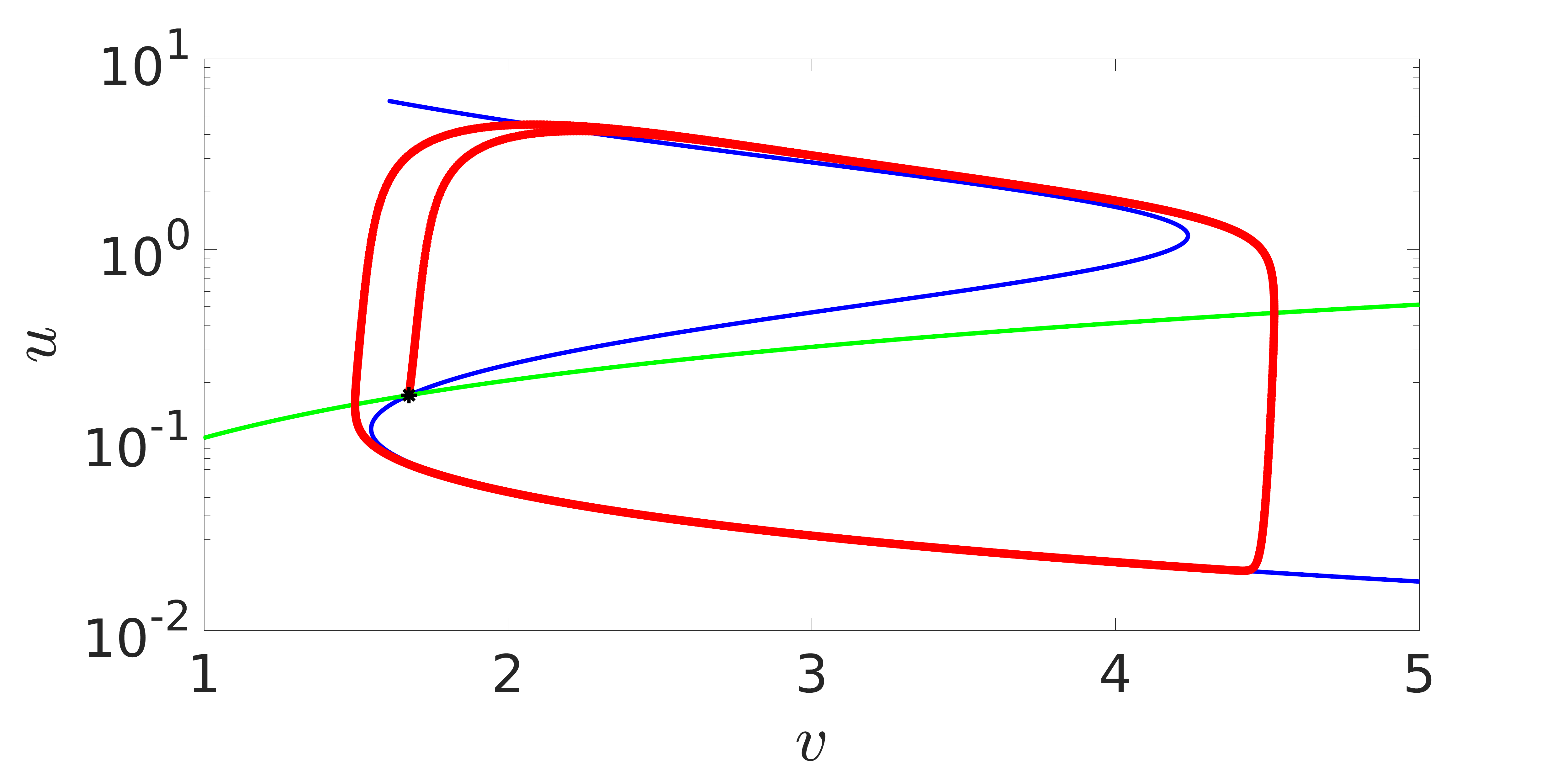
  \caption{}
    \label{fig:Sim:PeriodicPP}
\end{subfigure}
\caption{Cross-section of Figure~\ref{fig:Sim:Periodic1} at $x=0$, together with the corresponding phase plane. It is clear that the solution leaves the background state (marked by an asterisk), but does not return to it.}
\label{fig:Sim:PeriodicB}
\end{figure}

In the previous sections, it was essential that the background state of the system was stable, because this allowed the dynamics to return to the rest state after an activation event. When we increase the value of $c_1$, the background state becomes unstable through a Hopf bifurcation, see Figure~\ref{fig:bifdiag}. 
In the phase plane, this transition is characterised by the fact that the background state is no longer located on the lower branch of the $u$-nullcline, as in Figures~\ref{fig:Sim:PulsePP} and \ref{fig:Sim:TWPP}, instead, it lies on the middle branch of the $u$-nullcline, see Figure~\ref{fig:Sim:PeriodicPP}. Hence, after an excursion, the solution cannot return to the unstable background state and is exited again, resulting in time-periodic motion. When we start with a spatial homogeneous initial condition, the PDE simulation shows periodic oscillations in time, see Figures~\ref{fig:Sim:Periodic1} and \ref{fig:Sim:PeriodicB}. Both components still display slow-fast behaviour, however, this time not in the spatial variable $x$ but in the temporal variable $t$. In the case of nonhomogeneous initial conditions, it takes several oscillations before they are all synced up spatially (not shown). The observed behaviour has the characteristics of a relaxation oscillation as studied intensively for the Van der Pol equation~\cite{van1926}. This is not a surprise as the Van der Pol equation formed the foundation for the classic FitzHugh-Nagumo model and PDE~\sref{eq:int:MainPDE} can be seen as a variation on this classic model.

\paragraph{Time periodic solutions with noise}
For small values of the noise $\sigma$, the observed period is close to the deterministic version, but when the value of $\sigma$ increases, the period also decreases monotonically, as is expected. Note that after excitation, the inhibitor remains high preventing activation events. When the noise is too high no patterns are observed. We can investigate the relation between the reduction of the period and the intensity of the noise. In Figure\ref{fig:PeriodvSigma}, we plot the estimated period versus the noise intensity. We indeed see that the period decreases monotonically with the noise.

\begin{figure}
\begin{subfigure}{.49\textwidth}
  \centering
 		\def\svgwidth{\columnwidth}
    		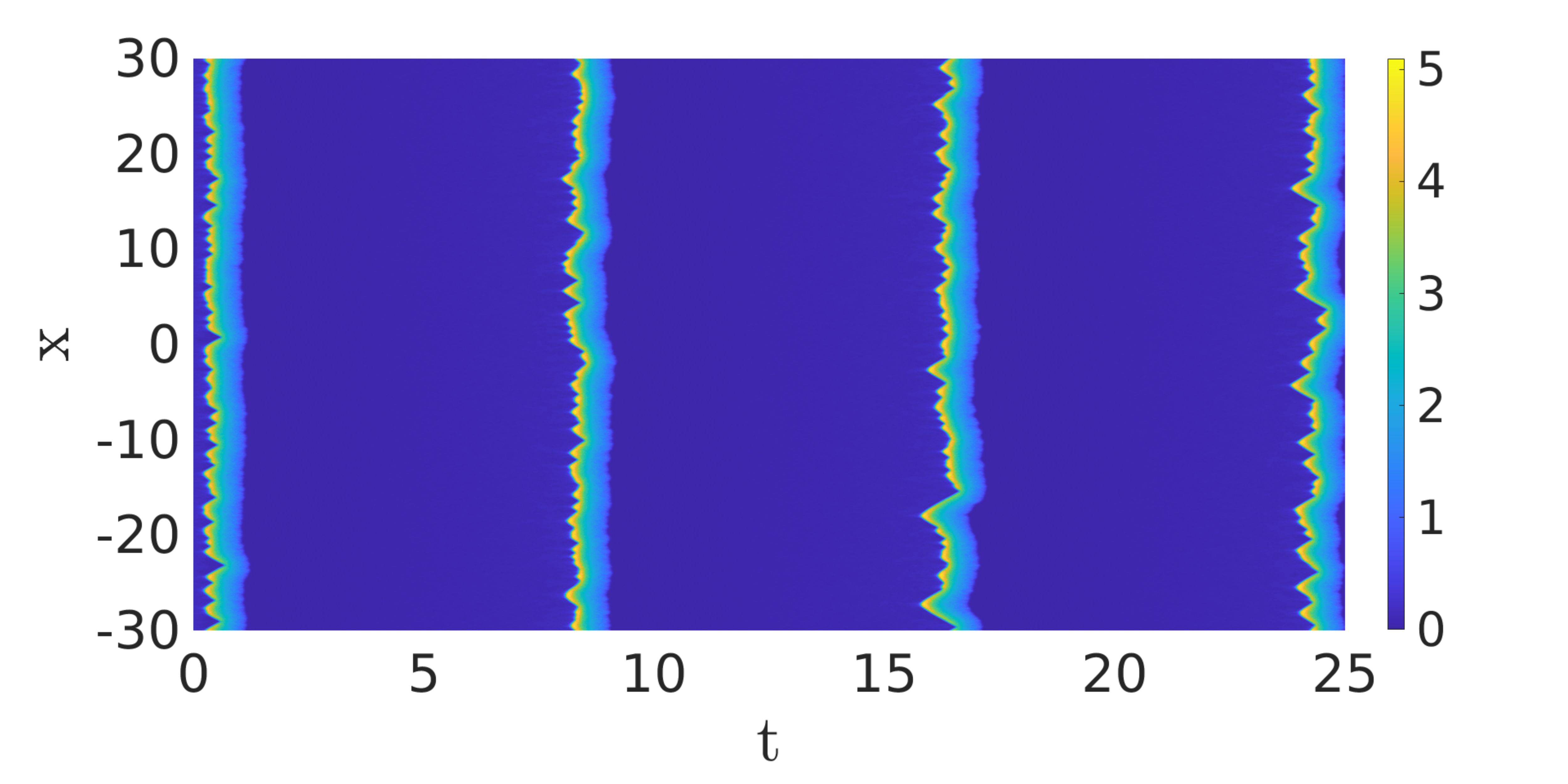
  \caption{}
    \label{fig:Sim:PeriodicNoiseU}
\end{subfigure}
\begin{subfigure}{.49\textwidth}
  \centering
 		\def\svgwidth{\columnwidth}
    		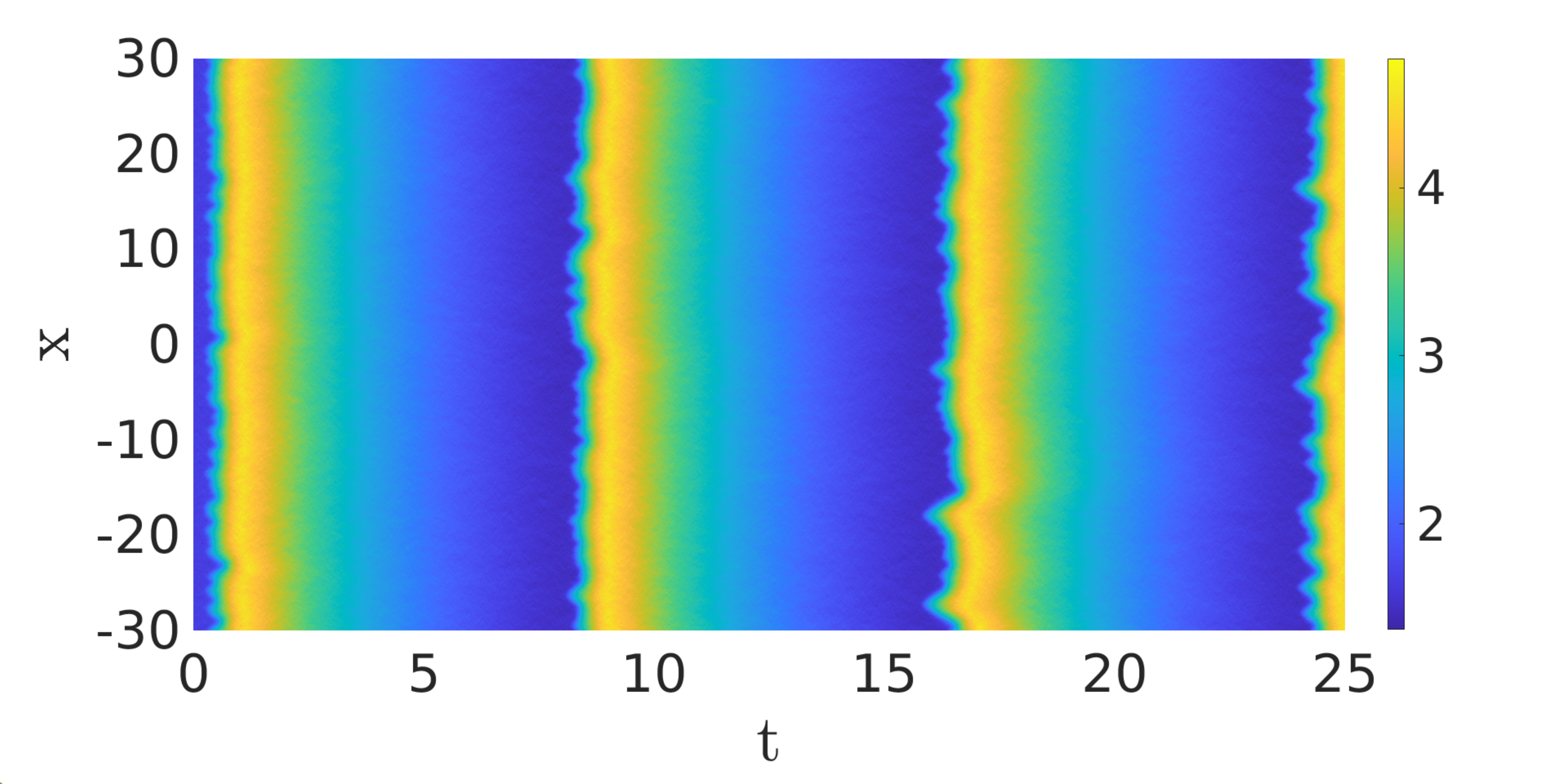
  \caption{}
    \label{fig:Sim:PeriodicNoiseV}
\end{subfigure}
\caption{Simulation of the SPDE~\sref{eq:FullSPDE}. Figure~(a) shows the activator $u$ and Figure~(b) the inhibitor $v$. When we average over the $x$-direction and measure the distance between the maxima, we find $T\approx7.87$. Same parameters as in Figure~\ref{fig:Sim:Periodic1} with $\sigma=0.01$.}
\label{fig:Sim:PeriodicNoise}
\end{figure}

\begin{figure}
\begin{subfigure}{0.49\textwidth}
    \centering
 		\def\svgwidth{\columnwidth}
    		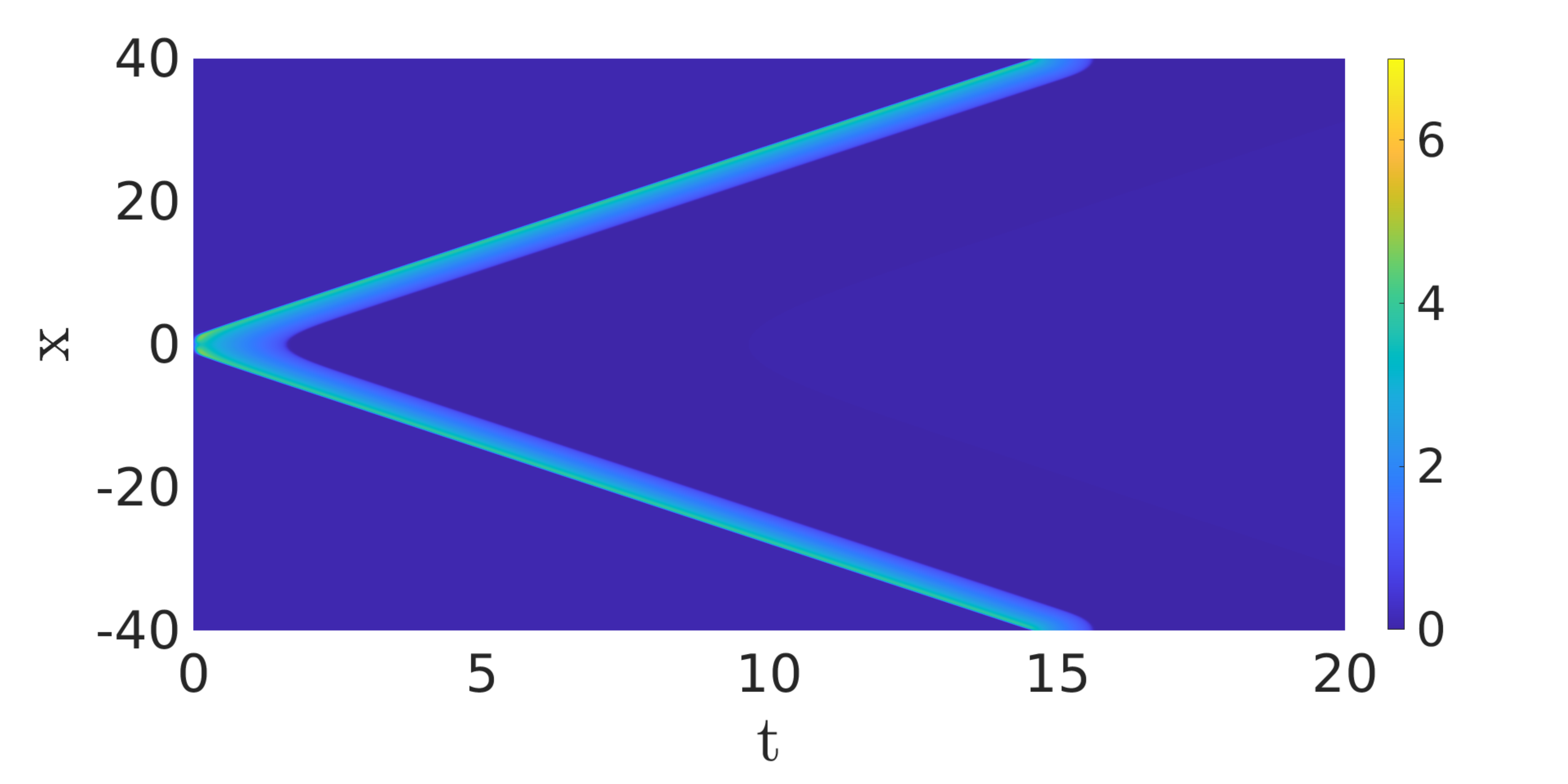
    \caption{Wild Type}
    \label{fig:WTdet}
\end{subfigure}
\begin{subfigure}{.49\textwidth}
    \centering
 		\def\svgwidth{\columnwidth}
    		\input{Figures/PTENdet.pdf_tex}
    \caption{PTEN-null}
    \label{fig:PTENdet}
\end{subfigure}
\caption{Two simulations of PDE~\sref{eq:int:MainPDE} with parameters as in~\cite{bhattacharya2020traveling}; $D_u=0.1$, $a_1=0.167$, $a_2=16.67$, $a_4=1.44$, $a_5=1.47$, $D_v=1$, $\e=0.4$, $c_1=0.1$ and, for Figure~(a),  $a_3=167$ and $c_2=2.1$, while $a_3=300.6$ and $c_2=3$ for Figure~(b). The initial condition is equal to those in the previous figures.}
\label{fig:wtvsptenDet}
\end{figure}

\subsection{Wild-Type versus PTEN-null Cells.}
\label{sec:WvsP}

Now that we have studied several different fundamental patterns, we can focus on understanding the different cell shapes. In~\cite{bhattacharya2020traveling}, two sets of parameters are compared, representing WT cells (i.e. healthy cells) and PTEN-null cells where the tumour-suppressing gene PTEN has been switched off\cite{chen2018pten}. 
First, we simulate the deterministic PDE~\sref{eq:int:MainPDE} for both sets of parameters, see Figure~\ref{fig:wtvsptenDet}. We observe that in both parameter regimes, there are two counter-propagating travelling waves but the specific profiles and speeds are different. Especially, note that the wave in Figure~\ref{fig:PTENdet} is significantly broader and higher than the wave in Figure~\ref{fig:WTdet}. 

\begin{figure}
\begin{subfigure}{0.49\textwidth}
    \centering
    \includegraphics[scale=0.27]{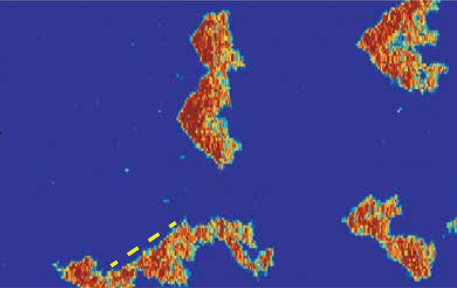}
    \label{fig:WTgil}
    \caption{}
\end{subfigure}
\begin{subfigure}{.49\textwidth}
    \centering
 		\def\svgwidth{\columnwidth}
    		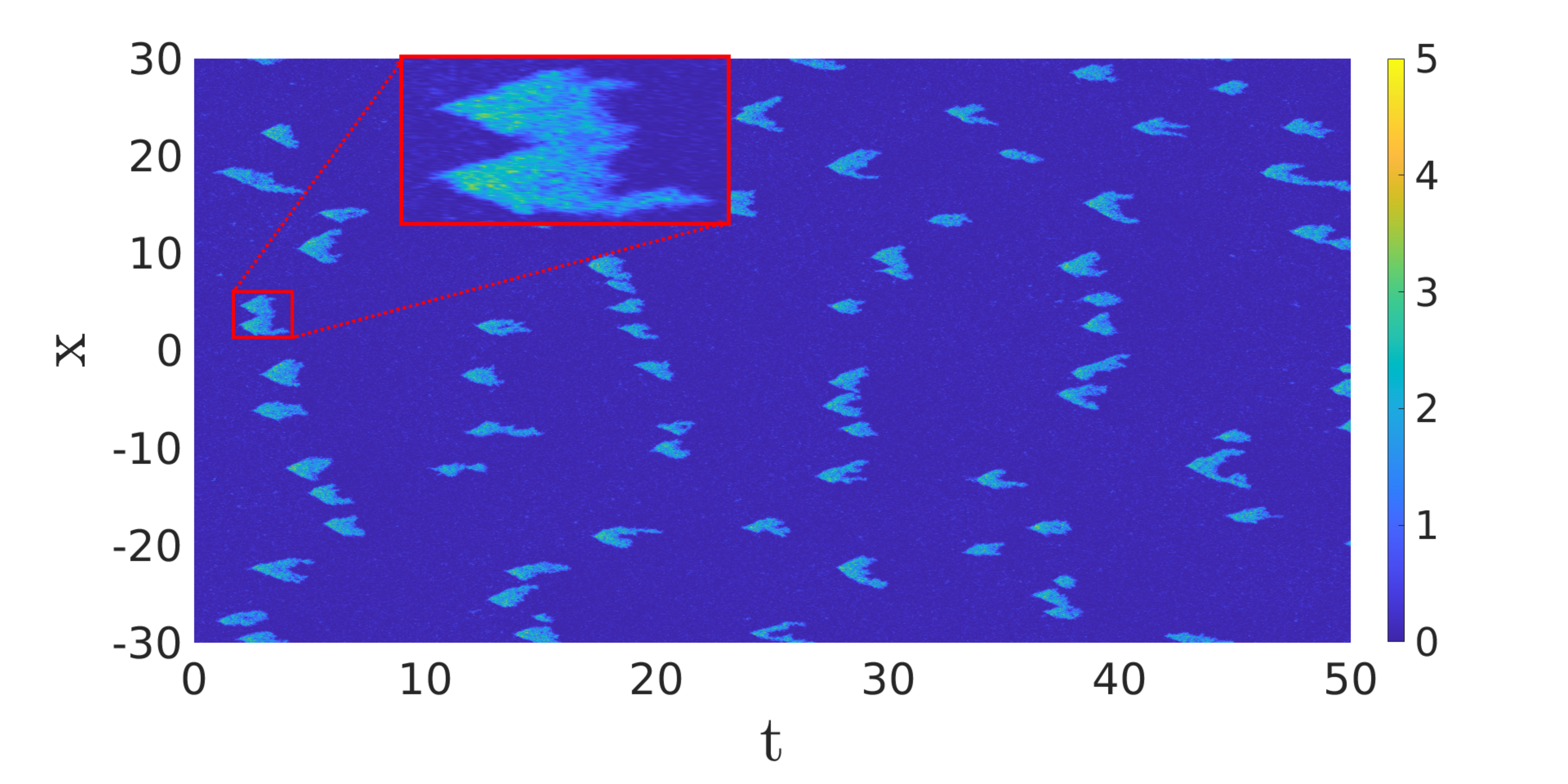
    \caption{}
    \label{fig:WTspde}
\end{subfigure}
\begin{subfigure}{0.49\textwidth}
    \centering
    \includegraphics[scale=0.27]{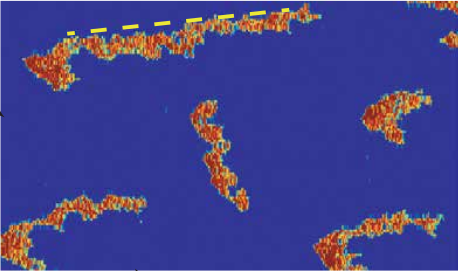}
    \label{fig:PTENgil}
    \caption{}
\end{subfigure}
\begin{subfigure}{.49\textwidth}
    \centering
 		\def\svgwidth{\columnwidth}
    		\input{Figures/PTENspde.pdf_tex}
    \caption{}
    \label{fig:PTENspde}
\end{subfigure}
    \caption{Comparison of the Gillespie model, Figures~(a) and (c) from~\cite{bhattacharya2020traveling}, versus the CLE approximation~\sref{eq:FullSPDE}, Figures~(b) and (d). The same parameters as in Figure~\ref{fig:wtvsptenDet}, with $\sigma=0.06$. The initial condition is $(u^*,2v^*)$. This can lead to an immediate excitation of the background state in Figure~(d), while in Figure~(b), the excitation of the background is more spread out. The zoom-boxes highlight the details of a single excitation.}
    \label{fig:wtvspten}
\end{figure}

When noise is applied, the statistics of the dynamics shows a clear difference. In Figure~\ref{fig:wtvspten}, we compare the SPDE simulations of~\sref{eq:FullSPDE} to the Gillespie simulations from~\cite{bhattacharya2020traveling}. Focusing on the typical shape of the excitations, there is a clear qualitative correspondence between the two types of simulation. Furthermore, in both types of simulation, the average pulse duration is longer in the case of the PTEN-null cell simulations. Note that we show the SPDE simulations on a larger spatio-temporal scale to get a better idea of the distribution in shapes and the zoom-boxes highlight the detailed structure of a typical single activation event. In the case of PTEN-null cells, the background state can be excited for much lower noise values ($\sigma\approx0.007$), while for WT cells, the noise needs to be twice as large ($\sigma\approx0.014$) as a result of the increased values of $c_2$ and $a_3$. Hence, in PTEN-null cells, an already existing pattern can more easily sustain itself, leading to the elongated shapes of Figure~\ref{fig:PTENspde}.

\section{Discussion $\&$ Outlook}
\label{sec:dis}
We set out to show how Stochastic Partial Differential Equations (CLE), or more specifically, Chemical Langevin Equations, can be used to gain more insight into the dynamics of models for cell motility. We have shown for an exemplary set of chemical reactions (see Tabel~\ref{tab:my_label}) that the CLE approach, combined with a basic analysis of the corresponding deterministic PDE, allows us to study the different possible patterns with relative ease, both qualitative and quantitative, while remaining close to the underlying chemical processes. To understand differences in cell behaviour, like the difference between wild-type and cancerous cells as in~\cite{bhattacharya2020traveling}, the study of the statistical properties of the observed dynamics is essential. For instance, an essential characteristic differentiating wild-type cells from cancerous cells is how long a pattern can survive after activation. The simulations in the previous section show that the answer not only depends on the parameters of the system but crucially on the interplay between the parameters and the noise. The CLE can be used to study this interplay. A natural question to ask is if all the stochastic terms introduced in the CLE~\sref{eq:cle2} are really necessary. Could we, for example, ignore the noise term coming from the diffusion or forget the derivation of the CLE altogether and just {\emph{naively}} add an additive white noise term to the equation for $u$? The histograms in Figure~\ref{fig:Sim:Hist} indicate that the effects of the terms that come from the diffusion are minimal (for the parameter values studied here) and therefore that these terms do not contribute meaningfully to our understanding of the cell dynamics. Note that this would solve the problem of the equation being ill-posed, see Remark~\ref{rem:ill}, and would open up the possibilities for more rigorous mathematical analysis based on the results in~\cite{hamster2020travelling}. We also noted that adding just additive white noise changes the statistics significantly, which indicates that completely abandoning the CLE approach throws away too much detail. 

In this paper, we studied a basic activator-inhibitor system with only a limited number of chemical reactions. However, the derivation of the CLE~\sref{eq:cle2} in \S\ref{sec:der} holds for any number of molecules and for any number of chemical reactions. As such, one can see this paper as a {\emph{proof of concept}} and the methodology of this paper can be directly applied to more complex regulating systems, such as the eight-component system designed in~\cite{biswas2021enhanced}. In subsequent work, we aim to work on these type of more complex model to better understand the stochastic dynamics that causes the cell to move robustly in one specific direction.  

Furthermore, as shown in detail in Appendix~\ref{sec:analysis}, the underlying deterministic RDE~\sref{eq:int:MainPDE} is amenable for rigorous mathematical analysis by using Geometric Singular Perturbation Theory~\cite[e.g.]{F79, hek2010geometric, J95, K99}. We derived a first-order approximation for the jump location where, under certain conditions, the standing wave has a sharp transition in its activator. This methodology could also be used to, for instance, further analyse the travelling waves to derive approximations for the speed of the waves. In other words, questions about the existence of localised solutions of \sref{eq:int:MainPDE} and bifurcations can thus be reduced to understanding relatively simple ODEs and the connections between them. The details of these computations are left as future work.

\bibliographystyle{klunumHJ}
\bibliography{ref}

\appendix
\renewcommand{\thefigure}{A.\arabic{figure}}
\renewcommand{\theequation}{A.\arabic{equation}}
\setcounter{figure}{0}

\section{Numerical Methods}
\label{sec:app:num}
\subsection{(S)PDE Simulations }
All the (S)PDE simulations in this paper were done using a semi-implicit Euler–Maruyama method from \cite[\S 10.5]{lord2014introduction}. For the spatial directions, a standard 2nd-order central difference is used and for the time stepping Euler-Maruyama. The deterministic linear part is evaluated at the next timestep, making it semi-implicit. To be concrete, we study an SPDE of the form
\begin{align*}
    du=[Lu+f(u)]dt+g(u)dW_t, 
\end{align*}
where $u$ is a vector, $L$ is a linear differential operator, $f, g$ are functions and $dW_t$ a white noise vector. In comparison with the main text, the vector $u$ equals $(u,v)$ and $L=D\partial_{xx}$. 

When we denote the numerical approximation of the linear part $L$ with $A$, and the spatial discretisation of $u$ at time $t$ with $\mathbf{u}(t)$, we find
\begin{align}
    \mathbf{u}(t+dt)=\mathbf{u}(t)+dt[A\mathbf{u}(t+dt)+f(\mathbf{u}(t))]+g(\mathbf{u}(t))d\mathbf{W}_t.
\end{align}
The white noise step $d\mathbf{W}_t$ is a vector where each random element is distributed as $\mathcal{N}(0,dt/h)$. Hence, the approximation for the new value $\mathbf{u}(t+dt)$
becomes
\begin{align}
       \left(I-dtA\right) \mathbf{u}(t+dt)=\mathbf{u}(t)+dtf(\mathbf{u}(t))+g(\mathbf{u}(t))d\mathbf{W}_t.
\end{align}
The equation for $\mathbf{u}(t+dt)$ is now a matrix equation and can be solved using standard solvers. It is important to realise that the algorithms from \cite{lord2014introduction} only work for Lipschitz noise terms. Hence, when the term under the square root in \sref{eq:FullSPDE} becomes close to zero, the algorithms become unstable. To correct this, we take after every timestep the maximum of $\mathbf{u}(t+dt)$ and 0. 

\begin{figure}[h]
\begin{subfigure}{.49\textwidth}
  \centering
 		\def\svgwidth{\columnwidth}
    		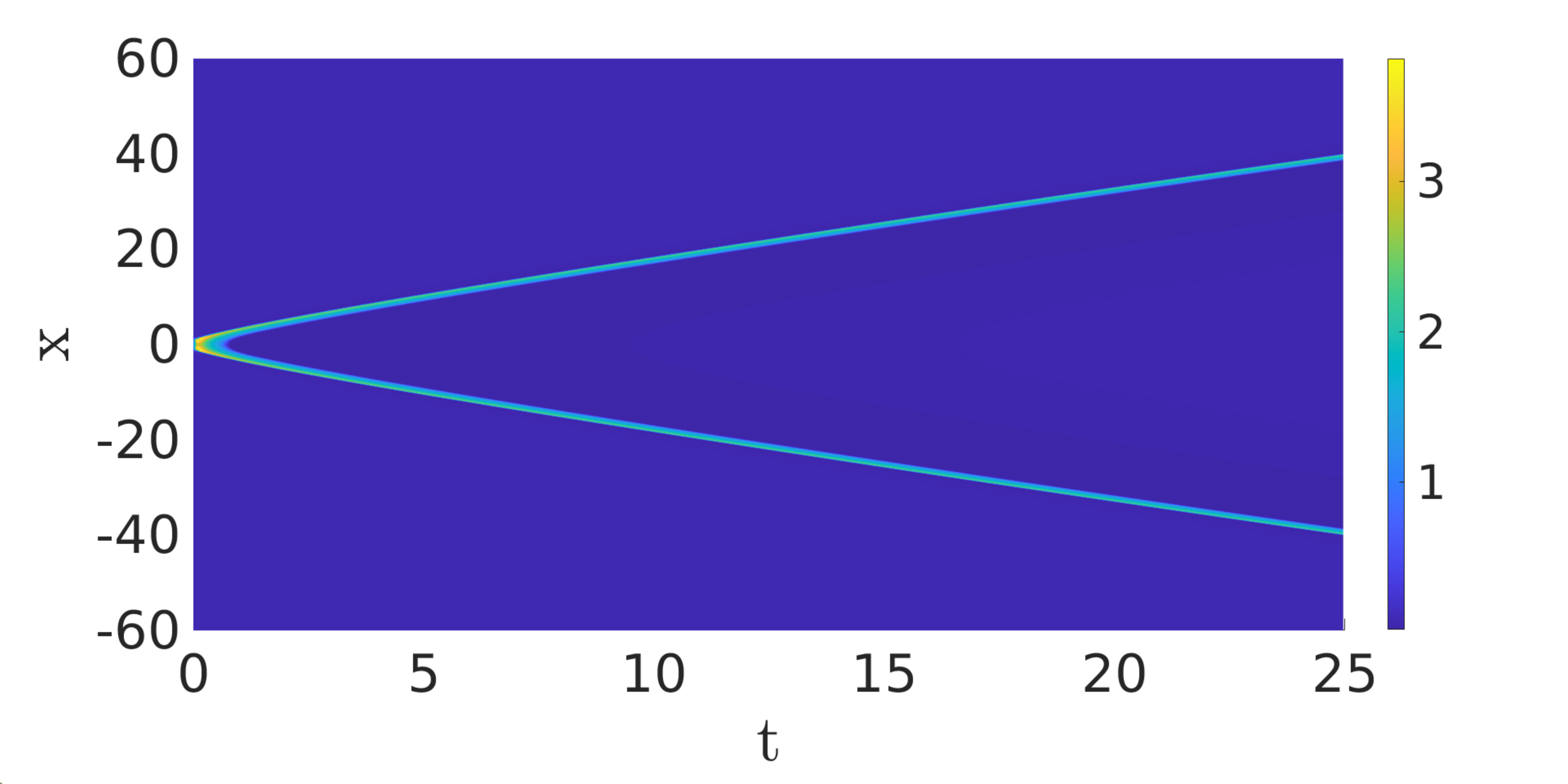
  \caption{}
\end{subfigure}
\begin{subfigure}{.49\textwidth}
  \centering
 		\def\svgwidth{\columnwidth}
    		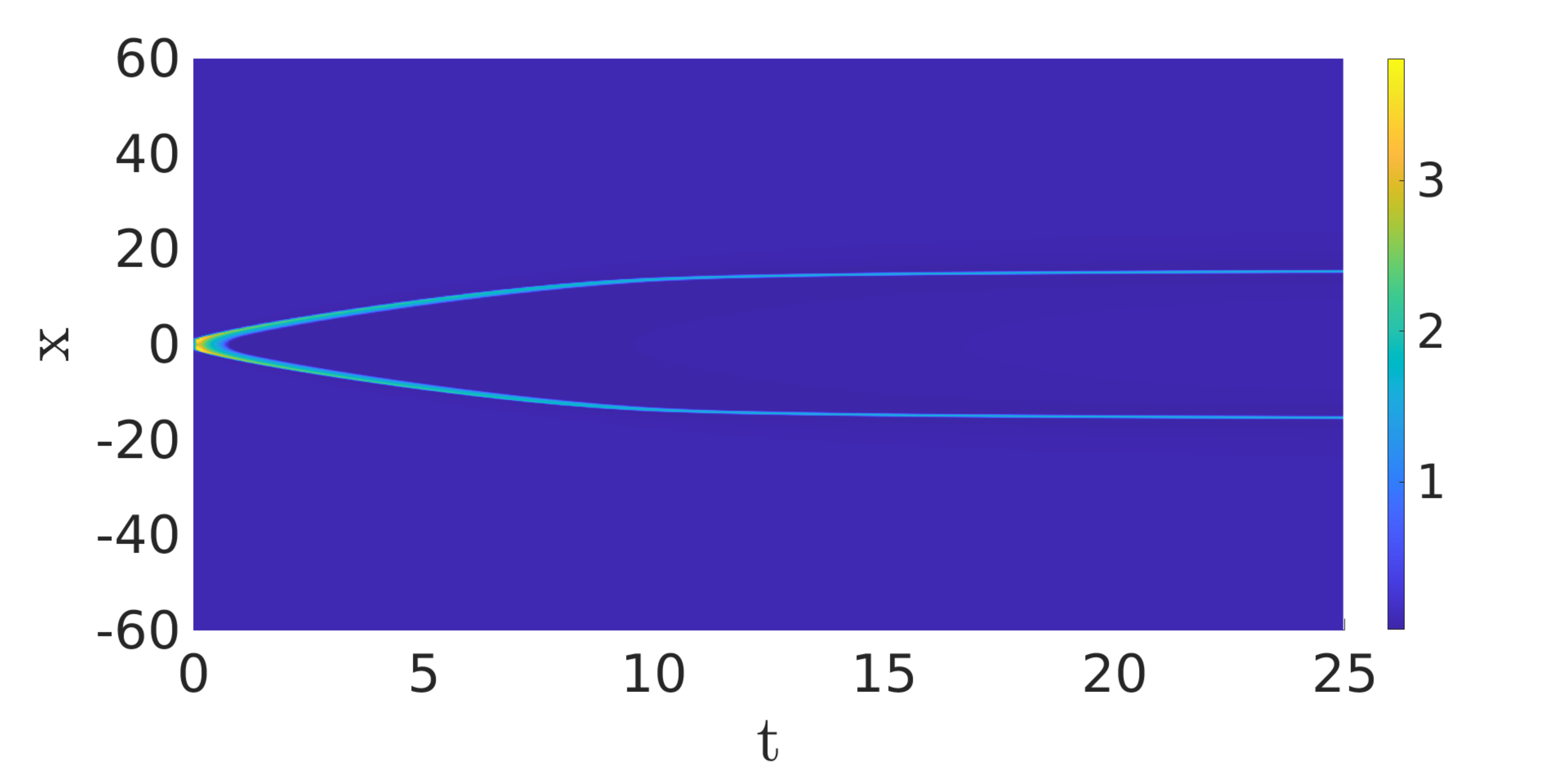
  \caption{}
\end{subfigure}
\caption{Simulation of the PDE~\sref{eq:int:MainPDE}. The spatial domain is $[-60,60]$ with $2^{12}$ gridpoints, the time interval is $[0,25]$. In Figure~(a), we used $dt=6.25\cdot 10^{-4}$ and in Figure~(b) $dt=0.0025$, i.e. 4 times larger. The parameters were set at $D_u=0.1$, $a_1=0.167$, $a_2=16.67$, $a_3=167$, $a_4=1.44$, $a_5=1.47$, $D_v=1$, $\e=0.52$, $c_1=0.15$ and $c_2=3.9$. In both cases, the initial condition is equal to the one in the main text.}
\label{fig:app:discr}
\end{figure}

The specific models studied in the main text, even the PDE \sref{eq:int:MainPDE} can be very sensitive to the size of the spatial discretisation $h$ and temporal discretisation $dt$ in certain parameter regimes. For example, when $c_1=0.15$ and the remaining parameters are equal to those in Figures~\ref{fig:Sim:Dpulse} and \ref{fig:Sim:TWext}, the dynamics can differ significantly depending on the chosen size of the discretisation, see Figure~\ref{fig:app:discr}. For the values used in the main text, $c_1=0.1$ and $c_1=0.2$, such a discrepancy was not observed for reasonable discretisations.  This is possibly related to the co-existence of travelling and standing waves in this regime. 

\begin{figure}
\begin{subfigure}{.49\textwidth}
    \centering
 		\def\svgwidth{\columnwidth}
    		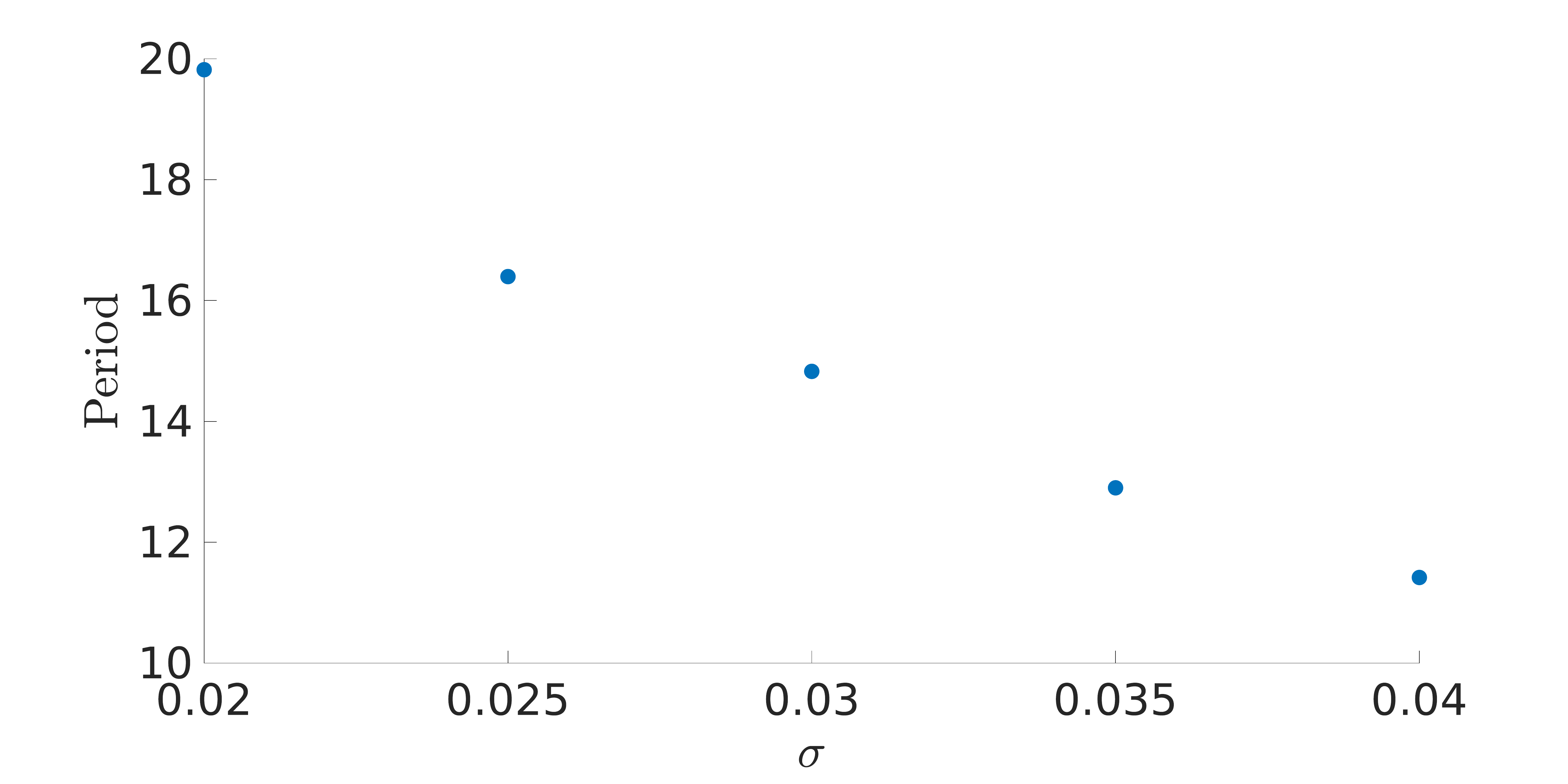
    \caption{$c_1=0.2$}
    \label{fig:QPeriodvSigma}
\end{subfigure}
\begin{subfigure}{.49\textwidth}
    \centering
 		\def\svgwidth{\columnwidth}
    		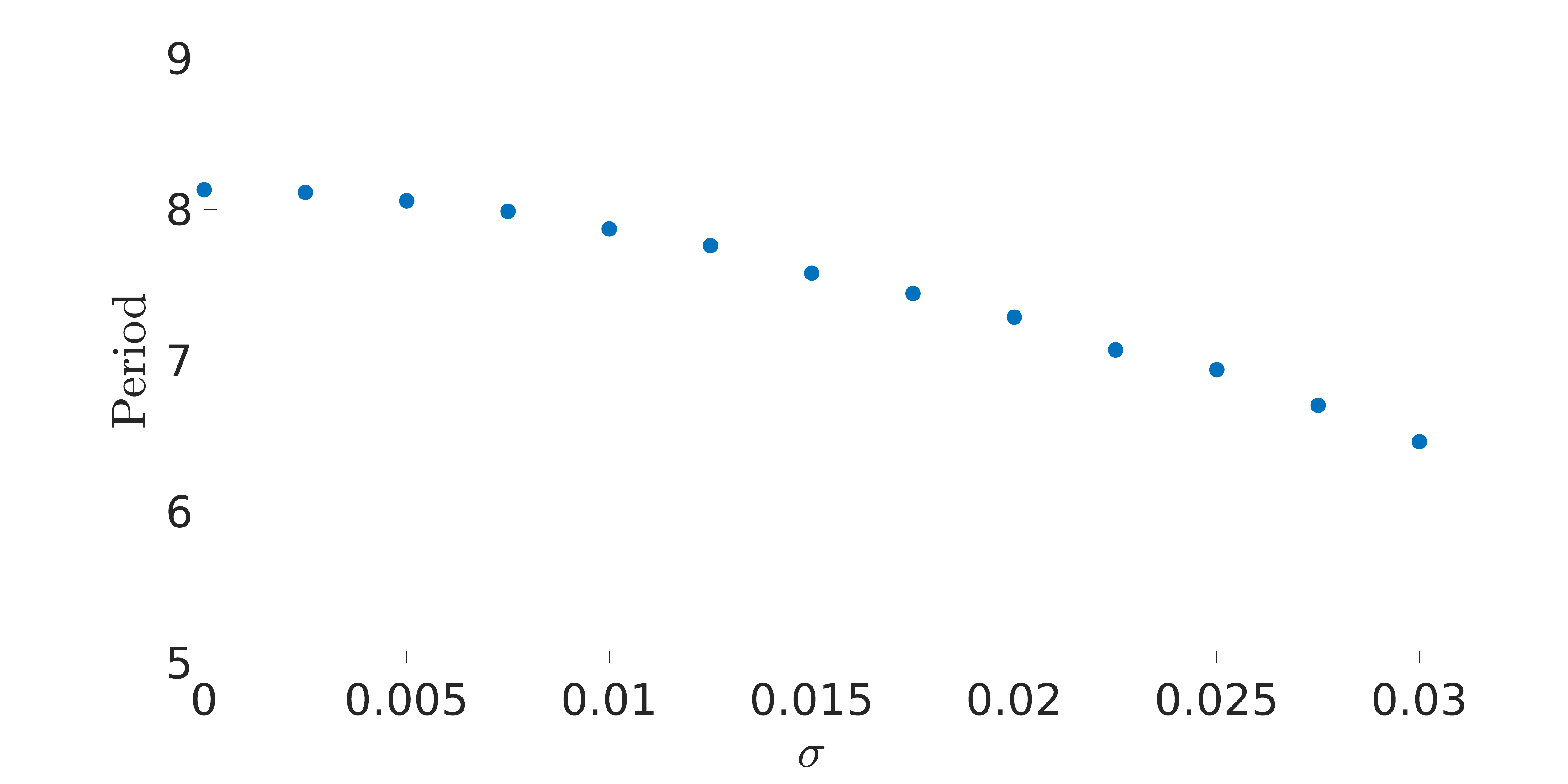
    \caption{$c_1=0.4$}
    \label{fig:PeriodvSigma}
\end{subfigure}
\caption{This figure shows the period of the dynamics of SPDE~\sref{eq:FullSPDE} in two different regimes. In Figure~(a), with $c_1=0.2$, we are in the regime of travelling waves with quasi-periodic movement as shown in Figure~\ref{fig:Sim:TWnoise}, while in Figure~(b), with $c_1=0.4$, we are in the regime of oscillations in time as shown in Figure~\ref{fig:Sim:PeriodicNoise}. In both figures, the period is estimated by computing the average in the spatial direction and subsequently computing the distances between the maxima in time.}
\end{figure}

\subsection{Pattern Recognition}
For Figures~\ref{fig:Sim:Box} and \ref{fig:Sim:Hist}, Matlab's \texttt{regionprops} algorithm is used to identify the activation events automatically. This proceeds in the following steps. First, we smooth the data using Matlab's Gaussian filter. Without smoothing, the algorithm detects multiple objects in a single event. Next, we transform the data to a binary value by comparing it with a certain threshold: we say that $u$ is \emph{activated} when it is five times its stationary value $u^*$. Then, the \texttt{regionprops} algorithm is applied with the option \texttt{BoundingBox}. 

One needs to take care of which initial condition to use. When we start with a spatial homogeneous initial condition $(u^*,v^*)$, there is a lot of activation in the first instances of the simulation, see Figure~\ref{fig:Sim:PulseNoise05}, and it is not possible to define and detect individual activation events. Therefore, we start not on the fixed point $(u^*,v^*)$, but on $(u^*,4v^*)$ plus a small perturbation. The result is that activation events only appear when $v$ has decayed enough for excitations to happen. As the decay is stochastic, and therefore not spatially homogeneous, the activation events start to appear more spread out, making it possible to determine individual events. 

\renewcommand{\thefigure}{B.\arabic{figure}}
\renewcommand{\theequation}{B.\arabic{equation}}
\setcounter{figure}{0}
\section{Analysis}
\label{sec:analysis}
The deterministic PDE~\sref{eq:int:MainPDE} has two components and ten parameters, making it difficult to  directly analyse mathematically, even for the simplest of localised structures simulated in the main text. However, these simulations do reveal that the profiles of the two components of the PDE evolve on a different spatial scale: the spatial changes in the {\emph{slow}} $v$-component are more gradual than these of the {\emph{fast}} $u$-component, see, for instance, Figures~\ref{fig:Sim:Pulse}, \ref{fig:Sim:TWprofile}, and \ref{fig:Sim:Periodic}. Furthermore, these simulations also revealed that a large part of the spatial dynamics centres around the lower and upper branch of the $u$-nullcline in the phase plane, with the $u$-profile making fast jumps in between. To amplify (and exploit) this scale separation, we set $D_v=1$ (as in~\cite{bhattacharya2020traveling}) and $D_u=\tilde\e^2$, where $\tilde\e$ is a small parameter that can be taken arbitrary small. Furthermore, we assume that our spatial domain is no longer periodic but instead unbounded\footnote{It is relatively straightforward to generalise the results for the unbounded domain to the periodic domain for the type of problems under consideration, see for example~\cite{PLOEG}}. This transforms the PDE model~\sref{eq:int:MainPDE} into 
\begin{align}
\label{eq:an:sing}
\begin{split}
    \partial_t u&=\tilde\e^2 \partial_{xx} u -a_1 u-a_2uv +\frac{a_3u^2}{a_4+u^2}+a_5\\
    \partial_t v&= \partial_{xx} v+\e(-c_1v+c_2u).
    \end{split}
\end{align} 
The small parameter $\tilde\e$ allows us to use Geometric Singular Perturbation Theory (GSPT) \cite[e.g.]{F79, hek2010geometric, J95, K99} to construct solutions that, to leading order in the small parameter, approximate the localised structures of the main text. 

In GSPT, the observation that the dynamics centres around the branches of $u$-nullcline is taken to the extreme and we construct solutions whose slow dynamics in the singular limit, i.e. in the limit of the small parameter $\tilde\e$ to zero, is confined to this nullcline, which we will refer to as the {\it slow or critical manifold}. In contrast, during the fast jump in $u$, the slow component will not change in this singular limit.
These assumptions simplify the computations and allow us to compute parts of the solution in the singular limit. The main theorems of GSPT \cite[e.g.]{F79, hek2010geometric, J95, K99}, sometimes called Fenichel Theorem 1-3, allow us to conclude that if the small parameter is small enough\footnote{Unfortunately, the theorems do not quantify what is meant by small enough.}, then there indeed is a {\emph{true}} solution of the PDE close to the one constructed in the singular limit.

Here, we only show the construction of the standing waves we found in \S\ref{sim:SW}. That is, we are interested in the fixed points of the PDE dynamics
\begin{align}
\label{eq:an:stat}
\begin{split}
    0&=\tilde\e^2 \partial_{xx} u -a_1 u-a_2uv +\frac{a_3u^2}{a_4+u^2}+a_5\\
    0&= \partial_{xx} v+\e(-c_1v+c_2u).
    \end{split}
\end{align} 
Upon defining $\tilde\e u'=p$ and $v'=q$, where $'$ denotes the derivative with respect to $x$, we can write this equation as a system of four ODEs:
\begin{align}
\label{eq:scaled:FullSlowODE}
    \begin{split}
        \tilde\e u'=&p\\
        \tilde\e p'=&a_1u+a_2uv-\frac{a_3u^2}{a_4+u^2}-a_5\\
        v'=&q\\
        q'=&\e(c_1v-c_2 u).
    \end{split}
\end{align}
Taking the scale separation to the extreme, i.e. setting $\tilde\e=0$, significantly simplifies the equation:
\begin{align}
\label{eq:scaled:0SlowODE}
    \begin{split}
        0=&p\\
        0=&a_1u+a_2uv-\frac{a_3u^2}{a_4+u^2}-a_5\\
        v'=&q\\
        q'=&\e(c_1v-c_2 u).
    \end{split}
\end{align}
We refer to this set of equations as the \textit{slow system}. This system should be understood in the following sense: on a large spatial scale, the dynamics of $(v,q)$ is approximated by lines 3 and 4 of the ODE above, and this approximation is valid in the region of the phase plane given by the algebraic equations in lines 1 and 2. We refer to the solution of these algebraic equations as the \textit{slow or critical manifold}. When we try to explicitly compute the critical manifold as a function $u(v)$, we encounter a third-order polynomial, which can be solved exactly. However, this is not practical as the graph $u(v)$ cannot be represented by a single function, but it has three branches, the upper, middle and lower branch. We will denote the upper branch with $u_+(v)$ and the lower branch with $u_-(v)$. Hence, system~\sref{eq:scaled:0SlowODE} now becomes
\begin{align}
\label{eq:scaled:SlowReducedODE}
    \begin{split}
        v'=&q\\
        q'=&\e(c_1v-c_2 u_\pm(v)).
    \end{split}
\end{align}
We refer to this equation as the \textit{reduced slow system}.

\begin{figure}
\begin{subfigure}{0.49\textwidth}
    \centering
 		\def\svgwidth{\columnwidth}
    		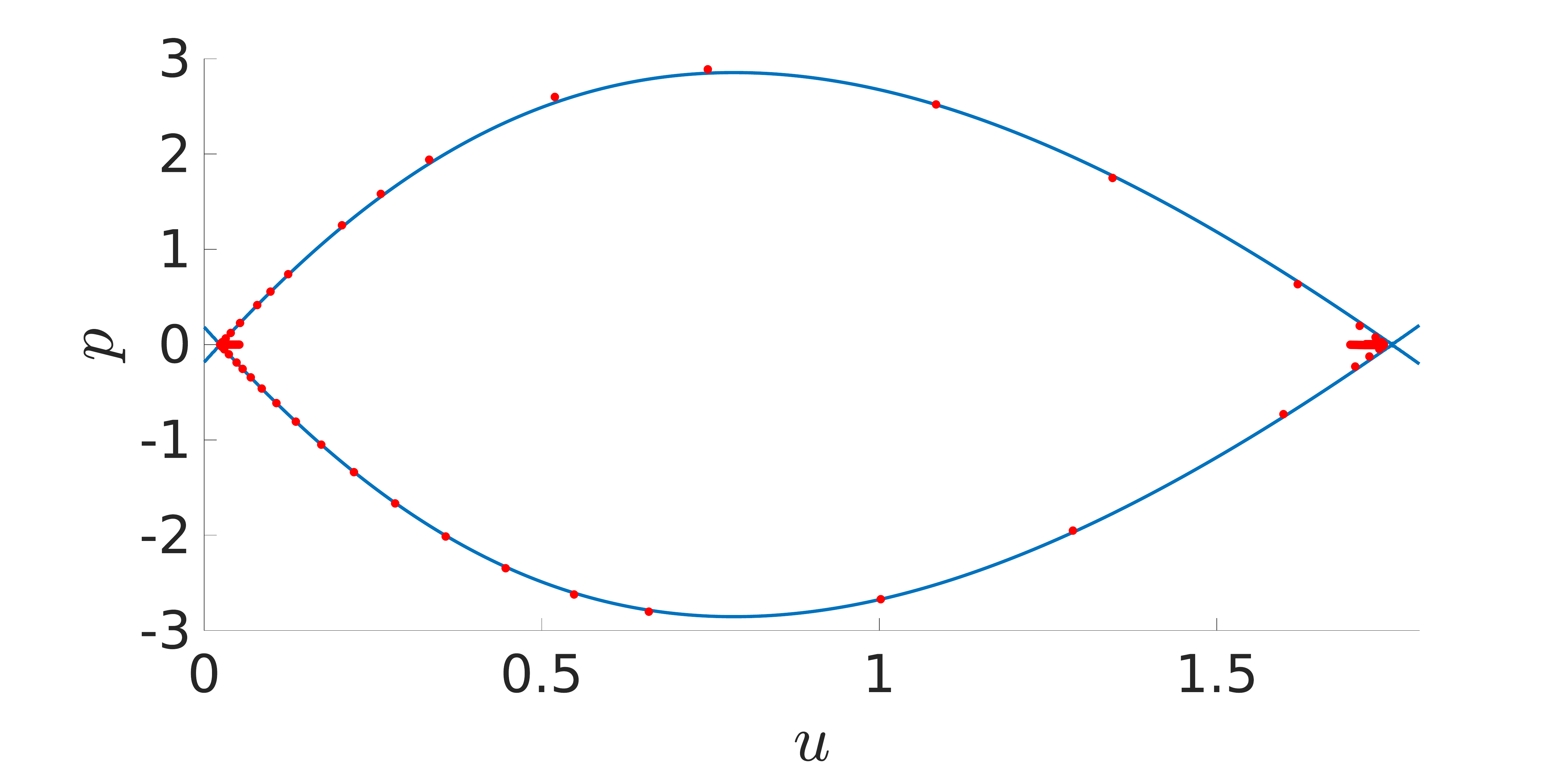
    \caption{}
    \label{fig:Hamil}
\end{subfigure}
\begin{subfigure}{0.49\textwidth}
      \centering
 		\def\svgwidth{\columnwidth}
    		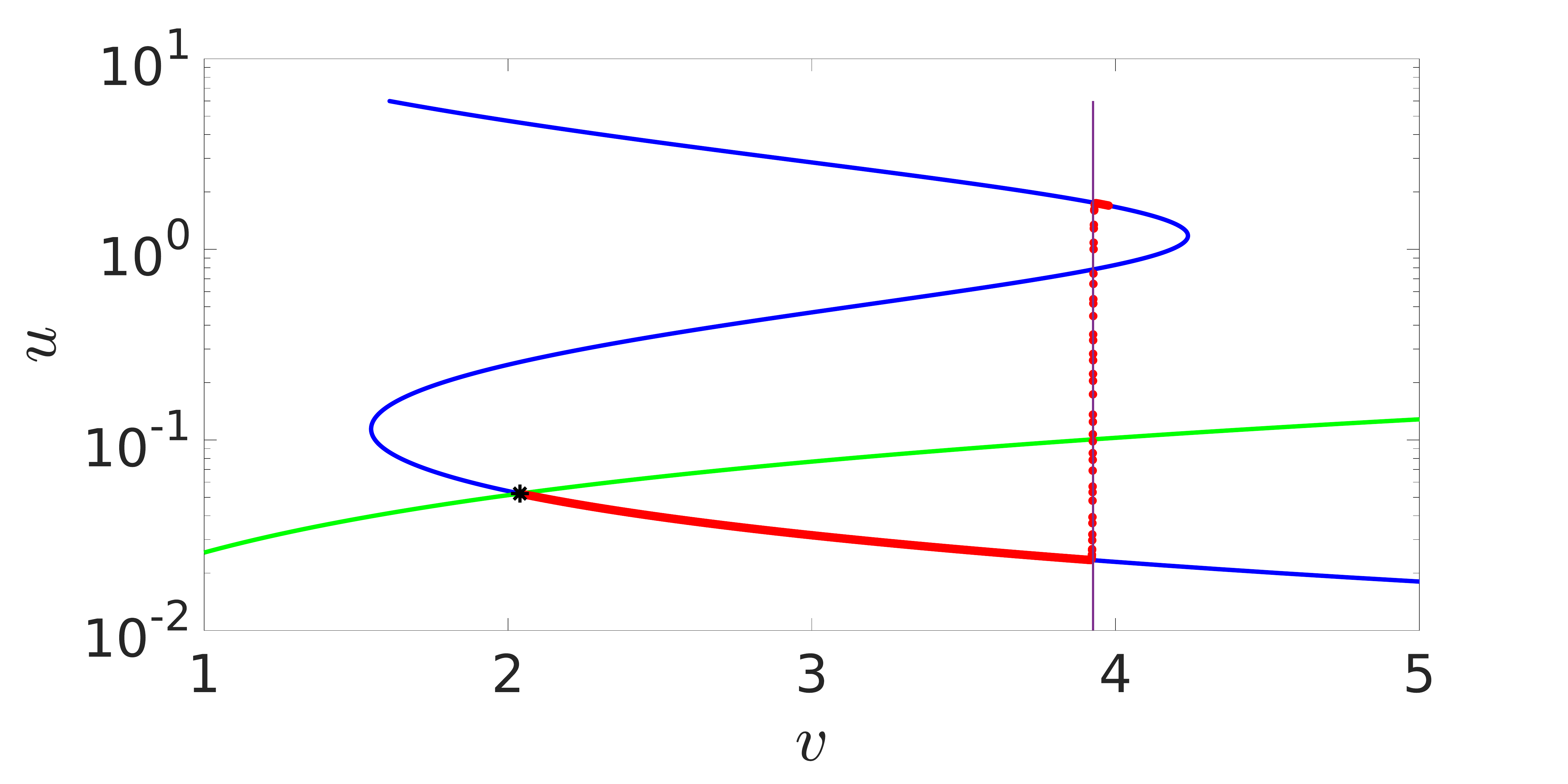
    \caption{}
    \label{fig:PulsePpSing}
\end{subfigure}
\caption{In both figures, the red dots are a solution of \sref{eq:an:stat}, found by using Matlab's \texttt{bvp4c} solver, with an initial condition coming from a PDE simulation and the small parameter $\tilde\e^2$ was set to $10^{-4}$. Note that this value corresponds to $D_u=0.01$, which is a factor ten smaller than in the main text.
In Figure~(a), we compare the fast dynamics, i.e. the jump for variables $u$ and $p$, with the predicted Hamiltonian~\sref{eq:an:Hamil}, the solid blue line. In Figure~(b), the purple line is given by $\bar v$, the value of the jump as predicted by~\sref{eq:an:barv}, while the green and blue curves denote the nullclines. The parameters are $a_1=0.167$, $a_2=16.67$, $a_3=167$,  $a_4=1.44$, $a_5=1.47$, $\e=0.52;$ $c_1=0.1$ and  $c_2=3.9$.}
\label{fig:rescaled1}
\end{figure}

When we are interested in the dynamics of $u$ instead of $v$, we must zoom in to a smaller length scale. Therefore, we define $\tilde\e\xi=x$ and use $\dot{\phantom{x}}$ to denote the derivative with respect to $\xi$. System~\sref{eq:scaled:FullSlowODE} now becomes
\begin{align}
\label{eq:scaled:FullFastODE}
    \begin{split}
        \dot u=&p\\
        \dot p=&a_1u+a_2uv-\frac{a_3u^2}{a_4+u^2}-a_5\\
        \dot v=&\tilde\e q\\
        \dot q=&\tilde\e\e(c_1v-c_2 u).
    \end{split}
\end{align}
This system is called the \textit{fast system} and is still equivalent to \sref{eq:scaled:FullSlowODE}, but when we set $\tilde\e=0$, it is no longer equivalent and reduces to
\begin{align}
\label{eq:scaled:0FastODE}
    \begin{split}
        \dot u=&p\\
        \dot p=&a_1u+a_2uv-\frac{a_3u^2}{a_4+u^2}-a_5\\
        \dot v=&0\\
        \dot q=&0.
    \end{split}
\end{align}
This shows that in the fast limit, the value of $v$ is constant. When we denote the unknown value by $\bar v$, the system reduces to
\begin{align}
\label{eq:scaled:ReducedFastODE}
    \begin{split}
        \dot u=&p\\
        \dot p=&a_1u+a_2u\bar v-\frac{a_3u^2}{a_4+u^2}-a_5.
    \end{split}
\end{align}
This system is known as the \textit{reduced fast system}. 

How can we use both reduced systems to understand the dynamics of Figure~\ref{fig:PulsePpSing}? We observe slow dynamics on the upper and lower branch of the critical manifold and a fast jump in between. The reduced fast system describes the fast jump between the upper and lower branch of the critical manifold. Therefore, a standing wave exists when this system has a heteroclinic orbit between the upper and lower branch. The reduced fast system~\sref{eq:scaled:ReducedFastODE} is a Hamiltonian system with Hamiltonian
\begin{align}
\label{eq:an:Hamil}
    H(u,p)=\frac{1}{2}p^2-\frac{1}{2}\left(a_1+a_2\bar v\right)u^2+ a_3(u-\sqrt{a_4}\arctan(u/\sqrt{a_4}))+a_5u.
\end{align}
Hence, a heteroclinic orbit exists when 
\begin{align}
\label{eq:an:barv}
H(u_-(\bar v),0)=H(u_+(\bar v),0).    
\end{align}
We cannot solve this algebraic equation exactly, but it is a straightforward numerical problem. Note that $\bar v$ only depends on the parameters $a_1,...,a_5$ and not on the parameters of the equation for $v$. For this value of $\bar v$, the Hamiltonian overlaps with the fast dynamics, as is shown in Figure~\ref{fig:Hamil}. Furthermore, from Figure~\ref{fig:PulsePpSing}, it is clear that the value for $\bar v$ is a good approximation for the location of the jump for $\tilde\e=10^{-2}$. 

Now we have all the ingredients to construct the standing wave. We start at $x=-\infty$ in the background state of the reduced slow system on the lower branch. We follow the dynamics of the reduced slow system~\sref{eq:scaled:SlowReducedODE} until we reach the value $\bar v$ where we jump to the upper branch following the reduced fast system \sref{eq:scaled:ReducedFastODE}. We will follow the slow $(v,q)$-dynamics on the upper branch until we return to the value $\bar v$, but with the opposite sign for the derivative, i.e. we trace a curve from $(\bar v,q(\bar v))$ to $(\bar v,-q(\bar v))$ in the reduced slow system. Then, we jump down again to the lower branch, which we now trace back to the background state. We implicitly assume here that the maximum value of $v$ remains below the fold of the critical manifold (which is not the case for travelling waves, see Figure~\ref{fig:Sim:TW2}). 

This example shows how GSPT can be used to construct localised solutions of \sref{eq:int:MainPDE} and also how to understand these solutions. Questions about the existence of localised solutions of \sref{eq:int:MainPDE} and bifurcations can thus be reduced to understanding relatively simple ODEs and the connections between them. 

\end{document}